\documentclass[12pt]{article}
\usepackage{amsfonts}
\usepackage{amsmath,amsthm}
\usepackage{amssymb}
\usepackage{mathrsfs}
\usepackage{authblk}
\usepackage{mathtools}
\usepackage[pdfpagelabels=true]{hyperref}
\usepackage{cleveref}
\usepackage{amscd}
\usepackage[utf8]{inputenc}
\usepackage{t1enc}
\usepackage{lmodern}
\usepackage{dsfont}
\usepackage[font=small,labelfont=sc]{caption}
\usepackage{xcolor}
\usepackage{stmaryrd} 
\usepackage{multicol}
\usepackage{tasks}
\usepackage{caption}
\usepackage{babel}
\usepackage{multirow}
\usepackage[mathscr]{eucal}
\usepackage{comment}
\usepackage{graphicx}
\usepackage[a4paper, margin=4cm]{geometry}
\usepackage{diagbox}
\usepackage{graphics}
\usepackage{tikz}
\usetikzlibrary{cd}
\usepackage{pict2e}
\usepackage{epic}
\numberwithin{equation}{section}

\DeclareMathOperator{\Id}{Id}

\DeclareMathOperator{\clarke}{\partial^* \!}
\DeclareMathOperator{\Nor}{\mathrm{Nor}}

\DeclareMathOperator{\Conv}{\mathrm{Conv}}
\DeclareMathOperator{\Cone}{\mathrm{Cone}}
\DeclareMathOperator{\Tan}{\mathrm{Tan}}

\DeclareMathOperator{\sff}{\mathrm{I\!I}}
\DeclareMathOperator{\dgm}{\mathrm{dgm}}
\DeclareMathOperator{\reach}{\mathrm{reach}}

\DeclareMathOperator{\interior}{\mathrm{int}}

\DeclareMathOperator{\tq}{\; | \;}

\newcommand{\Sphere}{\mathbb{S}}
\newcommand{\N}{\mathbb{N}}
\newcommand{\R}{\mathbb{R}}

\newcommand{\scal}[2]{ \left \langle #1 , #2 \right \rangle}

\newcommand{\module}[1]{\left\lvert #1 \right\rvert}
\newcommand{\norme}[1]{\left\lvert \left \lvert #1 \right \rvert \right \rvert}

\newcommand{\eps}{\varepsilon}

\newcommand{\complementaire}[1]{^\neg #1}

\newcommand{\normalise}[1]{\frac{#1}{\norme{#1}}}

\newtheorem{theorem}{Theorem}[section]
\newtheorem*{theorem*}{Theorem}
\newtheorem{corollary}[theorem]{Corollary}
\newtheorem{lemma}[theorem]{Lemma} 
\newtheorem{definition}[theorem]{Definition}
\newtheorem{proposition}[theorem]{Proposition}
\newtheorem{remark}[theorem]{Remark}

\title{Generalized Morse theory for tubular neighborhoods}

\date{\today}
\begin{document} 
\emergencystretch 3em

\author[1,2]{Antoine Commaret\thanks{antoine.commaret@inria.fr}}
\affil[1]{\textit{Centre INRIA d'Université Côte d'Azur, Valbonne, France}}
\affil[2]{\textit{Laboratoire Jean-Alexandre Dieudonné, Nice, France}}

\maketitle

\begin{abstract}
We define a notion of Morse function and establish Morse theory-like theorems on offsets of any compact set in a Euclidean space at regular values of their distance function. Using non-smooth analysis and tools from geometric measure theory,
we prove that the homotopy type of sublevel sets of these Morse functions changes at a critical value by gluing exactly one cell around each critical point.\end{abstract}



\section{Introduction} \label{section:intro}

In his celebrated book \textit{Morse Theory} \cite{milnor_morse}, Milnor describes the changes in topology of the closed sublevel sets $X_c \coloneqq f^{-1}(-\infty, c]$ when $c$ 
runs among  $\R$ for certain $C^2$ functions over a compact $C^2$ manifold $f: X \to \R$  satisfying generic conditions,
which he calls  \emph{Morse functions}. In this setting, Milnor shows that topological changes only happen at a finite number of reals called \emph{critical values}. These are the values the function $f$ takes at \emph{critical points}, which are the points where the differential of $f$ vanishes.
Furthermore, around a critical point $x$ with critical value $c = f(x)$, 
when $\eps$ is small enough, the 
homotopy type of the sublevel sets $X_{c+\eps}$ is obtained from $X_{c-\eps}$ by gluing a cell (i.e., a set homeomorphic to the closed unit ball of a Euclidean space) around $x$.

More precisely, a $C^2$ function $f :X \to \R$ is said to be Morse when its Hessian is non-degenerate at every critical point. In this case the previous considerations can be summarized by the two fundamental results of Morse theory, which we call \textit{Morse theorems}:  
\begin{itemize}
\item \textbf{Constant homotopy type.} If an interval $[a,b]$ does not contain any critical value of $f$, $X_a$ has the same homotopy type as $X_b$. 
\item \textbf{Handle attachment.} Around a critical value $c$ of $f$, when $\eps$ is small enough the homotopy type of $X_{c+\eps}$ is obtained from $X_{c-\eps}$ by gluing a cell around each critical point $x_i \in f^{-1}(c)$, where the dimension of the cell is the index of the Hessian of $f$ at $x_i$.
\end{itemize}

In the classical setting the handle attachment lemma is obtained by studying the second order approximation of $f$ around a critical point, taking the curvature of $X$ into account. Such considerations do not apply when $X$ is not a $C^2$ manifold.  
Several works aimed at adapting the constant homotopy and handle-attachment lemmas to other classes of sets or functions.  
A systematic analysis of Morse theory for continuous functions in metric spaces was given by Corvellec \cite{MorseCorvellec,DeformCorvellec} using the notion of \textit{weak slope} defined by Degiovanni \& Marzocchi \cite{WeakSlope} at a point $x$ as local maximal rate of decrease of a real valued map $f$ deformed by a homotopy; in this context a critical point is defined as a point of weak slope zero. In between critical values, the above constant homotopy type theorem is verified \cite{MorseCorvellec}, while the homotopy type changes at critical values by a gluing process, but not necessarily that of a cell - a weaker handle attachment lemma. 

During the eighties, 
significant
works extended Morse theory to $C^2$ functions restricted to a stratified subset of a Riemannian manifold, culminating in the monograph of Goresky and MacPherson \cite{stratified}. In the context of stratified Morse theory, the change of topology around a non-degenerate critical point is obtained by the gluing of the so-called local Morse data, which 
is a stratified set 
that might not be a cell. 
With a similar version of the handle-attachment lemma, Morse theory was extended to broader classes of functions such as so-called min-type functions on a manifold \cite{min_type}, which are functions that can be locally written as the minimum of a finite number of $C^2$ functions, and for 
functions of the form $(d_Z)_{|Y}$ where $Z, Y$ are two real algebraic subsets of $\R^d$ and $d_Z$ is the function distance to $Z$
\cite{MorseDistance}. An adaptation of the Morse lemmas to the combinatorial settings of functions on CW-Complexes can be found in the discrete Morse theory of Forman \cite{DiscreteMorse}. 

\begin{figure}[h!]
    \centering
    \includegraphics[width = 0.9\textwidth]{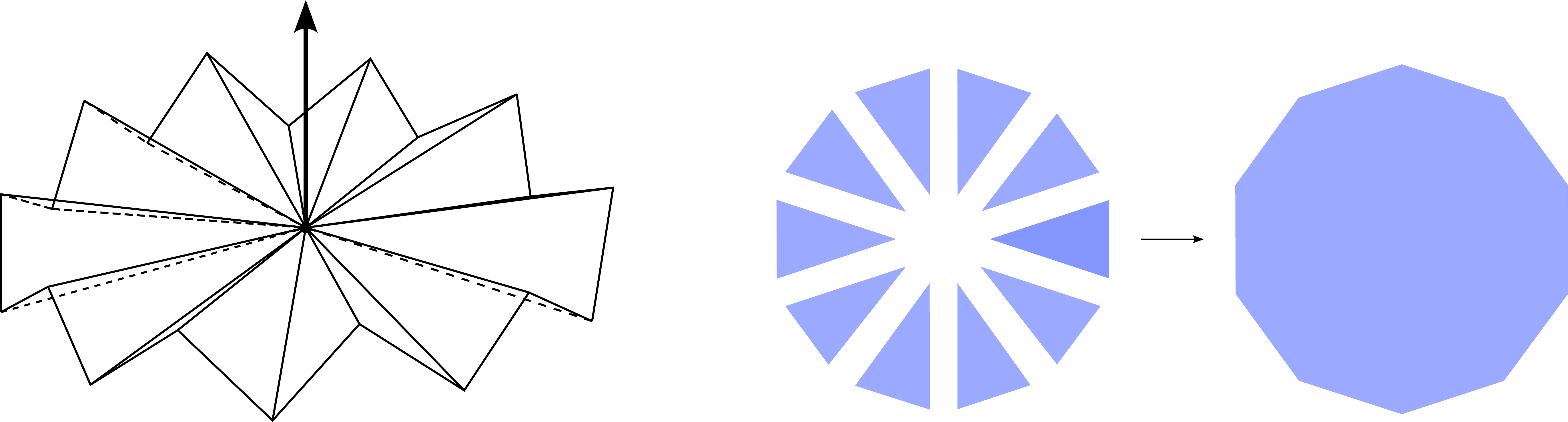}
    \caption{Consider the height filtration (indicated by the arrow) on a ruff as in the left figure. The topological event around the central point is not equivalent to the gluing of a cell as the number of connected components changes by more than one when the filtration encounters the critical point.}
    \label{fig:my_label}
\end{figure}

 In 1989, Fu \cite{fuCurvatureMeasuresGeneralized1989} proved Morse theorems (with the classical handle-attachment lemma) for smooth functions on any compact set $X$ with a $C^{1,1}$-hypersurface boundary and more generally to sets with positive reach in Euclidean spaces.
The reach of a subset of $\R^d$ was defined in 1959 by Federer \cite{curvature:federer} as the largest distance to this set under which any point has a unique closest point. The class of sets with positive reach contains submanifolds and convex subsets of $\R^d$, and intuitively consists of subsets of Euclidean spaces whose singularities are locally convex. A set with positive reach admits a unit normal bundle, which is a $(d-1)$-manifold generalizing the normal bundle of submanifolds of Euclidean spaces. From this object one can define curvatures in a generalized sense for sets of positive reach \cite{curvature:federer}, as well as a new notion for critical points used by Fu in \cite{fuCurvatureMeasuresGeneralized1989} to prove the Morse theorems. His reasoning is the main inspiration for the present article, as we adapt his proofs to our setting using non-smooth analysis. 

Our contribution is as follows. We prove that the Morse theorems (with the classical handle-attachment lemma) extend to smooth functions on any offset of a subset of $\R^d$ at a regular value of its distance function,
 that is, sets of the form $Y^{\eps}  = d^{-1}_Y(-\infty, \eps]$ where $\eps > 0$ is a regular value of the distance function $d_Y : x \mapsto \min_{y \in Y} \norme{x-y}$ in the sense of Lipschitz analysis.
Such sets are Lipschitz domains which are not necessarily smooth nor stratified, and we call them \emph{complementary regular sets}. 
Although the complement set of a complementary regular set has positive reach, in general one cannot infer Morse-like results on maps restricted to a subset $X$ of $\R^d$ from that of maps restricted to its complement set
\footnote{It is possible that the homotopy type of the sublevel-set filtration of $f_{|X}$ changes by the gluing of a cell around a critical point $x$, but not that of $(-f)_{|\overline{\R^d \setminus X}}$.  See \Cref{rk:complement} for an example where $X$ is a stratified subset of $\R^2$ and $f$ a linear form.}, 
which motivated the present work. In the process we obtain several results relating the distance function of a complementary regular set to its tangent cones.

In contrast to Morse theory for continuous functionals \cite{MorseCorvellec} and
 stratified Morse theory \cite{stratified}, which to the best of our knowledge were the only previous extensions of Morse theory to a class containing sets with possibly non-convex singularities, we prove that the classical handle attachment lemma does hold for complementary regular sets:
 topological changes in the sublevel set filtration of a Morse function are given by the attachment of a cell around each critical point - with appropriate definitions of critical points and Hessians (see \Cref{sub:critical}). While it is hopeless to extend the Morse theorems to general compact sets, which can be topologically wayward, this result shows that Morse theorems do typically hold when one replaces the compact set with an arbitrarily small tubular neighborhood.

\begin{theorem}[Main Result]
\label{th:boss}
Let $Y$ be a compact subset of $\R^d$ and $\eps > 0$ be a regular value of the distance function to $Y$. Let $X = Y^{\eps}$ be the $\eps$-offset of $Y$.  Let $f: \R^d \to \R$ be a smooth function such that $f_{|X}$ is Morse, i.e., it admits only non-degenerate critical points.

Then the set of critical points of $f_{|X}$ is finite, and for every regular value $c$ of $f_{|X}$, $X_c \coloneqq X \cap f^{-1}(-\infty,c]$ has the homotopy type of a CW-complex with one cell per critical point with value less than $c$, whose dimension is determined by the index of the Hessian of $f_{|X}$ and the curvatures of $X$. 
\end{theorem}
The precise statement of this result is given in \Cref{cor:MR}, where in particular the exact definition of the index of the glued cells is given.

 Along the way, we prove an independent result on approximate inverse flows of any real valued, locally Lipschitz map $\phi$ whose Clarke gradient is uniformly bounded away from zero on $\phi^{-1}(a,b]$. Under these assumptions, we build a deformation retraction of $\phi^{-1}(-\infty,b]$ onto $\phi^{-1}(-\infty,a]$ whose trajectories are rectifiable 
and have uniformly bounded length,
see \Cref{prop:flot}.
 Although our result was
most likely known amongst communities making use of Clarke gradients, we
were not able to find a reference for it and chose to provide a full proof. 
  In the more general context of continuous functionals on a metric space, a deformation retraction between sublevel sets had already been built by Corvellec in Theorem 2 of \cite{DeformCorvellec} under the more general assumption that the weak slope is uniformly bounded away from zero. However, this construction does not provide the rectifiability or bounded-length properties of the trajectories.

\section*{Outline}

In \Cref{section:def}, we define the objects used throughout this article.
\begin{itemize}
    \item \Cref{sub:preli} contains definitions and illustrations of the classical tools of our study. This includes the $\reach_{\mu}$ and the $\reach$ of a compact subset of $\R^d$, eroded sets $X^{-r}$ for any positive real $r$ and any $X \subset \R^d$, Clarke gradients of locally Lipschitz functions, normal and tangent cones of a compact set with positive reach.
    \item In \Cref{sub:normal} we recall the definition of the unit normal bundle of sets with positive $\reach$ and define the normal bundle of their complement set by adapting classical definitions. We describe how local curvatures of such sets are related to their normal bundle.
    \item \Cref{sub:critical} gives the definitions and notations of critical points and Hessian for a restricted function $f_{|X}$ for sets with positive reach used by Fu in \cite{fuCurvatureMeasuresGeneralized1989}. We will use the same definitions of critical points, Hessians and non-degeneracy for the class of \textit{complementary regular sets} defined in \Cref{sub:complementary_regular}.   
    \item \Cref{sub:clarke} focuses on properties of locally Lipschitz functions. More precisely, we build a retraction between sublevel sets $\phi^{-1}(-\infty, a]$, $\phi^{-1}(-\infty, b]$ when $a < b$ of any locally Lipschitz map $\phi : U \to \R$, where $U$ is an open subset of $\R^d$,  assuming a bound from below on the distance to zero of its Clarke gradient on $\phi^{-1}(a,b]$. Our result is obtained by interpolating Clarke gradients, yielding so-called approximate flows whose trajectories are Lipschitz and bounded in length. This result is independent of Morse theory.
    \item In \Cref{sub:normal_clarke} we establish a link between the normal bundle of a set $X$ and the Clarke gradient of its distance function $d_X$. This crucial step allows us to use results from non-smooth analysis on assumptions about critical points of $f_{|X}$.  
\end{itemize}

\Cref{section:morse_complement} articulates the previous results to establish the main theorem.
\begin{itemize}
    \item In \Cref{sub:complementary_regular} we define the class of \textit{complementary regular sets}, which are the sets verifying the assumptions needed in our reasoning through the remainder of the section to prove Morse theory results. We prove that $X \subset \R^d$
 is a complementary regular set if and only if it is an offset of some compact set $Y$ at a positive regular value of $d_Y$.
    \item In \Cref{sub:smoothed} we describe how to build functions $f_{r,c}$ such that the sublevel sets $X^{-r}_c = X^{-r} \cap f_{r,c}^{-1}(-\infty, c]$ and $X_c$ have the same homotopy type when $c$ is a regular value and $r > 0$ is small enough.
    To that end we consider some locally Lipschitz functions and prove that they verify the assumptions needed in the theorems of \Cref{sub:clarke}. The retractions obtained are used to build a homotopy equivalence between $X^{-r}_c$ and $X_c$.
    \item In \Cref{sub:isotopie} we show that between critical values, 
 the homotopy type of sublevel sets
does not change. This is done by applying once again \Cref{sub:clarke} using some computations from the previous section.
    \item In \Cref{sub:handle} we describe the topological changes happening around a critical value as long as it has only one corresponding critical point which is non-degenerate. We adapt the proof from Fu \cite{fuCurvatureMeasuresGeneralized1989} to our setting, circumventing the problem of considering sets with reach zero using non-smooth analysis. We then extend this result to critical values with a finite number of corresponding critical points which are all non-degenerate.
    \item In \Cref{sub:persistence}, we describe the connections of our work with persistent homology. We discuss the definitions of persistent homology transform and the Euler characteristic transform for complementary regular sets, and prove that the last one is an injective operation on complementary regular sets.
\end{itemize}

\section{Definitions and useful lemmas} \label{section:def}

\subsection{Preliminaries}
\label{sub:preli}

We fix $d \in \N$ to be the dimension of the Euclidean space in which our objects will be included. The word \textit{smooth} is meant as a synonym for $C^2$. The canonical scalar product over $\R^d$ will be denoted by $\scal{ \cdot }{ \cdot }$, and by $B(x,r)$ we will denote the closed ball of radius $r$ centered at $x \in \R^d$. The inclusion (proper or not) of a set into another will be denoted by $\subset$. $\mathcal{H}^k$ denotes the $k$-Hausdorff measure.
\vspace{1em}

 For any subset $X$ of $\R^d$, $\interior(X)$ denotes the interior of $X$ while $\overline{X}$ denotes its closure, both for the topology of $\R^d$ induced by the Euclidean distance. Throughout this paper, we define the \emph{complement set} of $X$ as the closure of the classical complement set and denote it by $\complementaire{X} \coloneqq \overline{\R^d \setminus X} = \R^d \setminus \interior(X)$. 
Taking the complement set is not an involutive operation, as we more precisely have $\complementaire{(\complementaire{X})} = \overline{\interior{(X)}}$, an identity used several times in the paper.
 The topological boundary $\overline{X} \setminus \interior(X) = \complementaire{\! X} \cap \overline{X}$ 
 will be denoted by $\partial X$. In particular if $X = \overline{\interior{(X)}}$, we have $\partial X = \partial (\complementaire{X})$, and we say that $X$ is a domain.
 \vspace{1em}

\begin{figure}[h!]
    \centering
    \captionsetup{justification = centering}
    \includegraphics[width = 0.92 \textwidth]{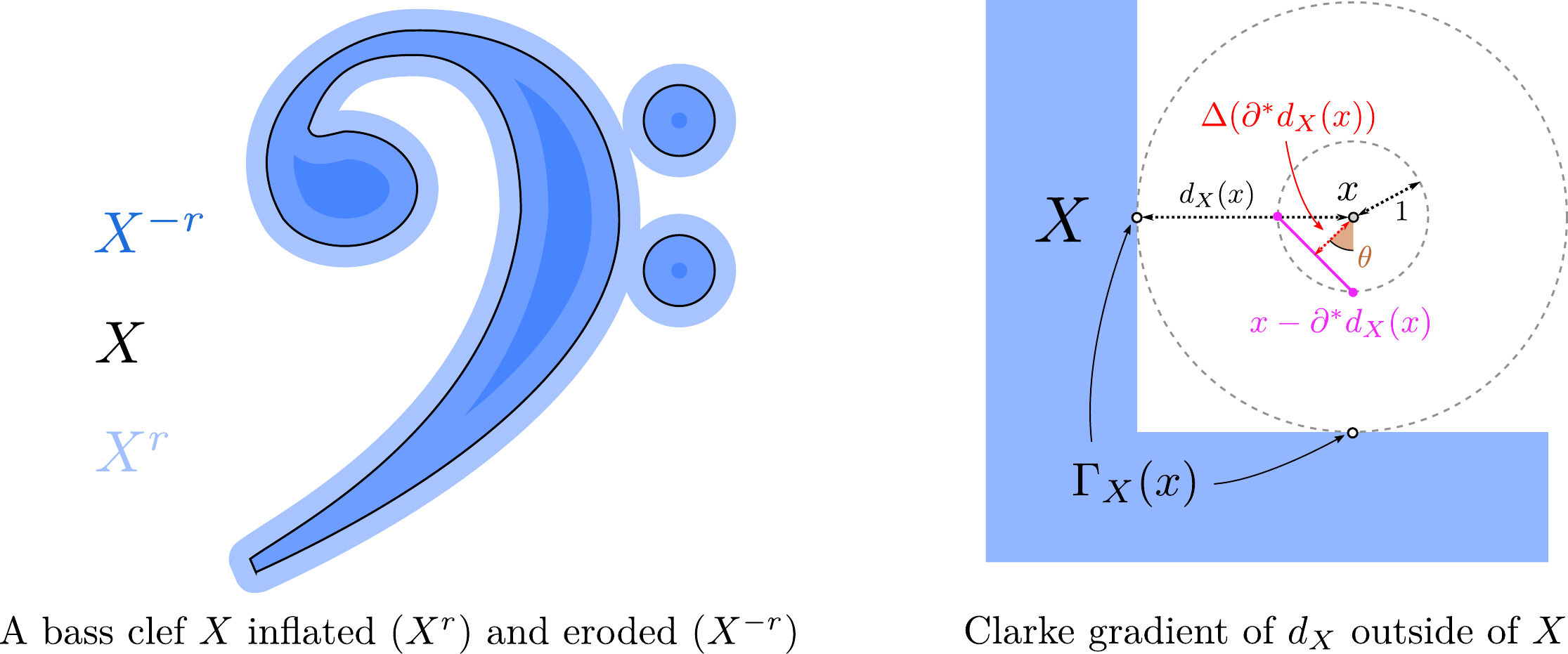}
    \caption{Offsets of $X$ and Clarke gradient of $d_X$ outside of $X$.}
    \label{fig:distance}
\end{figure}

Let $A$ be a subset of $\R^d$. Its \emph{distance function} is $d_A: x \mapsto \inf \{\norme{x - a} \tq a \in A \}$. Any distance function is 1-Lipschitz over $\R^d$. For any positive real $r$ and any subset $X$ of $\R^d$, define the $r$ and $-r$ tubular neighborhoods (respectively \emph{offsets} and \emph{counter-offsets}, also called \textit{eroded sets}) of $X$ (see (\Cref{fig:distance}), left) as follows: 
    \[ 
    \begin{array}{r l}
         X^{r} \coloneqq &  \left \{ x \in \R^d \tq  d_X(x) \leq r \, \right \}  \\
         X^{-r} \coloneqq &  \left \{ x \in \R^d \tq d_{\complementaire{X}}(x) \geq r \, \right \}. \\
    \end{array} \]

 The \emph{Hausdorff distance} $d_H(A,B)$ between two subsets $A, B$ of $\R^d$ is the infimum of the set of $t \in \R^+$ such that $B \subset A^{t}$ and $A \subset B^t$. It is also equal to $\norme{d_A - d_B}_{\infty} = \sup_{x \in \R^d} \module{d_A(x) - d_B(x)}$. Given $X$ a compact subset of $\R^d$, the set equality $\overline{\interior(X)} = X$ is equivalent to the Hausdorff convergence of the arbitrarily small counter-offsets to $X$ itself, i.e., $\displaystyle \lim_{\; r \to 0^+} X^{-r} = X$.
Considering sublevel sets of a continuous real-valued function $f$, one has 
$
\overline{f^{-1}(-\infty, t)} = f^{-1}(-\infty, t]$ if and only if there is no local minimum of $f$ with value $t$. In particular, for any compact set $Y \subset \R^d$ and any positive $t$, 
\begin{equation}
\label{eq:offset_domain}
\overline{\interior{(Y^t)}} = Y^t. 
\end{equation}

 A \emph{cone} $A$ in $\R^d$ is a set stable under multiplication by a positive number. Given any subset $B$ of $\R^d$, we denote by $\Cone B$ the smallest cone containing $B$, defined as the image of $[0, \infty) \times B$ by the scalar multiplication map $(\lambda, x) \mapsto \lambda x$.
 The image of $[a,b] \times B$ by this map will be denoted $[a,b] \cdot B$.
    We denote by $\Conv B$ the \emph{convex hull} of $B$. The \emph{dimension} of a cone or a convex set is the dimension of the vector space it spans. 
    The \textit{polar cone} or \textit{dual cone} of a set $B \subset \R^d$, denoted by $B^{\text{o}}$, is the closed convex cone defined by:
    \begin{equation*}
        B^{\text{o}} \coloneqq \{ u \in \R^d \tq \scal{u}{b} \leq 0 \quad  \forall \; b \in B \}.
    \end{equation*}
    When $B$ is a vector space, its polar cone $B^{\text{o}}$ is its orthogonal complement. 
    The polar cone operation is
    an involution 
     on closed convex cones, as it satisfies the identity $(B^{\text{o}})^{\text{o}} = \overline{\Conv \left ( \Cone B \right )}$ for every subset $B$ of $\R^d$. 
\vspace{1em}

 Given a locally Lipschitz function $f: \R^d \to \R$, the \emph{Clarke gradient of $f$ at $x$} is the convex hull of limits of the form $\lim_{n \to \infty} \nabla f (x + h_n)$ where $h_n$ is a sequence converging to $0$ such that the gradient of $f$ exists at $x + h_n$ for every $n \in \N$. We denote the Clarke gradient of $f$ at $x$ by $\clarke f(x)$. In particular, if $f = d_X$ and if $x$ lies outside of $X$, it is known (see for instance \cite{Clarke1975GeneralizedGA}) that $-\clarke d_X(x)$ is the convex hull of the directions to $x$ from the points $z \in X$  such that $d_X(x) = \norme{x-z}$:
    \begin{equation}
    \label{eq:clarke_distance}
         \clarke d_X(x) \coloneqq \Conv \left ( \left \{ \frac{x-z}{\norme{x-z}} \, \Biggr | \, z \in \Gamma_X(x)  \right \} \right ), 
         \end{equation}
    where $\Gamma_X(x)$ is the set of \emph{closest points} to $x$ in $X$ (\Cref{fig:distance}, right). Elements of $\Gamma_X(x)$ will be denoted by the letter $\xi$. In particular, we denote by $\xi_X(x)$ \textit{the} closest point to $x$ in $X$ when $\Gamma_X(x)$ is a singleton.
\vspace{1em}

 Given a subset $X$ of $\R^d$, its \emph{distance to 0} measures how far it is from intersecting $\{ 0\}$. It is defined by $\Delta(X) \coloneqq \inf \left \{ \, \norme{x} \tq  x \in X \right \} = d_X(0).$
We say that $x \in \R^d$ is a \textit{critical point} of a locally Lipschitz function $\phi : \R^d \to \R$ when $0 \in \clarke \phi(x)$, or equivalently when $\Delta(\clarke \phi(x)) = 0$. A number $c \in \R$ is called a \textit{critical value} of $\phi$ when $\phi^{-1} (\{c \})$ contains a critical point, and a \textit{regular value} of $\phi$ otherwise.
The definition for critical points of $\phi$ should not be confused with the critical points of $\phi_{|X}$.

\vspace{1em}

 Given $\mu$ in $(0,1]$, the \emph{$\mu$-reach} of a subset $X$ of $\R^d$ is defined by:
    \begin{equation*}
     \reach_{\mu}(X) \coloneqq \sup \left \{ s \in \R | \; 0 <   \, d_X(x) \leq s \implies \Delta(\clarke d_X(x)) \geq \mu \right \}.   
    \end{equation*}
    Equivalently, having $\reach_{\mu}(X) > 0$ means that for any point $x$ in a certain neighborhood of $X$, the cosines of the half-angle ($\theta$ in (\Cref{fig:distance}), right) between any two points in $\Gamma_X(x)$ are bounded from below by $\mu$. 
    This definition coincides with the classical one found in geometric inference as $\Delta(\clarke d_X(x))$ is the norm of the generalized gradient $\nabla d_X(x)$ defined by Lieutier in \cite{any_open}.
\vspace{1em}
    
    Throughout this article, when no value of $\mu$ has been fixed, for any closed $X \subset \R^d$, \emph{having a positive $\mu$-reach} means that there exists $\mu \in (0,1]$ with $\reach_{\mu}(X) > 0$. The class of sets having a positive $\mu$-reach is certainly broad, intuitively containing stratified sets without concave cusps. A corollary from Fu, Lemma 1.6 of \cite{fuCurvatureMeasuresSubanalytic1994}, is that for any subanalytic set $X \subset \R^d$, the set of values $r > 0$ such that $X^r$
    does not have
     a positive $\mu$-reach is finite.  
\vspace{1em}

The \emph{reach} of a subset $X$ of $\R^d$ is a quantity that was first studied by Federer in \cite{curvature:federer} and that coincides with $\reach_1(X)$. It is the
 supremum of the set of $t \in \R^+$ such that $d_X(x) \leq t$
  implies that $x$ has a unique closest point in $X$. The class of sets with positive reach notably contains convex sets and smooth submanifolds of Euclidean spaces.
Geometric properties of such sets have been studied for a 
long time, and we refer the reader to \cite{reach} for a 
broad overview or \cite{Lytchak} for more in-depth analysis.
In particular, such sets admit curvatures in a generalized sense, a fact that we will use in \Cref{section:morse_complement}.

\begin{figure}[h!]
    \centering
    \includegraphics[width = 0.95\textwidth]
    {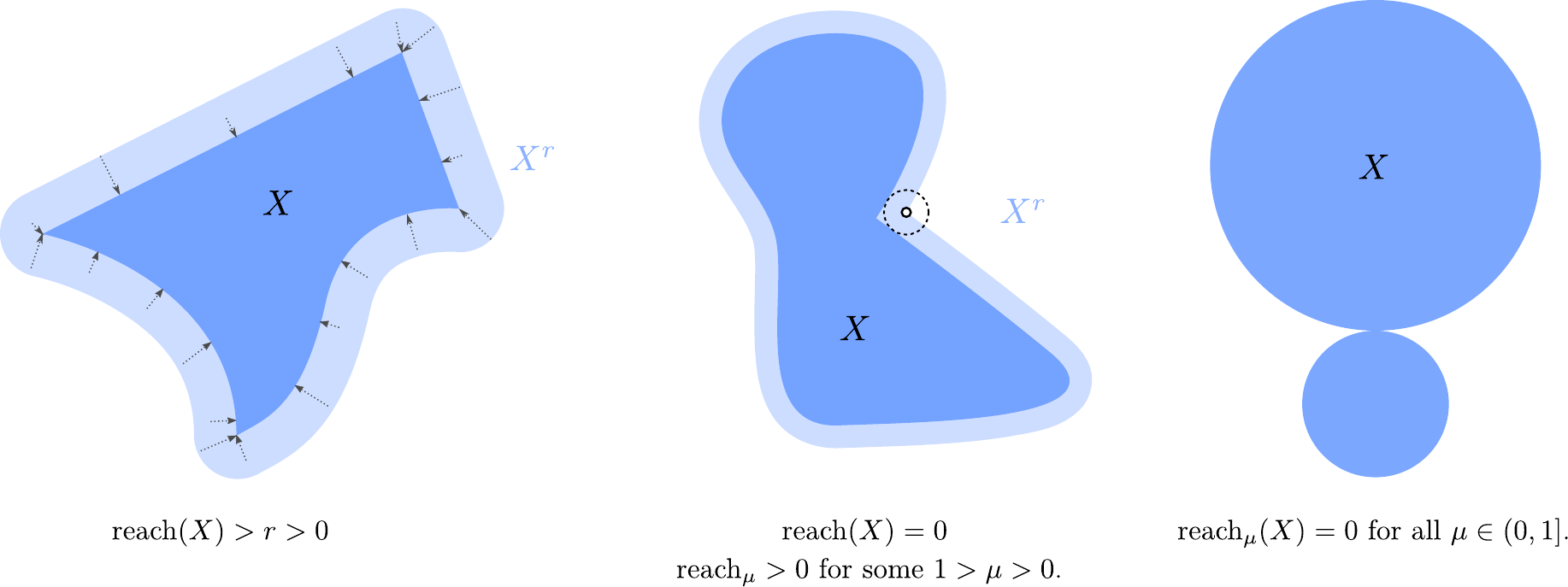}
    \captionsetup{justification=centering}
    \caption{Sets with particular $\reach_{\mu}$.}
    \label{fig:offset_clarke}
\end{figure}

A map $f : U \subset \R^d \to \R$ defined on an open subset $U$ is said to be \textit{semiconcave} when there exists a $K \in \R$ such that the map $x \mapsto f(x) - K \norme{x}^2$ is concave. 
Bangert \cite{Bangert} showed  that any set of the form $f^{-1}[c, + \infty)$ where $f$ is a semiconcave function and $c$ a regular value of $f$  (in the sense of Lipschitz analysis given above, i.e., $\Delta(\clarke f(x)) \neq 0$ for every $x \in f^{-1}(c)$) has positive reach. Since $x \mapsto d^2_X(x) - \norme{x}^2 = \inf_{z \in X} -2\scal{x}{z} + \norme{z}^2$ is concave as the minimum of a family of linear maps indexed over $X$ yields the following theorem. An elementary proof of this fact can be found in  \cite{chazalShapeSmoothingUsing2007}, 4.1.

\begin{theorem}[{Reach of complements of offsets \cite{Bangert}}]
\label{th:reach_complements}
Let $X$ be a compact subset of $\R^d$. If $r > 0$ is a regular value of $d_X$, then \[ \reach(\complementaire{(X^r)}) > 0. \]
    \end{theorem} 

\vspace{1em}

 The \emph{tangent cone} of $X$ at $x$, $\Tan(X,x)$ is defined as the cone generated by the limits $\lim_{n \to \infty} \frac{x_n - x}{\norme{x_n - x}}$, where the sequence $(x_n)_{n \in \N}$ belongs to $X$, converges to $x$ and never takes the value $x$. In that case, we say that $u$ is \emph{represented} by the sequence $(x_n)_{n \in \N}$. 
By a diagonal argument, tangent cones are always closed.

Equivalently, for any $x \in X$, a non-zero $u$ belongs to $\Tan(X,x)$ if and only if
\begin{equation}
\label{eq:zero_derivative}
	\liminf_{t \to 0^+} \frac{1}{t}d_X(x + tu) = 0. 
\end{equation}
We call this the zero-derivative characterization of tangent cones, as it means that the derivative of $d_X$ at $x$ in direction $u$  (if it were to exist) is $0$.

When $X \subset \R^d$ has positive reach, each tangent cone $\Tan(X,x)$ is also convex - see (12), Theorem 4.8 in \cite{curvature:federer}. Federer also showed (statement (11) in Theorem 4.8 of \cite{curvature:federer}) that for any non-zero $u$, $\mathrm{liminf}_{t \to 0} \frac{1}{t} d_X(x + tu) = 0$ implies $\lim_{t \to 0} \frac{1}{t}d_X(x + tu) = 0$, so that when $X$ has positive reach, 
\begin{equation}
\label{eq:zero_derivative_reach}
\Tan(X,x) = \left \{ u \in \R^d, \lim_{t \to 0^+} \frac{1}{t} d_X(x + tu) = 0 \right \}.
\end{equation}
\vspace{1.5em}
    
 When $X$ has positive reach, we define its \textit{normal cone} at $x \in \partial X$, denoted by $\Nor(X,x)$, as the dual cone to the tangent cone at $x$:
    \[ 
    \begin{array}{lr}
        \Nor(X,x) \coloneqq \Tan(X,x)^{\text{o}}.
    \end{array}\] 
   
    Federer also shows that $\Nor(X,x)$ contains at least one non-zero vector, and that it is related (see 4.8, \cite{curvature:federer}) to the projection to the closest point in $X$ function $\xi_X$  by the following characterisation, for any $0 < t < \reach(X)$:
    \begin{equation}
    \label{eq:proj-char}
     \emptyset \; \subsetneq \; \Nor(X,x) \cap \Sphere^{d-1} = \left \{ u \in \Sphere^{d-1} \,  \Bigr | \; \xi_X(x + tu) = x \right \}. 
\end{equation}    

    \begin{figure}[!htb]
       \begin{minipage}{0.46\textwidth}
         \centering
        \includegraphics[width=.68\linewidth]{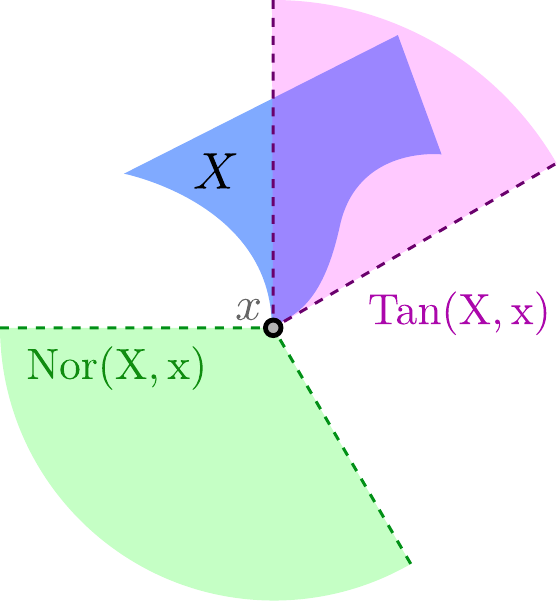}
        \caption{Tangent and normal cones of $X$ at $x$ when $\reach(X) >0$.}\label{Fig:Data2}
       \end{minipage}\hfill
       \begin{minipage}{0.46\textwidth}
         \centering
         \includegraphics[width=.7\linewidth]{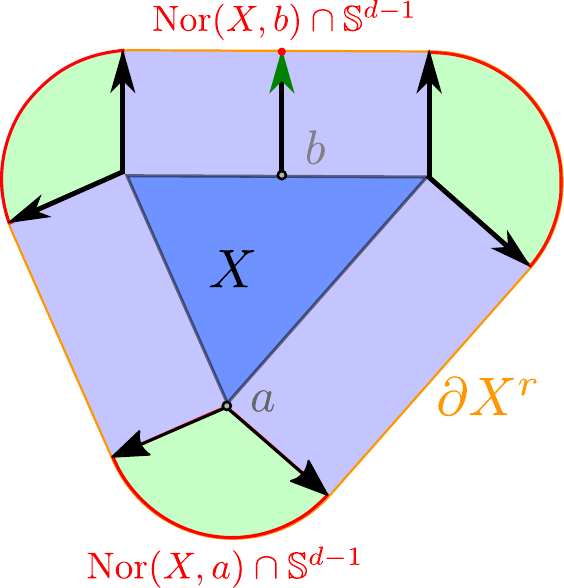}
         \caption{Some unit normal cones (in red) when $0 < r < \reach(X)$.}\label{Fig:Data1}
       \end{minipage}
    \end{figure}
    
\vspace{1em}
 If $X \subset \R^d$ has positive reach, we say that $X$ is \emph{fully dimensional}  when $\Tan(X,x)$ has dimension $d$ for every $x \in \partial X$. This is equivalent to having the set equality $\overline{\interior(\Tan(X,x))} = \Tan(X,x)$ for all $x \in \partial X$. In particular, a set with positive reach is fully dimensional if and only if it is a Lipschitz domain \cite{Lytchak}.

\subsection{Normal bundles} \label{sub:normal}

We are now in position to define the \emph{normal bundle} of sets with positive reach or whose complement sets have positive reach. 

    \begin{definition}[Sets admitting a normal bundle] Let $X \subset \R^d$.
When $\complementaire{X} \coloneqq \overline{\R^d \setminus X}$ has positive reach and is a Lipschitz domain, for all $x \in \partial X$ define 

    \[ \Nor(X,x) \coloneqq - \Nor(\complementaire{X},x). \]
    This definition is consistent in the case both $X, \complementaire{X}$ have positive reach - see \Cref{cor:consistency}.\\

 If either $\reach(X) > 0$ or both $\reach(\complementaire{X}) > 0$ and $X$ is a Lipschitz domain, we say that \emph{$X$ admits a normal bundle $\Nor(X)$} with
    \[\Nor(X) \coloneqq \bigcup_{x \in \partial X} \{ x \} \times (\Nor(X,x) \cap \Sphere^{d-1}). \]
\end{definition}

Normal bundles have intrinsic dimension $(d-1)$, in the following sense. 
\begin{proposition}[Normal bundles are Lipschitz submanifolds of $\R^d \times \Sphere^{d-1}$]
When either $X$ or $\complementaire{X}$ has positive reach, $\Nor(X)$ is a $(d-1)$-Lipschitz submanifold of $\R^d \times \Sphere^{d-1}$. 
\end{proposition}

\begin{proof}
When $X$ has positive reach, the result is already known (see for instance Lemma 4.21 in \cite{Ja-2019}).
Otherwise, $\Nor(X) = \rho(\Nor(\complementaire{X}))$ is the image of a Lipschitz submanifold under the bilipschitz map $\rho :(x,n) \mapsto (x,-n)$.
\end{proof}

    \begin{figure}[h!]
    \centering
    \includegraphics[scale = 0.65]{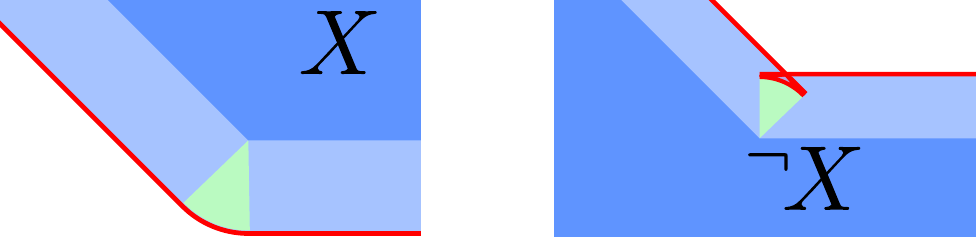}
    \caption{Normal bundles (in red) of a set of positive reach (left) and its complement set $\complementaire{X}$.}
    \end{figure}
    
Since $\Nor(X)$ is a $(d-1)$-Lipschitz submanifold of $\R^d \times \R^d$, its tangent cones are $\mathcal{H}^{d-1}$-almost everywhere vector spaces, prompting the definition of \textit{regular pairs}.

\begin{definition}[Regular pairs]
\label{def:regular}
Any pair $(x,n) \in \Nor(X)$ such that the tangent cone $\Tan(\Nor(X), (x,n))$ is a vector subspace of dimension $(d-1)$ in $\R^d \times \R^d$ is said to be a \emph{regular pair} of $\Nor(X)$. The set of regular pairs has full $\mathcal{H}^{d-1}$-measure in $\Nor(X)$.
\end{definition}

At a regular pair $(x,n)$, the tangent cone (in fact, tangent space) of $\Nor(X)$ enjoys a particular structure
stemming from the more general concept of \emph{normal cycle} of a set \cite{zahleCurvaturesCurrentsUnions1987,fuCurvatureMeasuresSubanalytic1994}.
  While we do not need to write our hypothesis using this more involved language, in our case the normal bundle $\Nor(X)$ is the support of a $(d-1)$-Legendrian cycle over $\R^d \times \Sphere^{d-1}$, whose tangent spaces' structure is already known. This fact will be used in computations involving curvatures in \Cref{section:morse_complement}.\\

\begin{proposition}[Tangent spaces of normal bundles \cite{Ja-2019}]
\label{prop:tangent_spaces}
Let $X$ be a compact set admitting a normal bundle $\Nor(X)$. Then for any regular pair $(x,n)$ in $\Nor(X)$, there exist 
\begin{itemize}
    \item A family
$\kappa_1, \dots , \kappa_{d-1}$ in $\R \cup \{ \infty \}$ called \emph{principal curvatures at $(x,n)$};
    \item A family $b_1, \dots, b_{d-1} \in \R^d$ of vectors orthogonal to $n$ called \emph{principal directions} at $(x,n)$ such that the family $\left ( \frac{1}{\sqrt{1 + \kappa^2_i}}b_i , \frac{\kappa_i}{\sqrt{1+\kappa_i^2}}b_i \right )_{1 \leq i \leq d-1}$ forms an orthonormal basis of $\Tan(\Nor(X), (x,n))$.
\end{itemize}
Moreover, the family of principal curvatures is unique up to permutations.
\end{proposition}

These principal curvatures coincide with the ones found in differential geometry as eigenvalues of the second fundamental form. Indeed, assume that $X \subset \R^d$ is bounded by a $C^{1,1}$-hypersurface, i.e the boundary of $X$ is a hypersurface such that the Gauss map $x \in \partial X \mapsto n(x) \in \Sphere^{d-1}$ is Lipschitz. The pair $(x,n(x)) \in \Nor(X)$ is regular if and only if $n$ is differentiable at $x$ 
(see 4.2 in \cite{fuCurvatureMeasuresGeneralized1989}).
 In that case, its differential is symmetric and its eigenvalues counted with multiplicity (resp. orthonormal basis of eigenvectors) are the principal curvatures (resp. principal directions) at $(x,n(x))$.

\subsection{Critical points and Hessians for $f_{|X}$} \label{sub:critical}

In \cite{fuCurvatureMeasuresGeneralized1989}, Fu defines a notion of Morse functions over sets of positive reach.
  Below in \Cref{def:crit}, we gather from (3.11) and (4.5) in \cite{fuCurvatureMeasuresGeneralized1989} the definitions of critical points of smooth functions, second fundamental forms, Hessians and non-degenerate critical points of restricted functions, all adapted to our notations. These definitions were originally given for sets with positive reach, but naturally extend to any set admitting a normal bundle. 

The projection $\R^d \times \R^d \to \R^d$ onto the first factor is denoted by $\pi_0$.
When $(x,n)$ is a regular pair, observe that 
while $\Tan(\Nor(X), (x,n))$ is a $(d-1)$-dimensional subspace of $\R^d \times \R^d$, $\pi_0(\Tan(\Nor(X), (x,n)))$ is a subspace of $\R^d$. From \Cref{prop:tangent_spaces}, its dimension is given by the number of finite principal curvatures of $\Nor(X)$ at $(x,n)$, which may be strictly less than $d-1$.
\\

\begin{definition}[Critical points and Hessian]
\label{def:crit}
Let $X$ be a compact set of $\R^d$ admitting a normal bundle, $f: U \to \R$ be a smooth map defined on an open set $U$ containing $X$.
\begin{itemize}
    \item Let $(x,n) \in \Nor(X)$ be a regular pair as in 
    \Cref{def:regular}. The second fundamental form $\sff_{x,n}$ of $X$ at $(x,n)$ is defined as the bilinear form on $\pi_0(\Tan(\Nor(X), (x,n)))$ such that for every pair $(u,v), (u',v')$ in $\Tan(\Nor(X), (x,n))$,
\begin{equation*}
    \sff_{x,n}(u,u') \coloneqq \scal{u}{v'}.
\end{equation*}
This map is well-defined and symmetric thanks to the structure of tangent spaces described in \Cref{prop:tangent_spaces}. Indeed,
    taking $(b_i)$ an orthonormal basis of $\pi_0(\Tan(\Nor(X), (x,n)))$ consisting of all principal directions with finite associated principal curvatures, this definition is equivalent to 
\begin{equation*}
\sff_{x,n}(b_i,b_j) \coloneqq \left \{
\begin{array}{cr}
     \kappa_i & \text{ if } i = j \\
     0 & \text{ if } i \neq j,
\end{array}
\right.
\end{equation*}  
and generalizes the classical second fundamental form obtained when $X$ has a smooth boundary;
    \item We say that $x \in X$ is a \emph{critical point} of $f_{|X}$ when 
	either $x \in \partial X$ and   
    $\nabla f(x) \in -\Nor(X,x)$, or $x \in \interior{(X)}$ and $\nabla f(x) = 0$;
    \item We say that $c \in \R$ is a \emph{critical value} of $f_{|X}$ when $f^{-1}(c)$ contains at least one critical point of $f_{|X}$. Otherwise, $c$ is a \emph{regular value} of $f_{|X}$;
    
    \item If $x \in \partial X$ is a critical point of $f_{|X}$ with $\nabla f(x) \neq 0$, let $n \coloneqq -\frac{\nabla f(x)}{\norme{\nabla f(x)}}$.
    When $(x,n)$ is a regular pair, \emph{the Hessian of $f_{|X}$ at $x$} denoted by $H_x f_{|X}$ is defined over  $\pi_0(\Tan(\Nor(X), (x,n)))$  by: 
\[ H_x f_{|X}(u,u') \coloneqq H_x f(u,u') + \norme{\nabla f(x)} \sff_{x,n}(u,u'). \]
	If $x \in \interior{(X)}$ is a critical point of $f_{|X}$, simply let $H_x f_{|X} \coloneqq H_x f$, defined on $\R^d$.

\item The \emph{index} of $H_x f_{|X}$ is the dimension of the largest subspace on which it is negative definite;

\item We say that a critical point $x \in \partial X$ of $f_{|X}$ is \emph{non-degenerate} when $\nabla f(x) \neq 0$, $(x,n)$ is a regular pair of $\Nor(X)$, and the Hessian $H_x f_{|X}$ is not degenerate. 
A critical point $x \in \interior(X)$ is non-degenerate when $H_x f{|X}$ is non-degenerate.

\item $f_{|X}$ is said to be \emph{Morse} when its critical points are non-degenerate.
\end{itemize}
\end{definition}

Using the definitions of \Cref{def:crit},
 Fu proved the Morse theorems for sets with positive reach, as stated in Theorem 4.8 of \cite{fuCurvatureMeasuresGeneralized1989} and its following corollary. 
These results can be summarized as follows.
\begin{theorem}[\textit{Generalized Morse theory for sets with positive reach }{\cite{fuCurvatureMeasuresGeneralized1989}}]
\label{th:fu}
 Let $X$ be a compact subset of $\R^d$ with positive reach and let $f: \R^d \to \R$ be a smooth function such that $f_{|X}$ is \textit{Morse} with at most one critical point per level set.\\

 Then for any regular value $c \in \R$, $X_c$ has the homotopy type of a $CW$-complex with one $\lambda_p$ cell for each critical point $p$ such that $f(p) < c$, where 
 \[ \lambda_p = \text{Index of } H f_{|X} \text{ at } p.\]
\end{theorem}

\begin{remark}[Motivation behind the definition of critical points]
The definition of critical points of $f_{|X}$ in $\partial X$ as points $x$ satisfying $-\nabla f(x) \in \Nor(X,x)$ is motivated as follows.

When $X$ is a smooth submanifold of $\R^d$ this definition coincides with the traditional definition of differential geometry. Indeed, viewing $X$ as a manifold and $f_{|X} : X \to \R$ as a map between manifolds, the (manifold) gradient of $f_{|X}$ at $x$ is the orthogonal projection onto the tangent space $\Tan(X,x)$, which is zero exactly if $\nabla f(x)$ belongs to the space of vectors at $x$ orthogonal to $\Tan(X,x)$, which is $\Tan(X,x)^{\text{o}} = \Nor(X,x)$.

When $X \subset \R^d$ is merely a compact subset of $\R^d$, it is a priori not clear what defines a critical point for $f_{|X}$. In stratified Morse theory, a point $x$ is critical for $f_{|X}$ if it is a critical point of $f_{|S}$ for some stratum $S$ of $X$ containing $x$. However, this definition allows for the absence of topological events at critical values: take for instance $X$ the unit ball of $\R^2$ and $f : \R^2 \to \R$ the height function, as illustrated in \Cref{fig:motivation_critique} below. $X$ can be decomposed in two strata, the open ball and the unit circle. There is only one change in topology among the height filtration restricted to $X$, yet the point $x_{\max}$ is critical for the restriction of $f$ to the circle.\\
   
    \begin{figure}[!htb]
             \centering
        \includegraphics[width=0.6\textwidth]{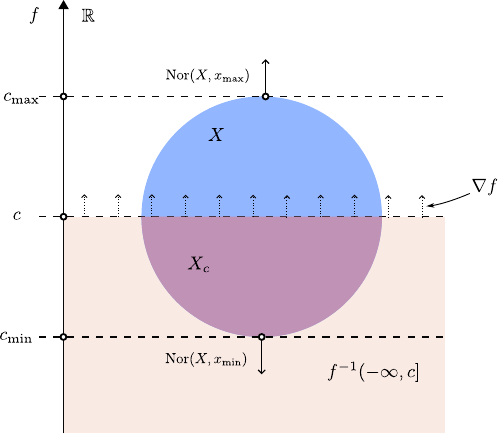}
        \caption{Sublevel set filtration $X_c = f^{-1}(-\infty,c]$ where $X \subset \R^2$ is a ball and $f$ a height function. The only point of $f_{|X}$ around which a change of homotopy type happens is $x_{\min}$, where the minimum of $f_{|X}$ is reached. There is no change of homotopy type around the maximum value $c_{\max}$ of $f_{|X}$, where $\nabla f(x_{\max})$ belongs to $\Nor(X,x_{\max})$ but not $\Nor(X, x_{\min})$.
}
	\label{fig:motivation_critique}
    \end{figure}

When $X$ is a complementary regular set and $c$ is a regular value of $f_{|X}$ in our sense, then there is no topological change of homotopy type for the filtration $(X_t)_{t \in \R}$ at filtration value $c$. 
The intuition behind this fact (proved in \Cref{sub:isotopie}) is the following:
when $x$ belongs to the boundary of $X$, we show that $-\nabla f(x) \notin \Nor(X,x)$ is equivalent to $0 \notin \Cone(\clarke d_X(x)) + \nabla f(x)$. Intuitively, this means that there exist directions around $x$ making both $d_X$ and $f$ decrease. This is reflected mathematically by
studying the non-negative function $\phi = d_X + \max(0, f - c)$. Assuming $c$ is a regular value, we are able to construct a deformation retraction of some sublevel set $\phi^{-1}[0, K]$ onto $X_c = \phi^{-1}( 0 )$, and from this, a deformation retraction of $X_{c +\eps}$ onto $X_c$ when $\eps$ is small enough.

\end{remark}

\subsection{Clarke gradients and approximate flows}\label{sub:clarke}
We use a classical tool in the analysis of Lipschitz functions called the \textit{Clarke gradient}. We recall its definition and study some of its properties. We refer the reader to the original article from Clarke \cite{Clarke1975GeneralizedGA} for the properties we do not prove.

\begin{definition}[Clarke gradients of locally Lipschitz functions]
\label{def:clarke}
Let $\phi: U \to \R$ be a locally Lipschitz function 
on an open set $U$.
  Its \emph{Clarke gradient}
 $\clarke \phi$ is a set-valued map defined at $x$ as the convex hull of limits of gradients taken in an arbitrarily close neighborhood of $x$, that is,
\[ 
\clarke \phi (x) \coloneqq \Conv \left ( \lim_{i \to \infty} \nabla \phi(x_i) \tq x_i \in \R^d \to x, \phi \text{ differentiable at } x_i \text{ for all } i \right ). 
\]
\end{definition}

Each $\clarke \phi(x)$ is a closed convex subset of $\R^d$.
Whenever we refer to the explicit definition of the Clarke gradient, the fact that $\phi$ needs to be differentiable at any $x_i$ will be implied. \\

Here are some properties of the Clarke gradient.

\begin{proposition}[Some properties of the Clarke gradient]
\label{prop:clarke_prop}
Let $\phi : U \to \R$ be a locally Lipschitz function defined
 on an open subset $U$ of $\R^d$.
\begin{itemize}
\item By Rademacher's theorem, $\clarke \phi(x)$ is non-empty for all $x$, and $\clarke \phi(x) \subset B(0,R)$ if $\phi$ is $R$-Lipschitz in a neighborhood of $x$; 
\item When $\phi$ is 
continuously differentiable at $x$, we have
\[ \clarke \phi (x) = \{ \nabla \phi(x) \}; \]
\item 
If $A \subset \R^d$ has density $1$ at $x$, then
\begin{equation}
\label{eq:clarke_densite}
\clarke \phi (x) = \Conv \{ \nabla \phi(x_i) \tq x_i \to x, x_i \in A  \};
\end{equation} 
\item
If $\phi = \max(f, g)$ where $f,g$ are locally Lipschitz maps $\R^d \to \R$ such that the set {$\{y \tq  f(y) =  g(y) \}$} has density 0 at $x$, we have:
\begin{equation}
\label{eq:clarke_max}
\mathmakebox[1.15\linewidth][c]{
    \displaystyle
        \clarke \phi(x) = \Conv \left (\{ \lim_{i \to \infty} \nabla f(x_i) \tq  f(x_i) > g(x_i) \} \cup \{ \lim_{i \to \infty} \nabla g(x_i) \tq g(x_i) > f(x_i) \} \right ),
}
\end{equation}
where $x_i$ are sequences in $\R^d$ converging to $x$;
\item If $\phi = \phi_1 + \phi_2$ where $\phi_1, \phi_2$ are locally Lipschitz maps $\R^d \to \R$, then we have for every $x \in \R^d$:
\[
\clarke \phi(x) \subset \clarke \phi_1(x) + \clarke \phi_2(x),
\]
where the addition between sets is $A + B = \{ a + b, a \in A, b \in B \}$.
Moreover, if $\phi_2$ is smooth around $x$, we have $\clarke \phi(x) = \clarke \phi_1(x) + \{ \nabla \phi_2(x) \}$.

\item When there exists a positive real number $C$ such that $x \mapsto \phi(x) - C \norme{x}^2$ is concave (i.e., $\phi$ is $C$-semiconcave), the Clarke gradient of $\phi$ at $x$ is, for any open neighborhood $V \subset U$ of $x$:

\begin{equation}
\label{eq:semi-concave}
\clarke \phi(x) = \left \{ p \in \R^d, \forall \,  y \in V, \phi(y) \leq \phi(x) + \scal{p}{y - x} + \frac{C}{2}\norme{x - y}^2  \right \}. 
\end{equation}

\item The map $\clarke \phi$ is upper semicontinuous in the following sense: for each $x \in \R^d$ and $\eps > 0$, there exists a radius $\delta > 0$ such that for every $y \in B(x, \delta)$, $\clarke \phi(y) \subset (\clarke \phi(x))^{\eps}$ (the $\eps$-offset of $\clarke \phi(x)$). 

\end{itemize}
\end{proposition}

\vspace{0.2cm}
For our purposes, an important corollary of the last point is the lower semicontinuity of maps of the form $x \mapsto \Delta(\clarke \phi(x))$.\\
\begin{proposition}[Semicontinuity of the Clarke gradient]\label{prop:openness}
Let $\phi: U \to \R$ be a locally Lipschitz function
defined on an open subset $U$ of $\R^d$.
If a sequence $(x_i)_{i \in \N}$ in $U$ converges to some $x \in U$, we have
\[ 
\liminf_{i \to \infty} \Delta \left ( \clarke \phi(x_i) \right ) \geq \Delta \left ( \clarke \phi(x) \right).
\]
\end{proposition}

\vspace{0.25cm}

 Assuming $\clarke \phi(x)$ stays uniformly away from $0$, we are able to build deformation retractions between the sublevel sets of $\phi$ using approximations of what would be the normalized flow of $-\phi$ had it been smooth.

\vspace{0.3cm}
\begin{proposition}[Approximate inverse flow of a Lipschitz function]
\label{prop:flot}
 Let $a < b \in \R$. Let $\phi: U \to \R$ be a function defined on an open subset of $\R^d$ and locally Lipschitz on an open neighborhood
of  $\phi^{-1}(a,b]$. Assume that \[ \inf \{ \Delta( \clarke \phi(x)) \tq x \in \phi^{-1}(a,b] \} = \mu > 0. \] Then for every 
 $\varepsilon \in (0, \mu)$, there exists a continuous function
 \[C_{\phi}: \left \{
 \begin{array}{ccc}
      [0,1] \times \phi^{-1}( -\infty, b]& \to & \phi^{-1}(-\infty, b]  \\
      (t,x) & \mapsto & C_{\phi}(t,x) 
 \end{array} \right.
 \]
 such that
 \begin{itemize}
     \item For any $s > t$ and $x$ such that $C_{\phi}(s,x) \in \phi^{-1}(a,b]$, we have
 \begin{equation*}
\phi(C_{\phi}(s,x)) - \phi(C_{\phi}(t,x)) \leq -(s-t)(b-a);
 \end{equation*} 
 \item  For any $t \in [0,1]$ and $x \in \phi^{-1}(-\infty, a]$, we have $C_{\phi}(t,x) = x$;
 \item For any $x \in \phi^{-1}( -\infty,b]$, the map $s \mapsto C_{\phi}(s,x)$ is $\frac{b-a}{\mu - \eps}$-Lipschitz.
 \end{itemize}
 
 In particular, the sublevel set $\phi^{-1}(-\infty, a]$ is a deformation retract of $\phi^{-1}(-\infty,b]$.
\end{proposition}

\begin{proof}
A weaker form of this claim can be found in Section D of \cite{kimHomotopyReconstructionCech2020}, and the method can be traced back to \cite{Grove}. Here the constants have been optimized and the proposition generalized to Lipschitz functions. The idea is (1) to determine the best direction $W$ to locally make $f$ decrease, (2) to interpolate these directions to build a vector field $V$ from which a flow is obtained, (3) to show that the trajectories of this flow satisfy the wanted rate of decrease for $f$, (4) to show that the extended flow is continuous. For the sake of completeness, we provide a full proof.\\

Let $\varepsilon > 0$. By the upper semicontinuity of the Clarke gradient, for any $x \in \phi^{-1}(a,b]$ we can consider $B_x \subset \phi^{-1}(a, + \infty)$ an open ball centered in $x$ such that $\clarke \phi(y) \subset \clarke \phi(x)^{\varepsilon}$ (the $\eps$-offset of $\clarke \phi(x)$) for any $y \in B_x$.
Since $\clarke \phi (x)$ is a closed convex set, there is a unique point $W(x)$ in $\clarke \phi(x)$ realising the distance to $0$, i.e., $\norme{W(x)} = \Delta(\clarke \phi(x))$. 
Since $W(x)$ is the closest point to 0 in $\clarke \phi(x)$, it follows from the convexity of $\clarke \phi(x)$ that we have:
\begin{equation*}
\forall \, u \in \clarke \phi(x), \scal{u}{W(x)} \geq \norme{W(x)}^2.    
\end{equation*}

The family $\{ B_x \}_{x \in \phi^{-1}(a,b]}$ is an open covering of the open set $B = \bigcup_{x \in \phi^{-1}(a,b]} B_x$ which contains $\phi^{-1}(a,b]$.
As an open set of $\R^d$, $B$ is paracompact, so that there exists a locally finite smooth partition of unity $(\rho_i)_{i \in I}$ subordinate to the previous family, i.e., such that the support of each $\rho_i$ is included in one of the balls $B(x_i)$ with $x_i \in \phi^{-1}(a,b]$. Use them to define the vector field $V$ as a smooth convex interpolation of normalized $-W$:
\begin{equation}
\label{eq:def_V}
V(y) \coloneqq -\sum_{i \in I} \rho_i(y)\frac{W(x_i)}{\norme{W(x_i)}}.  
\end{equation}
 By construction, $\norme{V}_{\infty} \leq 1$ and $V$ is locally Lipschitz. 
Take $C$ to be the maximal flow of $V$, i.e., a function satisfying the ordinary differential equation $\frac{\partial}{\partial t} C(t,x) = V(C(t,x))$ defined on a maximal open domain $\mathbb{D}$ in $\R^+ \times U$. For any $y \in \phi^{-1}(a,b]$, the only non-zero terms in the sum of \Cref{eq:def_V} are from points $x_i$ such that $y \in B_{x_i}$, so that $\clarke \phi(y) \subset (\partial \phi(x_i))^{\eps}$.
 We thus get
\begin{equation}
\label{eq:clarke_vs_flow}
\forall \zeta \in \partial \clarke \phi(x), \; \scal{V(x)}{\zeta} \leq - \sum_{i \in I} \rho_i (x) \left ( \, \norme{W(x_i)} - \eps \, \right ) \leq -\mu + \eps.
\end{equation} 
The maximal subset of $\R^+$ for which the flow 
of $V$
starting at $x$ is defined is open and connected in $\R^+$, so that it is of the form $[0, s_x)$, where $s_x \in (0, \infty]$.
 Since the trajectory $s \mapsto C(s, x)$ has gradient $V(C(s,x))$ of norm at most 1, it
  is 1-Lipschitz.
Assuming $s_x$ is finite,
  we can define by completeness $C(s_x,x)$ as the endpoint of this curve, that is, $ \displaystyle C(s_x,x) = \lim_{s \to s_x} C(s,x)$. 

Now fix $x \in \phi^{-1}(a, b]$. We want to study how the map $\psi : s \mapsto \phi(C(s,x)), [0, s_x) \to (a,b]$ decreases. Since $s \mapsto C(s,x)$ is a 1-Lipschitz curve, $\psi$ is Lipschitz and thus differentiable almost everywhere in $[0, s_x)$. Assume $\psi$ is differentiable at some $t > 0$; since $s \mapsto C(s,x)$ satisfies a first-order differential equation, write the first-order expansion around $t$:
\begin{align*}
\psi(t+h) = \phi(C(t + h,x)) = &  \quad \phi( C(t,x) +  h \frac{\partial}{\partial t} C(t,x) + o(h)) \\
= & \quad  \phi( C(t,x) + h V(C(t,x))) + o(h).
\end{align*}
Since $\psi$ is differentiable at $t$, $\phi$ admits a directional derivative at $C(t,x)$ in direction $V(C(t,x))$ which we denote by $\phi'(C(t,x),V(C(t,x)))$; seeing the first order-expansion above this directional derivative coincides with $\psi'(t)$.
From Proposition 1.4 in \cite{Clarke1975GeneralizedGA}, we know that the support function associated to the Clarke gradient is an upper bound for directional derivatives:
\begin{align*}
\psi'(t) = & \; \phi'(C(t,x),V(C(t,x))) \leq \max \left \{ \scal{\zeta}{V(C(t,x))} \tq \zeta \in \clarke \phi(C(t,x)) \right \}.
\end{align*}

By \Cref{eq:clarke_vs_flow}, we thus have for almost every $t \in (0, s_x)$
\[ \psi'(t) \leq - \mu + \eps. \]

When $s \leq t$ both belong to $[0, s_x)$ we can integrate the previous inequality to obtain:
\begin{equation}
	\label{eq:decrease_bound}
    \phi(C(t,x)) - \phi(C(s,x)) \leq -(\mu - \varepsilon)(t-s).
\end{equation}
In particular, $s \mapsto \phi(C(s,x))$ is a decreasing map on $[0, s_x)$, so that trajectories starting at a point $x \in \phi^{-1}(a, b]$ stay in $\phi^{-1}(a,b]$. The bound on the rate of decrease \Cref{eq:decrease_bound} implies that $s_x \leq \frac{\phi(x)-a}{\mu - \eps} < +\infty$ for all $x \in \phi^{-1}(a,b]$, and consequently $\phi(C(s_x,x)) = a$.  \\

Consider $C^*$ the restriction of $C$ to points in $\phi^{-1}(a,b]$ and extended by the identity on $\phi^{-1}(-\infty, a]$:
\[
C^*(t,x) \coloneqq \left \{ 
\begin{array}{l l}
    C(\min(t,s_x),x) & \text{ when } a < \phi(x) \leq b,\\
    x & \text{ else.} \\
\end{array}
\right. \]
We will now show that $C^*$ is continuous at every point $(s,x) \in \R^+ \times \phi^{-1}(-\infty, b]$. $C^*$ is continuous inside $\R^+ \times \phi^{-1}(a,b]$, where it is even smooth. $C^*$ is continuous inside the open set $\R^+ \times \phi^{-1}(-\infty, a)$ since in this set $C^*(t,x) = x$.
 We now turn our attention to the other points,
 which we divide into either (1) the set of pairs $(s,x)$ with $\phi(x) > a$ and $s \geq s_x$, or (2) those with $\phi(x) = a$.
 Begin by studying a pair $(s,x)$ belonging to the set (1), recalling that $s_x > 0$. The idea is to aggregate the facts that $C$ is continuous inside the maximal set $\mathbb{D}$, and that each trajectory $s \mapsto C(s,y)$ leading to the endpoints $C(s_y,y)$ is 1-Lipschitz.
  For every $\delta \in (0, \phi(x) - a)$, consider the time parameter $s^{\delta}_x > 0$ defined by $\phi(C(s^{\delta}_x, x)) = a + \delta$. Since $\phi$ is continuous and strictly decreasing on flow trajectories, $s^{\delta}_x$ is well-defined. Now since $C$ is continuous inside $\mathbb{D}$, there exists a radius $\rho_{\delta} > 0$ such that for every $y \in B(x, \rho_{\delta})$, we have 
  \[
  \begin{cases}
  a + \delta/2 \leq \phi(C(s^{\delta}_x, y)) \leq a + 3\delta/2 \\
  \delta \geq \norme{C(s^{\delta}_x, y) - C(s^{\delta}_x, x)}.
  \end{cases} \]
   Thanks to the first inequality and the bound on the rate of decrease of $\phi$ (\Cref{eq:decrease_bound}), for every $y \in B(x, \rho_{\delta})$ we have $0 < s_y - s^{\delta}_x \leq \frac{3\delta}{2(\mu - \eps)}$ and $s_x \leq s^{\delta}_x + \frac{\delta}{\mu - \eps}$. 
   Now let $(t,y)$ be such that $\module{s - t} \leq \module{s_x - s^{\delta}_{x}}$,
    $\norme{x - y} \leq \rho_{\delta}$. Since $C^*(s,x) = C^*(s_x,x)$ and since maps of the form $t \mapsto C^*(t,z)$ are 1-Lipschitz we end up with
   \begin{align*}
    \norme{C^*(t,y) - C^*(s,x)} \leq \; & 
     \norme{C(\min(t, s_y),y) - C(s^{\delta}_{x},y)} \\ & + \norme{C(s^{\delta}_{x},y) - C(s^{\delta}_{x},x)} \\ &+ \norme{C(s^{\delta}_{x},x) - C(s_x,x)}  \\
     \leq \; & \; \module{\min(t, s_y) - s^{\delta}_{x}} + \delta + \module{s^{\delta}_{x} - s_x}.
\end{align*} 
Since we have taken $t$ to be greater than $s^{\delta}_{x}$, we have $\module{\min(t, s_y) - s^{\delta}_{x}} \leq \module{s_y - s^{\delta}_x} \leq \frac{3\delta}{2(\mu - \eps)}$. Also aggregating the bound on $s_x - s_x^{\delta}$ we end up with
\begin{equation*}
\norme{C^*(t,y) - C^*(s,x)} \leq \frac{2\delta}{(\mu - \eps)} + \delta,
\end{equation*}
thereby proving the continuity of $C^*$ at $(s,x)$. 

Now assume that $\phi(x) = a$, i.e., $x$ lies in set (2). Equivalently $C^*(s,x) = x$ for every $s \geq 0$. Since $\phi$ is continuous at $x$ and greater than or equal to $a$, for every $\delta > 0$, 
there exists $\rho'_{\delta}$ such that for any point $y$ lying at distance $\rho_{\delta}$ to $x$, either $\phi(y) \leq a$ or $a < \phi(y) < a + \delta$, in which case $s_y \leq \frac{\delta}{\mu - \eps}$, so that $\norme{C^*(t,y) - y} \leq \frac{\delta}{\mu - \eps}$ since the trajectory is 1-Lipschitz and constant after $s_y$. Finally, for every $y$ in a ball of radius $\min(\rho_{\delta}, \delta)$ we have for every $t \geq 0$
\begin{equation*}
\norme{C^*(t,y) - C^*(s,x)} \leq \norme{C^*(t,y) - y} + \norme{x - y} \leq \frac{\delta}{\mu - \eps} + \delta, 
\end{equation*}
thereby proving the continuity of $C^*$ at $(s,x)$.
To reduce the time domain to $[0,1]$, since each $s_x$ is at most $(b-a)/(\mu -\eps)$ we let $C_{\phi}(t,x) = C^*(t(b-a)/(\mu-\eps), x)$ and obtain the desired map.

\end{proof}

    \begin{figure}[!htb]
       \begin{minipage}{0.46\textwidth}
         \centering
        \includegraphics[width=.68\linewidth]{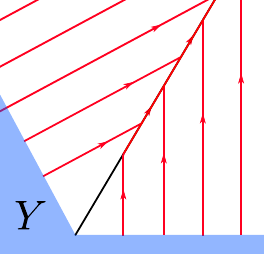}
        \caption{Exact flow of the locally semiconcave map $d_Y$.}\label{Fig:generalise}
       \end{minipage}\hfill
       \begin{minipage}{0.46\textwidth}
         \centering
         \includegraphics[width=.7\linewidth]{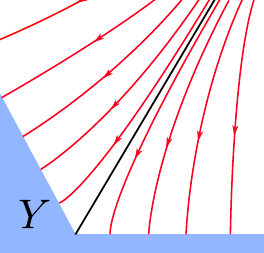}
         \caption{Approximate inverse flow of $d_Y$ outside of $Y$.}\label{Fig:inverse}
       \end{minipage}
    \end{figure}
    
\begin{remark}[On inverse flows of non-smooth maps]
It is a well-known fact in Alexandrov's geometry \cite{Alexandrov} and in the optimal control community \cite{Cannarsa} that any semiconcave map $f : \mathrm{dom} \, f \subset \R^d \to \R$ admits a continuous flow $\mathrm{dom} \, f \times \R^+ \to \mathrm{dom}\, f$ whose trajectories make $f$ increase at a rate $\Delta( \clarke f)$ when parametrized by the arc-length. Semiconcave maps notably include functions of the form $(d_K)^2 + \phi$ where $K$ is a compact subset of $\R^d$ and $\phi$ any $C^2$ map, or $(d_K)_{|\R^d \setminus V}$, where $V$ is an open neighborhood of $K$. However, when $-f$ is not semiconcave this exact flow is not bijective and cannot be reversed to obtain a continuous flow making $f$ decrease at rate $\Delta( \clarke f)$, as in \Cref{Fig:generalise}. \Cref{prop:flot} shows that under some hypothesis on the Clarke gradient, such "ideal inverse flows" can be approximated as illustrated in \Cref{Fig:inverse}.
\end{remark}

\subsection{Relating normal cones to Clarke gradients of distance functions}
\label{sub:normal_clarke}

We prove several results on tangent cones of compact sets of $\R^d$ satisfying weak regularity assumptions, leading to \Cref{th:clarke_cone} which relates normal cones to the Clarke gradient of the distance function. These assumptions are verified by all \textit{complementary regular sets} defined in \Cref{sub:complementary_regular}, which is the class for which we will prove the Morse theorems.  \\

Recall from the paragraphs above \Cref{eq:zero_derivative_reach} that for any $A \subset \R^d$ and point $x \in A$, the tangent cone $\Tan(A,x)$ is closed, and that $u \in \Tan(A,x)$ is said to be \textit{represented} by a sequence $x_n$ in $A$ when this sequence never takes value $x$ and such that there exists a sequence $\eps_n$ of positive numbers converging to zero such that $x_n = x + \eps_n(u + o(1))$. Observe that $\frac{\eps_n}{\norme{x - x_n}}$ converges to $\norme{u}$. \\

Our first lemma relates the tangent cones of the boundary $\partial X$ to those of $X$ and $\complementaire{X}$.
\begin{lemma}[Tangent cone of the boundary]
\label{lem:partial_tangent}
Let $X \subset \R^d$. Then for every $x \in \partial X$,
\[ \Tan(\partial X,x) = \Tan(X,x) \cap \Tan(\complementaire{X},x).\]
\end{lemma}
\begin{proof}
The cone $\Tan(\partial X,x)$ being included in both $\Tan(X,x)$ and $\Tan(\complementaire{X},x)$, we have to prove that $\Tan(X,x) \cap \Tan(\complementaire{X},x)$ is included in $\Tan(\partial X,x)$.

Let $u \in \Tan(X,x) \cap \Tan(\complementaire{X},x)$ be of norm 1. Take a sequence $x_n$ (resp. a sequence $\complementaire{x_n}$) in $X$ (resp. $\complementaire{X}$) representing $u$, i.e., such that
\begin{align*}
    x_n & = x + \norme{x_n -x}(u + o(1) \,), \\
    \complementaire{x_n} & = x + \norme{\complementaire{x_n} -x}(u + o(1)\,).
\end{align*}
The segment $[x_n, \complementaire{x_n}]$ has to intersect $\partial X$, which means that there exists a $\lambda_n \in [0,1]$ such that $\partial x_n = \lambda_n x_n + (1 - \lambda_n)\complementaire{x_n}$ belongs to $\partial X$. This yields
\begin{align*}
\partial x_n - x = &  \left (\lambda_n \norme{x_n - x} + (1 - \lambda_n) \norme{\complementaire{x_n} -x} \, \right )(u + o(1)\,).
\end{align*}
Taking the norm of this equality yields 
\[\norme{\partial x_n - x} =  \left (\lambda_n \norme{x_n - x} + (1 - \lambda_n) \norme{\complementaire{x_n} -x} \, \right ) + o (\norme{\partial x_n - x}).\]
This quantity is strictly positive when $n$ is sufficiently large, and we have
\[\partial x_n - x =  \norme{\partial x_n -x}(u + o(1)\,), \]
meaning that $u$ is represented by the sequence $\partial x_n$, which lies in $\partial X$.
\end{proof}

The following lemma uses the zero-derivative characterization of Federer of tangent cones (\Cref{eq:zero_derivative}) to ultimately obtain in \Cref{cor:complementary} that tangent cones of complement sets are complement to tangent cones of original sets under assumptions on the reach of the complement set. This result holds assuming that the tangent cones have full dimension, or, equivalently, non-empty interior in $\R^d$. This condition is not satisfied by lower-dimensional sets; for instance when $X$ is a submanifold of $\R^d$ of positive codimension $\Tan(X,x)$ is a vector subspace of positive codimension in $\R^d$. Nevertheless, this condition will be satisfied by the complement of a \textit{complementary regular set}, which are the sets we study in \Cref{section:morse_complement}.
\begin{lemma}
Let $A \subset \R^d$ be a closed set with positive reach. Then for any $x$ in $\partial A$ and any non-zero $u \in \interior(\Tan(A,x))$, we have 
\begin{equation}
    \label{eq:distance_complement}
    \liminf_{t \to 0^+}\frac{1}{t} d_{\complementaire{\!A}}(x + tu) > 0.
\end{equation}
As a consequence, for any $A \subset \R^d$ with positive reach, we have:
\begin{equation}
\label{eq:tan_complement}
\interior(\Tan(A,x)) \cap \Tan(\complementaire{\!A},x) = \emptyset. 
\end{equation} 
\end{lemma}

\begin{proof}
First, observe that since $\Nor(A,x)$ is non-empty, 0 does not belong to the interior of its dual cone $\Tan(A,x)$.
We want to show that a vector of $\interior(\Tan(A,x))$ cannot satisfy the zero-derivative characterization of elements of $\Tan(\complementaire{A},x)$, given in \Cref{eq:zero_derivative_reach}.
Let $\lambda > 0$ be such that $\Tan(A,x)$ contains the closed ball of radius $\lambda$ centered at $u$. 
From the zero-derivative characterization of tangent cones for sets with positive reach we have for any $0 \leq \lambda' \leq \lambda$ and any unit vector $v$:
\begin{equation}
    \label{eq:federer_carac}
    d_A(x + t(u+\lambda'v)) =  o(t).
\end{equation}
For any $t \geq 0$, let $x(t)$ be a point in $\complementaire{\! A}$ realizing $d_{\complementaire{\! A}}(x + tu) = \norme{x+tu - x(t)} \eqqcolon d(t)$. 
For $t > 0$ small enough, the point $x + tu$ lies well within the reach of $A$, so that we can take $\delta > 0$ small enough and a unit $v(t) \in \Nor(A,x(t))$ (with the convention that $\Nor(A,x) = \R^d$ if $x \notin A$) 
such that ${d_{A}(x(t) + \delta' v(t)) \geq \delta'}$ for every $0 \leq \delta' \leq \delta$. Distance functions being $1$-Lipschitz, we have obtained
\begin{equation*}
    d_{A}(x + tu + t\delta v(t)) \geq d_{A}(x(t) + t\delta v(t)) - \norme{x + tu - x(t)} \geq t\delta - d(t).
\end{equation*}
Now consider a positive sequence $t_n$ converging to $0$. Extracting a subsequence we can assume that $v(t_n)$ converges to a unit vector $v$, leading to 
\begin{equation*}
    d_{A}(x + t_n(u + \delta v)) \geq t_n\delta - d(t_n) - o(t_n).
\end{equation*}
By \Cref{eq:federer_carac}, the right-hand side must be $o(t_n)$ as $n$ tends to $+\infty$, which shows that $\liminf_{n \to \infty} \frac{d(t_n)}{t_n} \geq \delta$. Since this applies to any sequence $t_n$ going to $0$, we have obtained the desired result \Cref{eq:distance_complement}.

Now to obtain \Cref{eq:tan_complement} let $y_n = x + a_n(u + o(1))$ be a sequence representing $u \in \interior(\Tan(A,x))$. By \Cref{eq:distance_complement} there is a $\delta > 0$ such that for $n$ large enough, $d_{\complementaire{\!A}}(x + a_n u) \geq \delta a_n$, implying that $d_{\complementaire{\!A}}(y_n) \geq a_n( \delta + o(1))$. The right hand-side quantity is strictly positive when $n$ is large enough. As such, $u$ cannot be represented by a sequence in $\complementaire{A}$.
\end{proof}

Taking $A = \complementaire{\!X}$ in the previous result yields the following.

\begin{corollary}[Complement of tangent cones are tangent cones of complements]
\label{cor:complementary}
Let $X \subset \R^d$ be a closed set such that $\complementaire{X}$ has positive reach. Then for any $x \in \partial (\complementaire{X})$, we have
\begin{equation*}
    \complementaire{\Tan(\complementaire{X}, x)} = \Tan(\complementaire{(\complementaire{X})},x).
\end{equation*}
In particular, 
if $X = {\complementaire{(\complementaire{X})}}$, or equivalently $X = \overline{\interior(X)}$,
 we have 
\begin{equation*}
    \complementaire{\Tan(\complementaire{X}, x)} = \Tan(X,x).
\end{equation*}
\end{corollary}

\begin{proof}
We write $Y = {\complementaire{(\complementaire{X})}}$ to ease notation.
Since $\complementaire{X} \cup Y = \R^d$, we have $\Tan(\complementaire{X},x) \cup \Tan(Y,x) = \R^d$.
Since tangent cones are closed, this implies
 $\complementaire{\Tan(\complementaire{X},x) \subset \Tan(Y,x)}$. The reverse inclusion is equivalent to $\Tan(Y,x) \cap \interior(\Tan(\complementaire{X},x))$ being empty, something we obtained in the previous lemma.
\end{proof}

Using this corollary, we are able to prove that the definition of $\Nor(X,x)$ as $- \Nor(\complementaire{X},x)$ is consistent when both $X, \complementaire{X}$ have positive reach. 

\begin{corollary}[Consistency of the normal bundle definition]
\label{cor:consistency}
When both $X \subset \R^d$, $\complementaire{X}$ have positive reach, for any $x \in \partial (\complementaire{X})$ we have
\begin{equation*}
\Nor(X,x) = - \Nor(\complementaire{X},x).
\end{equation*}
\end{corollary}
\begin{proof}
Recall from \Cref{eq:proj-char} that $\Nor(X,x), \Nor(\complementaire{X},x)$ each contain at least one non-zero vector.
By \Cref{cor:complementary}, we know that $\complementaire{\Tan(\complementaire{X},x)} \subset \Tan(X,x)$. Now for any non-zero $n \in \Nor(\complementaire{X},x)$, $\Tan(\complementaire{X},x) \subset \{ n \}^{\text{o}} = \{ z \in \R^d, \scal{z}{n} \leq 0 \}$. The complement set of $\{ n \}^{\text{o}}$ being $\{-n \}^{\text{o}}$, combining these last two inclusions yield $\{-n\}^{\text{o}} \subset \Tan(X,x)$. Taking in turn the dual cone gives $\Nor(X,x) \subset \Cone(-n)$. Since $\Nor(X,x)$ contains at least one unit vector, necessarily $\Nor(X,x) = \Cone(-n)$. Since this holds for any unit vector $n$ in $\Nor(\complementaire{X},x)$, these $n$ are all equal and $\Nor(X,x) = - \Nor(\complementaire{X},x)$.
\end{proof}

Using a characterization of the Clarke gradient similar (but distinct) to the zero-derivative characterization of tangent cones, we obtain the following result.
\begin{lemma}[Tangent cone stability under addition with $\clarke d_X(x)^{\text{o}}$]\label{lem:iterated}
Let $X \subset \R^d, x \in \partial X$  and $u \in \clarke d_X(x)^{\text{o}}$. 
Then for all $h \in \Tan(X,x)$, $u + h \in \Tan(X,x)$.
\end{lemma}
\begin{proof}
We use Clarke's \cite{Clarke1975GeneralizedGA} characterization of the dual cone to the Clarke gradient
of $d_X$ at a point $x \in \partial X$:
\begin{equation} \label{eq:dual d_x}
\clarke d_X(x)^{\text{o}} = \left \{ u \;  \Biggr |  \; \lim_{\substack{x_h \to x\\ x_h \in X}} \lim_{\delta \to 0^+} \frac{1}{\delta} d_X(x_h + \delta u) = 0 \right \}. 
\end{equation} 

Consider the following modulus of continuity:
\[ \omega_u(\eps, \lambda) \coloneqq \sup_{\substack{x_h \in X \\ \norme{x - x_h} \leq \eps}} \sup_{0 < \delta \leq \lambda} \frac{d_X(x_h + \delta u)}{\delta}.\]
When $u$ belongs to $\clarke d_X(x)^{\text{o}}$, by Clarke's characterization \ref{eq:dual d_x} we have 
\begin{equation*}
\lim_{\substack{\eps \to 0^+ \\ \lambda \to 0^+}}\omega_u(\eps, \lambda) = 0.
\end{equation*}

Let $h$ be a unit vector in  $\Tan(X,x)$ and let $x_i \to x$ be a sequence representing $h$. Put $\eps_i = \norme{x - x_i}$ and consider the sequence $x_i + \eps_i u$. Take $\xi_i$ in $\Gamma_X(x_i + \eps_i u)$, that is, a point in $X$ realizing the distance of $x_i + \eps_i u$ to $X$. By the definition of $\omega_u$, we have:
\[ \norme{\xi_i - x_i - \eps_i u} = d_X(x_i + \eps_i u) \leq \eps_i \omega_u(\eps_i, \eps_i). \]
Thus we can write
\begin{align*}
 \xi_i - x = & \, (\xi_i - x_i) + (x_i - x) \\
 = & \, \eps_i ( u + O(\omega_u(\eps_i, \eps_i)) + \eps_i(h + o(1))) \\
 = & \, \eps_i( u + h + o(1)),
 \end{align*}
which shows that $\xi_i$ is a sequence in $X$ representing $u+h$.
\end{proof}

Aggregating the previous lemma with \Cref{cor:complementary}, we obtain the following.

\begin{lemma}[Relationship between normal cones and Clarke gradients]\label{lem:conclusion_inclusion}
Let $X$ be a closed subset of $\R^d$ such that $\reach(\complementaire{X}) > 0$ and $\overline{\interior{(X)}} = X$, and let $x \in \partial X$. Then if $\Tan(\complementaire{X},x)$ has full dimension we have:
\[ \clarke d_X(x)^{\text{o}} \subset -\Tan(\complementaire{X},x).\]
\end{lemma}
\begin{proof}
Let $u \in \clarke d_X(x)^{\text{o}}$.
By \Cref{lem:iterated} we know that
\begin{equation*}
 \{ u \} + \Tan(X,x) \subset \Tan(X,x),
\end{equation*}
which is equivalent to 
\begin{equation}
\label{eq:inclusion_tan}
 \{ u \} + \left ( \R^d \setminus \Tan(X,x) \right ) =   \R^d \setminus \left ( \{ u \} +  \Tan(X,x) \right )   \supset \R^d \setminus \Tan(X,x).
\end{equation}
By \Cref{cor:complementary} we have $\complementaire{\Tan(X,x)} = \Tan(\complementaire{X},x)$. Thanks to the full dimensionality of $\Tan(X,x)$, and to the fact that tangent cones are closed, taking the closure of the inclusion of \Cref{eq:inclusion_tan} yields:
\begin{equation*} \{ u \} + \Tan(\complementaire{X},x) \supset \Tan(\complementaire{X},x).
\end{equation*}
Since $\Tan(\complementaire{X},x)$ contains zero, $u$ belongs to $-\Tan(\complementaire{X})$.
\end{proof}

We are now in position to relate normal cones of a set $X$ to the Clarke gradient of $d_X$ under weak regularity assumptions.

\begin{theorem}[Normal cones and the Clarke gradient of the distance function]
\label{th:clarke_cone}
Let $X \subset \R^d$ be such that $\reach(\complementaire{X}) > 0$, $\overline{\interior{(X)}} = X$ (or equivalently $\complementaire{(\complementaire{X})} = X$) and such that $\complementaire{X}$ is fully dimensional. Let $x \in  \partial \left ( \complementaire{X} \right )$.
Then the normal cone of $X$ at $x$ is determined by the Clarke gradient of $d_X$:  
\[ \Nor(X,x) = \Cone \left ( \clarke d_X(x) \right).  \]
\end{theorem}

\begin{proof}
Let $\mathrm{reach}(\complementaire{X}) > r > 0$.
First note that
by the characterization of the Clarke gradient in terms of closest points (\Cref{eq:clarke_distance}) 
as well as the characterization of normal cones in terms of  \Cref{eq:proj-char}, we have
\begin{align*}
\clarke d_{X^{-r}}(x) = & - \Conv  \left \{ \frac{x-z}{\norme{x-z}} \; \Biggr |\; z \in X^{-r} \text{ with } d_X^{-r}(x) = \norme{z-x} \right \}\\
= & - \Conv \{ u \in \Sphere^{d-1} \tq d_{\complementaire{X}}(x +ru) = r \})\\
= & - \Conv \left (\Nor(\complementaire{X},x) \cap \Sphere^{d-1}\right )\!.
\end{align*}
On the other hand, by definition the Clarke gradient of $d_{X^{-r}}$ at $x$ is determined locally by the gradients around $x$ in every direction: \\
\[ 
\clarke d_{X^{-r}}(x) = \Conv \left \{ \lim_{i \to \infty} \nabla d_{X^{-r}}(x_i) \; \Biggr | \;   x_i \underset{i \to \infty}{\to} x\right \}. \]
Now compare to the Clarke gradient of $d_X$ at $x$, for which the gradient contributing only come from directions outside of $X$ \cite{Clarke1975GeneralizedGA}: 
\[ \clarke d_X(x) = \Conv \left (\{ 0 \} \cup \left \{ \lim_{i \to \infty} \nabla d_X(x_i) \; 
\Biggr | \; x_i \underset{i \to \infty}{\to} x \text{ with } d_X(x_i) > 0 \;\right \} \right ).
\]
 Note that in both definitions we implicitly require $x_i$ to be points where $d_X$ is differentiable. On those points $x_i$ the gradients of $d_X$ and $d_{X^{-r}}$ coincide. 
Since cones contain $0$, the previous fact yields the inclusion
\begin{equation*}
\Cone \left (\clarke d_X(x) \right )  \subset \Cone \left ( \clarke d_{X^{-r}}(x) \right )= - \Nor(\complementaire{X},x).
\end{equation*}

The other inclusion $-\Nor(\complementaire{X},x) \subset \Cone \left( \clarke d_X(x) \right )$ 
is a consequence of 
\Cref{lem:conclusion_inclusion}: under our assumptions $\complementaire{X}$ is fully dimensional and $\complementaire{(\complementaire{X})} = X$ so that we have $(\clarke d_X(x))^{\text{o}} \subset - \Tan(\complementaire{X},x)$. The desired inclusion is obtained by applying the dual cone operation to the previous inclusion.
\end{proof}

\section{Morse theory for complementary regular sets} \label{section:morse_complement}

In this section, we use the previous tools and propositions to derive the two Morse theorems when $X$ is \textit{complementary regular} (cf. \Cref{sub:complementary_regular}) and $f$ is Morse (in the sense of \Cref{def:crit}). In this setting, the eroded sets $X^{-r}$ converge to $X$ in the Hausdorff sense when $r$ tends to 0, and they are also $C^{1,1}$-domains for all $ r \in (0, \reach(\complementaire{X}))$.\\

Our approach is as follows. Let $c \in \R$. Consider a family of functions $f_{r,c}$ converging to $f$  as $r$ tends to $0$, in a way we will later specify. When $r = 0$, our notations are consistent with $f_{0,c} = f$. Consider the sublevel sets: 
\[
\begin{array}{lcr}
 X_c = X \cap f^{-1}(-\infty, c]  & \text{ and }&
 X^{-r}_c \coloneqq X^{-r} \cap f_{r,c}^{-1}(-\infty,c].
\end{array}
\]
They are the zero sublevel sets of the following functions: 
\[
\begin{array}{l c r}
     \phi = d_X + \max(f-c,0) & \text{ and } & \phi_r = d_{X^{-r}} + \max(f_{r,c} - c, 0). 
\end{array}
\] 

\begin{itemize}
    \item In \Cref{sub:complementary_regular}, we define the regularity condition required on sets $X \subset \R^d$ for which we prove the Morse theorems. Such sets are called \emph{complementary regular}. We describe some of their properties and show that these sets are exactly sets of the form $Y^{\eps}$, where $Y$ is a compact subset of $\R^d$ and the map $\Delta \circ \clarke d_Y$ is strictly positive over $d_Y^{-1}( \{ \eps \})$, i.e, $X$ is the offset of a set at a regular value of its distance function.
    \item In \Cref{sub:smoothed}, we take $f_r$ as $f$ precomposed with a uniformly bounded smooth function. If $c$ is a regular value of $f_{|X}$, we prove that there exists a $K > 0$ such that 
    for any $r > 0$ small enough there exists a retraction of any tubular neighborhood $(X^{-r}_c)^K$ onto $X_c$. This is done by proving a technical lemma which ensures that we can build an approximate inverse flow of $\phi^c_r$ using \Cref{prop:flot} on the sets $(X^{-r}_c)^K$.  
    \item In \Cref{sub:isotopie} we study the case $r = 0$ and prove that for $\eps > 0$ small enough, the set $X_{c + a}$ is a deformation retract of $X_{c+b}$ for any $a \leq b \in [-\eps, \eps]$ when $c$ is a regular value, also using \Cref{prop:flot}. As a consequence, we obtain the constant homotopy lemma.
    \item In \Cref{sub:handle} we let $c$ be a critical value and assume there is only one critical point $x$ in $f^{-1}(c)$, which is non-degenerate. We show that for any $\eps >0$ small enough the change in topology between $X_{c+\eps}$ and $X_{c-\eps}$ is determined by the curvature of $X$ at the pair $\left (x, -\nabla f(x)/\norme{\nabla f(x)} \right )$ and the Hessian of $f_{|X}$ at $x$. We prove this by considering $f_{r,c}$ to be $f$ translated with magnitude $r$ in the direction $- \nabla f(x)$. 
    We extend this result to the case where the level set $f^{-1}(c)$ has a finite number of critical points by considering a modified, more involved $f_{r,c}$ which depends on the different critical points of $f^{-1}(c)$.
    \item In \Cref{sub:persistence}, we discuss the connections of our work to persistent homology as well as the Euler Characteristic Transform. We prove that the Euler Characteristic Transform is injective on the class of complementary regular sets.   
    \end{itemize}

\begin{remark}[Changes of topology in the sublevel set filtrations of $f_{|X}$ and $(-f)_{|\complementaire{\!X}}$]
\label{rk:complement}
Following Fu's definitions \cite{fuCurvatureMeasuresGeneralized1989}, when $X$ is a domain bounded by a $C^{1,1}$-hypersurface, the maps $f_{|X}$ and $(-f)_{|\complementaire{\!X}}$ share the same critical points. One is Morse if and only if the other is, and the change of topology around a critical point $x$ in one sublevel-set filtration determines the corresponding change in the other, with
\begin{equation}
    \label{eq:cellule_complementaire}
    \mathrm{dim cell}(f_{|X},x) = d - 1 - \mathrm{dim cell}((-f)_{|\complementaire{\!X}},x),
\end{equation}
where $\mathrm{dim cell}(f_{|X},x)$ (resp. $\mathrm{dim cell}((-f)_{|\complementaire{\!X}}, x)$) denotes the dimension of the cell glued around $x$ in the filtration $(f^{-1}(-\infty, t] \cap X)_{t \in \R}$ (resp. $(f^{-1}[-t , \infty) \cap (\complementaire{\!X}))_{t \in \R}$). This equality follows from the fact that around $x$, the second order approximation of $f$	 in $\partial X$ goes down in $\mathrm{dim cell}(f_{|X},x)$ directions (whence the gluing of a cell of dimension $\mathrm{dim cell}(f_{|X},x)$ around $x$ in the classical handle-attachment lemma) and up in $d - 1 - \mathrm{dim cell}(f_{|X},x)$ directions. It is the opposite for the second-order approximation of $(-f)_{|\complementaire{X}}$, yielding \Cref{eq:cellule_complementaire}. 

However, for general $X \subset \R^d$ the changes of topology of the sublevel set filtration of $f_{|X}$ can have a Morse behavior - that is, a change of the form of the gluing of a cell around $x$ - while not that of $(-f)_{|\complementaire{X}}$. This phenomenon is illustrated in \Cref{fig:contre_morse} below with a stratified set $X \subset \R^2$ such that $\overline{\interior(X)} = X$.

    \begin{figure}[!htb]
             \centering
        \includegraphics[width=\textwidth]{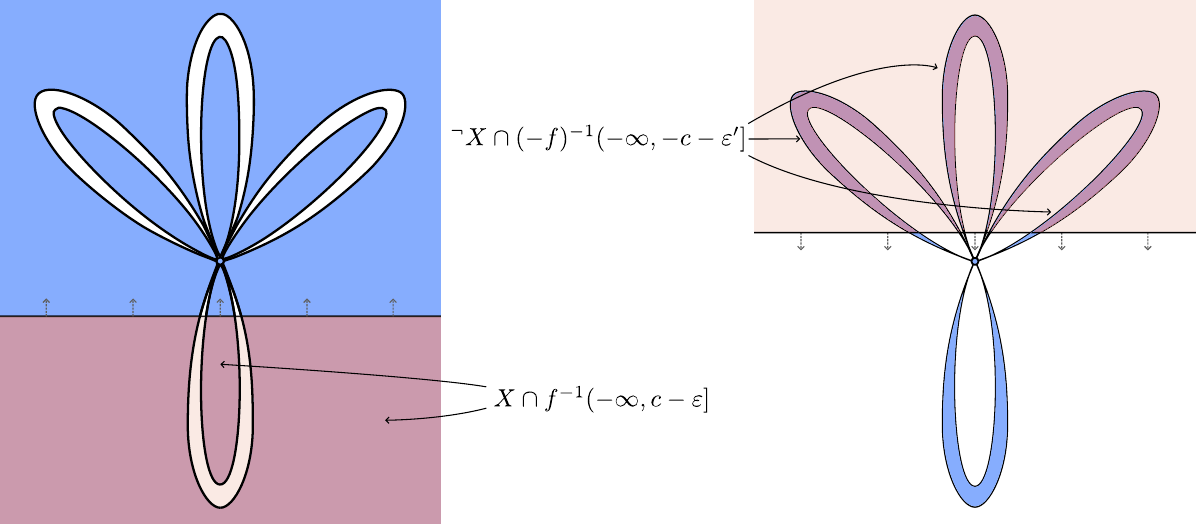}
        \caption{$\complementaire{X}$ is a bouquet of cuspy handles around a central point $x$ and $X$ is its complement set. Below (resp. above) the horizontal line of the left (resp. right) figure is a half-space of the form $f^{-1}(-\infty, c - \eps]$ (resp. $(-f)^{-1}(-\infty, -c - \eps']$) where $f$ is a vertical height function and $c = f(x)$.}\label{fig:contre_morse}
    \end{figure}
    
Indeed, when encountering the central point, the sublevel set filtration of $f_{|X}$ evolves by gluing one cell around the central point, going from 2 connected components (in purple in the left-hand side of \Cref{fig:contre_morse}) to 1. In that of $(-f)_{| \complementaire{X}}$ there are 5 changes in homology: going from 3 connected components (in purple in the right-hand side of \Cref{fig:contre_morse}) to 1 connected component and 3 cycles.
Note that the number of cell attachments describing the changes of homotopy type of the sublevel set filtration of $(-f)_{|\complementaire{X}}$ around $x$ can be arbitrarily high while maintaining a Morse behavior for $f_{|X}$ around $x$, as one can consider another set $X$ with more cuspy handles above the central point's height.
\end{remark}

\subsection{Complementary regular sets and their properties}
In this section, we define the class of \textit{complementary regular sets} which are the subsets of $\R^d$ for which we will prove the Morse theorems. We describe some of their properties and prove that they are exactly offsets of compact subsets of $\R^d$ at a regular value.  \\

\label{sub:complementary_regular}
\begin{definition}[Complementary regular sets]
\label{def:complementary_regular}
We say that a compact subset $X$ of $\R^d$ is a \emph{complementary regular set} when it satisfies the following three conditions: 
\begin{itemize}
    \item[$(A_1)$] $\overline{\interior(X)} = X$;
    \item[$(A_2)$] $\exists \, \mu \in (0,1]$ such that $\reach_{\mu}(X) > 0$;
    \item[$(A_3)$] $\reach(\complementaire{X}) > 0$.
\end{itemize}
\end{definition}

Intuitively, a complementary regular set is a compact domain (condition $(A_1)$), whose singularities are concave (condition $(A_3)$), but never cusps (condition $(A_2)$). More precisely, the following lemma shows that any number $\mu$ such that $\reach_{\mu}(X) > 0$ bounds the uniform thickness of the tangent cones of the complement set.

\begin{lemma}[Tangent cones of complements of complementary regular sets contain a ball]
Let $\mu \in (0,1]$ and let $X$ be complementary regular with $\reach_{\mu}(X) > 0$. Let $x \in \partial X$. 
Then $\Tan(\complementaire{X},x)$ contains a ball of radius $\mu$ centered at a unit vector.
\end{lemma}

\begin{proof}
By $(A_2)$, $\reach_{\mu}(X)$ is positive. Thanks to this fact, 
for each $r \in (0,\reach_{\mu}(X))$ there exists a point $x_r$ such that $d_X(x_r) = r$ and $\norme{x_r -x} \leq \frac{r}{\mu}$ (Theorem 3.2 in \cite{chazalShapeSmoothingUsing2007}).
Let $r_n$ be any sequence converging to $0$ and consider a sequence $x_n$ such that $\norme{x_n - x} \leq \frac{r_n}{\mu}$ and $d_X(x_n) = r_n$. 
Extracting a subsequence we can assume that $\frac{x_n -x}{\norme{x_n -x}}$ converges to a unit vector $u \in \Tan(\complementaire{X},x)$, i.e., we have
\begin{equation*}
x_n = x + \eps_n(u + o(1)).
\end{equation*} 
where $\eps_n = \norme{x_n - x} \to 0^+$. 
Now let $v \in \R^d$ be in the unit ball. 
Since $d_X(x_n) = r_n \geq \mu \eps_n$,
 the sequence $x_n + \mu \eps_n v$ belongs to $\complementaire{X}$ for every $n \in \N$. Moreover we have
\begin{equation*}
x_n + \eps_n \mu v = x + \eps_n(u + \mu v + o(1)), 
\end{equation*}
which implies that $u + \mu v$ belongs to $\Tan(\complementaire{X},x)$.
\end{proof}

\begin{corollary}[Normal cones of $\complementaire{X}$ are thin]
\label{cor:d_zero_mu}
Let $\mu \in (0,1]$ and let $X$ be complementary regular with $\reach_{\mu}(X) > 0$. Let $x \in \partial X$. Then 
\[ \Delta(\Conv(\Nor(\complementaire{X},x) \cap \Sphere^{d-1})) \geq \mu,\]
and for any unit vector $n \in \Nor(X,x) = - \Nor(\complementaire{X},x)$, there exists $u \in \clarke d_X(x)$ with $\norme{u} \geq \mu$ such that $n = \frac{u}{\norme{u}}$.

\end{corollary}
\begin{proof}
By the previous lemma, take a unit vector $u$ such that $B(u, \mu) \subset \Tan(\complementaire{X},x)$.
This yields the opposite inclusion on their dual cones $\Nor(\complementaire{X},x) \subset B(u,\mu)^{\text{o}}$. Take any unit vector $w \in B(u,\mu)^{\text{o}}$. For any $v \in \Sphere^{d-1}$, we have
\begin{equation}
\label{eq:half_space}
    0 \geq \scal{w}{u + \mu v} = \scal{u}{w} + \mu \scal{w}{v}.
\end{equation}
In particular, any such $w$ lies in the half space ${H^{-\mu}_{u} = \{ u' \in \R^d \tq \scal{u}{u'} \leq -\mu \} }$ which is a convex set with $\Delta(H^{-\mu}_u) \geq \mu$.

For the second point, recall that $\Nor(X,x) = \Cone(\clarke d_X(x))$ by \Cref{th:clarke_cone}
, and that by definition $\Nor(X,x) = - \Nor(\complementaire{X},x)$. Moreover, since $x \in \partial X$, we have by Corollary 2.5 in \cite{Clarke1975GeneralizedGA}
\begin{equation*}
\clarke d_X(x) = \Conv \left \{0, \lim \frac{x_i - e_i}{\norme{x_i - e_i}} \right \}, 
\end{equation*}
where we consider all sequences $x_i, e_i$ such that $x_i$ is not in $X$ and has closest
point $e_i$ in $X$, and $x_i \to x$ as $i \to + \infty$. Letting $V$ be the set of $\R^d$ consisting of all these limits, any vector in $\clarke d_X(x)$ can be expressed as a convex combination of elements of $V$ multiplied by a dilatation factor in $[0,1]$. Moreover, by \Cref{eq:half_space} and its immediate consequences,  $V \subset - H^{-\mu}_u$. Thus if $n$ is a unit vector in $\Cone( \clarke d_X(x))$, it can be written as $\normalise{u}$ where $u = t \sum \lambda_i v_i$ for some $t \in [0,1]$ and $u' = \sum_i \lambda_i v_i$ is some convex combination of elements of $V$ and thus belongs to $\clarke d_X(x)$, and also to $- H^{-\mu}_u$ so that it has norm at least $\mu$. The desired result follows from the equality $\normalise{u} = \normalise{u'}$.
\end{proof}

\begin{figure}[h!]
        \centering
        \includegraphics[scale = 2.5]{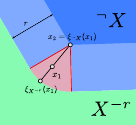}
        \caption{Illustration of $X^{-r}$ and $\complementaire{X}$ when $r \in (0, \reach(\complementaire{X}))$. The normal cone $\Nor(\complementaire{X},x_2)$ is depicted in red.}
        \label{fig:eroded}
    \end{figure}

Recall that the eroded set $X^{-r}$ is the set of points of $\R^d$ lying at distance at least $r$ to the complement $\complementaire{X}$, that is, the complement set to $(\complementaire{X})^r$, and that the projection map onto the closest point of $\complementaire{X}$ is denoted by $\xi_{\complementaire{X}}$.
Assuming that $\reach (\complementaire{X})$ is positive, the following properties on $\clarke d_{X^{-r}}$ and $\Nor(X^{-r})$ hold. They are illustrated in \Cref{fig:eroded}. 
\begin{lemma}[On eroded sets]
\label{lem:eroded}
Let $X$ be a closed subset of $\R^d$. Then for every $r \in (0, \reach(\complementaire{X}))$, 
$X^{-r}$ is a $C^{1,1}$ domain of $\R^d$ such that
\begin{itemize}
\item[(1)] $\reach(X^{-r}) \geq r$,
\item[(2)] $(X^{-r})^r = \complementaire{(\complementaire{X})}$,
\item[(3)] The unit normal bundle of $X^{-r}$ is the image of $\Nor(\complementaire{X})$ under the map $(x,n) \mapsto (x + rn, -n)$. In particular for any $x \in \partial ( X^{-r} )$,
 \[ \Nor(X^{-r},x) = \Cone(\xi_{\complementaire{X}}(x) - x). \]  
\item[(4.1)] For any $x$ such that $0 < d_{X^{-r}}(x) < r$, we have
\[ \clarke d_{X^{-r}}(x) =
 \left \{ \frac{ \xi_{\complementaire{X}}(x) - x }{\norme{\xi_{\complementaire{X}}(x) - x}} \right \}. \]
\item[(4.2)]
If $x$ belongs to $\partial (X^{-r})$, we have  
\[ \clarke d_{X^{-r}}(x) = [0, 1] \cdot \left \{ \frac{ \xi_{\complementaire{X}}(x) - x}{\norme{\xi_{\complementaire{X}}(x) - x}} \right \} = \Conv \left ( 0, \frac{ \xi_{\complementaire{X}}(x) - x}{\norme{\xi_{\complementaire{X}}(x) - x}} \right ),\]
which is a subset of $[0, 1] \cdot \clarke d_X(\xi_{\complementaire{X}}(x))$.
\item[(4.3)]
Finally, if $x \in \partial(\complementaire{X})$,
\begin{align*}
    \clarke d_{X^{-r}}(x) & = - \Conv(\Nor(\complementaire{X},x) \cap \Sphere^{d-1})  = \Conv(\Nor(X, x) \cap \Sphere^{d-1}) \\
    & = \Conv( \clarke d_{X}(\xi_{\complementaire{X}}(x)) \cap \Sphere^{d-1}).
\end{align*} 
\end{itemize}
\end{lemma}

\begin{proof}
First, observe that $X^{-r}$ is the upper-level set of the map $d_{\complementaire{X}}$ at a regular value, where it is differentiable with gradient equal to the opposite or the direction of the projection onto $\complementaire{X}$. At distance to $\complementaire{X}$ strictly less than $\reach(\complementaire{X})$, we know from Federer \cite{curvature:federer} that the projection map is Lipschitz, so that $X^{-r}$ is a $C^{1,1}$ domain.\\
 
For point $(1)$, let $x$ be a point lying at distance $0 < t < r$ to $X^{-r}$. Then it is within the reach of $\complementaire{X}$, so that it has a unique closest point $z$ in $\complementaire{X}$. By definition, the open ball $\interior{(B(z,r))}$ has empty intersection with $X^{-r}$. Since the ball $B(x,t)$ is included in $\interior{(B(z,r))}$ except for one point, which is $z + \frac{t}{r-t}(x - z)$, the closed ball $B(x,t)$ intersects $X^{-r}$ only at this point. \\

Since each point in $X^{-r}$ lies at distance greater or equal to $r$ to $\complementaire{X}$, points in $(X^{-r})^{r}$ are either outside of $\complementaire{X}$ or in its boundary, so that $(X^{-r})^r \subset \complementaire{(\complementaire{X})$. Conversely, any point lying at distance $0 < t < r$ to $\complementaire{X}$ lies at distance at most $r -t$ to $X^{-r}$, so that $\{ x \in \R^d, d_{\complementaire{X}}(x) > 0} \}$ is included in $(X^{-r})^r$. Taking the closure, we end up with $\complementaire{(\complementaire{X})} \subset \overline{X^{-r})^r} = (X^{-r})^r$ so that point $(2)$ is obtained.\\

Point $(3)$ follows from the definition $X^{-r} =\complementaire{((\complementaire{X})^r)}$: we have $\Nor(\complementaire{(X^r)}) = \{(x + rn, n), (x,n) \in \Nor(\complementaire{X}) \}$, and the normal unit bundle of the complement of $\complementaire{(X^r)}$ is obtained as its image by $(x,n) \mapsto (x, -n)$.\\

For points $(4.1), (4.2), (4.3)$, recall that for any closed set $Y \subset \R^d$ and $x$ outside of $Y$, when it exists the gradient $\nabla d_Y(x)$ is the opposite to the direction to the unique closest point of $x$ in $Y$, so that by definition the Clarke gradient $\clarke d_Y(x)$ is the convex hull of such directions.

For point $(4.1)$, if $d_X(x) \in (0, r)$, $x$ has one closest point in $X^{-r}$, and the proof of point $(1)$ shows that the closest point is in direction opposite to the projection onto the closest point in $\complementaire{X}$. 
Similarly, for $(4.2)$, let $x \in \partial (X^{-r})$ and look at the gradients of $d_{X^{-r}}$ around $x$. Outside of $X^{-r}$, we have gradients of $d_{X^{-r}}$ point to the closest point in $\complementaire{X}$ in a Lipschitz way, converging to the direction of $\xi_{\complementaire{X}}(x) -x$ as the neighborhood gets arbitrarily close to $x$. 
Inside $X^{-r}$, the gradients of $d_{X^{-r}}$ are 0.
Taking the convex hull yields the result.

Finally, point $(4.3)$ stems from the fact that closest points to a point $x \in \partial(\complementaire{X})$ in $X^{-r}$ are exactly $\{ x + rn, n \in \Nor(\complementaire{X},x) \}$.

\end{proof}

\begin{lemma}[Characterization of complementary regular sets]
Let $X$ be a compact subset of $\R^d$ and let $\mu \in (0,1]$. Then the 
simultaneous satisfaction of the three conditions
\begin{itemize}
    \item[$(A_1)$] $\overline{\interior(X)} = X$;
    \item[$(A'_2)$] $\reach_{\mu}(X) > 0$;
    \item[$(A_3)$] $\reach(\complementaire{X}) > 0$;
\end{itemize} is equivalent to the existence of $\eps, \delta > 0$ and of a compact subset $Y$ of $\R^d$ such that 
$X = Y^{\eps}$ with $\inf \{ \Delta(\clarke d_Y(x)) \tq d_Y(x) \in [\eps, \eps + \delta] \} \geq \mu$. The quantity $\reach_{\mu}(X)$ is the supremum of $\delta$ such that the previous inequality holds.
\end{lemma}
\begin{proof}
On the one hand, assume that the conditions
$(A_1), (A'_2), (A_3)$
are satisfied.
Then for any $0 < r < \reach(\complementaire{X})$ we have $(X^{-r})^r = X$ thanks to $(A_1)$ and \Cref{lem:eroded}.
Further assuming that $r < \reach_{\mu}(X)$, any such $X^{-r}$ will provide a suitable $Y$ with $\eps = r$. For any $x \in \R^d$ such that $d_X(x) > 0$ we have $d_{X^{-r}} = d_X + r$ on a neighborhood of $x$. Thus we have for any $\delta \in ( 0, \reach_{\mu}(X) )$:
\begin{equation*}
 \mu \leq \inf \{ \Delta(\clarke d_{X^{-r}}(x)) \tq d_{X^{-r}}(x) \in (r, r + \delta] \}. 
\end{equation*} It remains to bound $\Delta(\clarke d_{X^{-r}}(x))$ from below for points $x$ such that $d_{X^{-r}}(x) = r$. Those points form the set $\partial(\complementaire{X})$. For such $x$, we have $\clarke d_{X^{-r}}(x) = -\Conv(\Nor(\complementaire{X},x) \cap \Sphere^{d-1})$
by point (4.3) in \Cref{lem:eroded}.
\Cref{cor:d_zero_mu} then yields the desired bound ${\Delta (\clarke d_{X^{-r}}(x)) \geq \mu}$.

On the other hand, if there exist $\eps > 0, Y \subset \R^d$ such that $X = Y^{\eps}$ with $\inf \{ \Delta(\clarke d_Y(x)) \tq d_Y(x) \in [\eps, \eps + \delta] \} \geq \mu$, 
then since $X$ is an offset we have $\overline{\interior(X)} = X$ by \Cref{eq:offset_domain}, (condition $(A_1)$). Since $d_X = d_Y - \eps$ around any point at distance to $Y$ strictly greater than $\eps$, by definition of the $\mu$-reach we have $\reach_{\mu}(X) \geq \delta$, implying condition $(A'_2)$, and that $\reach_{\mu}(X)$ is equal to the supremum of such $\delta$. Finally,
 $Y$ is compact and $\eps$ is a positive regular value of $d_Y$, so that $\complementaire{(Y^{\eps})} = \complementaire{X}$ has positive reach by \Cref{th:reach_complements}.
\end{proof}

\begin{theorem}[Complementary regular sets are offsets at a regular value of the distance function]
\label{th:offsets}
A set is complementary regular if and only if it is the offset of a compact set at a regular value of its distance function.
\end{theorem}
\begin{proof}
This is a consequence of the previous lemma along with the lower semicontinuity of $\Delta \circ \clarke d_Y$: if $\reach_{\mu}(X) > 0$ and $X = Y^{\eps}$, there is a $\sigma > 0$ such that for every $x$ in $d_Y^{-1}[\eps - \sigma, \eps + \sigma]$, $\Delta(\clarke d_Y(x))$ is greater than $\frac{\mu}{2}$ and thus positive. From the set equality $d_Y^{-1}(\eps, \eps + \sigma] = d_X^{-1}(0, \sigma]$ and the fact that in this set $\clarke d_Y$ and $\clarke d_X$ coincide, we have the desired result.
\end{proof}

\begin{remark}[Approximating sets with their offsets]
In the previous theorem we established that complementary regular sets are exactly sets of the form $Y^\eps$ where $Y$ is a compact subset of $\R^d$ and $\eps > 0$ is a regular value of $d_Y$. One might wonder whether any compact $Y$ can be approximated up to arbitrary precision by sets of the form $Y^t$ which are complementary regular. This question reduces to the study of critical values of distance functions $d_Y$. In \cite{FuTube} it was proved that for ambient dimension $d=2,3$ the set of critical values of $d_Y$ is locally finite for any compact set $Y \subset \R^d$, but that this property is false for ambient dimension $4$ or above. Assuming that $Y$ is definable in an o-minimal structure (for instance, if $Y$ is semi-algebraic or subanalytic), this property holds nonetheless. In particular, the $\eps$-offset of a point cloud or a polyhedron is complementary regular for almost every positive $\eps$. Furthermore, the $\eps$-offset of a point cloud can be used to recover the homology of unknown underlying sets when $\eps$ is chosen among an interval prescribed by \cite{Sampling} under weak regularity assumptions. Our results show that such an offset can almost surely be used as an approximation on which the Morse theorems apply.
\end{remark}

The following proposition gives a characterization of complementary regular sets as regular sublevel sets of semiconcave functions.

\begin{proposition}[Characterization of complementary regular sets as sublevel sets of semiconcave functions]
A compact set $X \subset \R^d$ is complementary regular if and only if it is the sublevel set of a semiconcave function $\R^d \to \R$ at a regular value.
\end{proposition}

\begin{proof}
Assume that $X$ is complementary regular. By \Cref{th:offsets}, $X = Y^{\eps}$ for some positive $\eps$ and some compact set $Y \subset \R^d$, such that $\eps$ is a regular value of $d_Y$. Then $\eps^2$ is a regular value of $(d_Y)^2$, and $X = (d_Y^2)^{-1}(-\infty, \eps^2]$.

Conversely, let $f: \R^d \to \R$ be a semiconcave map with $c \in \R$ a regular value, and let $X = f^{-1}(-\infty, c]$. Since $c$ is a regular value, there is no local minimum of $f$ with value $c$ so that the closure of $\interior(X) = f^{-1}(-\infty, c)$ is $f^{-1}(-\infty, c] = X$, and assumption  $(A_1)$ is verified. Assumption $(A_3)$ is verified, since \Cref{th:reach_complements}, $\complementaire{X} = f^{-1}[c, + \infty)$ is a set with positive reach. There remains to prove $(A_2)$, i.e., that there exists $\mu > 0$ such that $\reach_{\mu}(X) > 0$. This is exactly the following \Cref{lem:semi-concave-sublevel}.
\end{proof}

\begin{lemma}[Regular sublevel sets of semiconcave maps have positive $\mu$-reach]
\label{lem:semi-concave-sublevel}
Let $f : \R^d \to \R$ be a semiconcave, $K$-Lipschitz map such that there exist $c \in \R$, $\rho > 0$ such that $\Delta( \clarke f(x) ) \geq K\rho$ on $f^{-1}(c)$ and such that $X \coloneqq f^{-1}(-\infty, c]$ is compact. Then  for every $\mu \in (0, \rho)$,
\begin{equation*}
\reach_{\mu}(X) > 0.
\end{equation*}
As a consequence, any compact sublevel set of a semiconcave map at a regular value has a positive $\mu$-reach for some $\mu \in (0,1)$.
\end{lemma}

\begin{proof}
Let $X = f^{-1}(-\infty, c]$ and assume that $\reach_s(X) = 0$. Then there exists a sequence $x_i$ with $d_X(x_i) > 0$ converging to zero such that $\Delta( \clarke d_X(x_i)) < s$. Consider $u_i$ the element of smallest norm in $\clarke d_X(x_i)$.
 We know that any convex combination of an arbitrary number of points in $\R^d$ is a convex combination of $(d+1)$ of these points (Carathéodory's Theorem). Using the description of $\clarke d_X(x_i)$ given in \Cref{eq:clarke_distance}, there exist closest points  $(x_{i,j})_{1 \leq j \leq d+1}$ of $x_i$ in $X$ such that
\[
u_i =  \sum_{j} \lambda_{i,j} \frac{x_i - x_{i,j}}{\norme{x_i - x_{i,j}}},
\] 
where $(\lambda_{i,j})_{j}$ are coefficients of a convex combination. Note that $x_{i}$ might have less than $(d+1)$ distinct closest points in $X$: we do not require the $x_{i,j}$ to be distinct. 
To ease notation, let $n_{i,j} = \frac{x_i - x_{i,j}}{\norme{x_i - x_{i,j}}}$.

On the one hand, observe that $\Nor(\complementaire{X},x_{i,j}) = - \Cone(n_{i,j})$. Indeed, the ball $B$ centered at $x_i$ of radius $\norme{x_i - x_{i,j}}$ is included in $\complementaire{X}$ so that the tangent cone $\Tan(\complementaire{X}, x_{i,j})$  contains the set $\Tan(B, x_{i,j}) = \{u \in \R^d, \scal{u}{n_{i,j}} \geq 0 \} = (-n_{i,j})^{\text{o}}$, which yields
\begin{equation*}
\Nor(\complementaire{X}, x_{i,j}) = \Tan(\complementaire{X}, x_{i,j})^{\text{o}} \subset (\{ -n_{i,j} \}^{\text{o}})^{\text{o}} = -\Cone(n_{i,j}).
\end{equation*}
From \Cref{eq:proj-char}, since $\complementaire{X}$ has positive reach and $x_{i,j}$ lies on the boundary of $\complementaire{X}$, the normal cone $\Nor(\complementaire{X},x_{i,j})$ must contain at least one non-zero vector, so that the  previous inclusion is necessarily an equality.

On the other hand, at any $x \in f^{-1}(c)$, the tangent cone $\Tan(\complementaire{X}, x)$ is included in $-(\clarke f(x))^{\text{o}}$. Indeed, since $f$ is semiconcave, by \Cref{eq:semi-concave} there exists a positive $C$ such that for any $u \in \clarke f(x)$, for any $y \in \R^d$ with $f(y) \geq c$, we have
\begin{equation*}
f(x) \leq f(y) \leq f(x) + \scal{u}{y -x} + \frac{C}{2} \norme{x-y}^2,
\end{equation*} 
so that $\scal{u}{\frac{y - x}{\norme{y - x}}} \geq -\frac{C}{2} \norme{x -y}$. Any unit vector $v \in \Tan(\complementaire{X},x)$ can be represented by a sequence $y_i$ converging to $x$ with $f(y_i) \geq c$, so $\scal{u}{v} \geq 0$ for all $u$. Since these inequalities hold for any $v$ in the tangent cone, we effectively have proved that $\Tan(\complementaire{X},x) \subset - \clarke f(x)^{\text{o}}$. 

Combining these two observations, we obtain that $n_{i,j} = \frac{v_{i,j}}{\norme{v_{i,j}}}$ for some $v_{i,j}$ in $\clarke f(x_{i,j})$ - note that any element of $\clarke f(x_{i,j})$ has norm at least $K\rho$ by assumption. 

Now up to extractions of subsequences, when $i$ goes to infinity, we can assume that each $n_{i,j}$ converges to a unit vector $n_j = \frac{v_j}{\norme{v_j}}$ where by upper semicontinuity of the Clarke gradient $v_j$ is some element of $\clarke f(x)$; and we can also assume that $\lambda_{i,j}$ converges to some $\lambda_j$ when $i \to + \infty$. The limit family $(\lambda_{j})_{j}$ is still that of a convex combination.
Since $u_i = \sum_{j} \lambda_{i,j} n_{i,j}$ has norm at most $s$, the vector 
\[
u = \sum_j \lambda_j \frac{v_j}{\norme{v_j}},
\]
which is the limit of $u_i$, has norm at most $s$.  Letting $\alpha = \sum_{j} \lambda_j \norme{v_j}^{-1}$, observe that $u/\alpha$ is a convex combination of the $v_j$ and as such belongs to $\clarke f(x)$, so that $\norme{u} \geq \alpha K \rho$. Since $f$ is $K$-Lipschitz, each $v_j$ has norm smaller or equal to $K$, so that $\alpha \geq 1/K$, whence the desired result $s \geq \rho$.
\end{proof}

\subsection{Building a homotopy equivalence between $X_c$ and its smooth surrogate}
\label{sub:smoothed}
For the remainder of this section, we let $X \subset \R^d$ be a complementary regular set and $f : U \to \R$ be a $C^2$ map defined on some open neighborhood $U$ of $X$. Following Fu's technique, we want to apply classical Morse theory to the eroded sets $X^{-r} = \{ x \in \R^d \tq d_{\complementaire{X}}(x) \geq r\}$, in a way that the topological events of the filtration 
$(X^{-r}_t)_{t \in \R} = \left ( X^{-r} \cap f_r^{-1}(-\infty, t] \right )_{t \in \R}$ converge to that of $(X_t)_{t \in \R} = (X \cap f^{-1}(-\infty, t])_{t \in \R}$ 
when $r \to 0$ and $f_r$ a function close to $f$, in a sense that we will later specify. To that end, we will consider the restriction to $X^{-r}$ of the sublevel sets filtration of $f$ slightly translated by a smooth function $\eta : \R^d \to \R$ with $\norme{\eta}_{\infty} \leq 1$.

\begin{definition}[Smooth surrogates for sublevel sets of $X$]
\label{def:surrogate}
 Let $c$ be a regular value of $f_{|X}$ and let $f_r$ be $f$ precomposed with the translation by $r \eta$:
\[ f_r : x \mapsto  f(x + r \eta(x)), \]
 where $\eta : \R^d \to \R^d$ is a smooth function with $\norme{\eta}_{\infty} \leq 1.$
We define the smooth surrogates for $X_c$ and non-negative $r$ as:
\[ X^{-r}_c \coloneqq X^{-r} \cap f_r^{-1}(-\infty,c],\]
and non-negative, locally Lipschitz functions
\[ 
\begin{array}{cccc}
   \phi^c \coloneqq d_X + \max(f-c, 0) & \quad \quad & \phi^c_r \coloneqq d_{X^{-r}} + \max(f_r -c, 0).
\end{array}
\]
verifying $X_c = (\phi^c)^{-1}(0)$ and $X^{-r}_c = (\phi^c_r)^{-1}(0)$. When the value of $c$ is clear from the context, we write $\phi_r$ instead of $\phi^c_r$ for notational convenience.
\end{definition}

When $r > 0$ is smaller than the reach of $X$, the set $X^{-r}$ is a $C^{1,1}$ domain. 
We aim at proving that $X^{-r}_c$ and $X_c$ have the same homotopy type when $c$ is a regular value and $r > 0$ is small enough, explaining the terminology \emph{smooth surrogate} for $X^{-r}_c$. We begin by showing that the following Hausdorff convergence holds.

\begin{lemma}[Hausdorff convergence of sublevel sets]
\label{lem:hausdorff_convergence}
Let $X \subset \R^d$ be a complementary regular set, let $f : U \to \R$ be a smooth map defined on an open neighborhood of $X$ and let $c$ be a regular value of $f_{|X}$. Then in the Hausdorff topology we have:
\begin{equation*}
\lim_{r \to 0} X^{-r}_c = X_c.
\end{equation*}  
\end{lemma}

\begin{proof}
Since $\norme{\eta} \leq 1$, we have $X^{-r}_c \subset (X_c)^r$ for any $r > 0$.
 Assume that there is no Hausdorff convergence of $X^{-r}_c$ to $X_c$ as $r$ tends to zero and that $c$ is a regular value. We want to obtain a contradiction.
If there is no Hausdorff convergence, there is a point 
$x \in X_c$
 and a positive real number $t$ such that $B(x,t)$ and  $X^{-r}_c$ have empty intersection for any $r > 0$ small enough. Since $c$ is a regular value, $\nabla f(x)$ is non-zero.
  Let $u$ be in $\Tan(X,x)$ and consider a sequence $x_n \in \interior(X)$ (possible since $\overline{\interior(X)} = X$) representing $u$, i.e., such that $x_n = x + \eps_n (u +o(1))$ with the sequence $\eps_n$ in $\R^+ \setminus{0}$ converging to 0.
 For $n$ large enough, $x_n$ lies in $B(x,t)$ and $x_n$ belongs to $X^{-r}$ for any $0 < r < d_{\complementaire{X}}(x_n)$. However, by assumption it does not belong to $X^{-r}_c$, so that
\[ f_r(x_n) > c.\]
Letting $r$ go to zero yields $c \leq f(x_n)$. Since $f(x) \leq c$, the first order expansion of $f$ at $x$ yields $\scal{\nabla f(x)}{u} \geq 0$.
This holds for any $u$ in $\Tan(X,x)$, which amounts to the following inclusion in a half-space:
\begin{equation}
	\label{eq:inclusion_cone}
  \Tan(X,x) \subset -\nabla f(x)^{\text{o}}.  
\end{equation}
By \Cref{cor:complementary},
taking the complements of both sides of \Cref{eq:inclusion_cone} yields the opposite inclusion $\nabla f(x)^{\text{o}} \subset \Tan(\complementaire{X},x)$. Applying the polar cone operation to both sides then yields
\begin{equation}
  \Nor(\complementaire{X},x) \subset \Cone( \nabla f(x)).  
\end{equation}
Recall from \Cref{eq:proj-char} that $\Nor(\complementaire{X},x)$ contains at least one non-zero vector, so that
$\Cone( \nabla f(x)) = \Nor(\complementaire{X},x) = - \Nor(X,x)$. This contradicts the fact that $c$ is a regular value of $f_{|X}$.
\end{proof}

When $c$ is a regular value, the following lemma gives a uniform lower bound on $\Delta \circ \clarke \phi_{r}$ over neighborhoods of $X^{-r}_c$ of fixed size when $r$ tends to $0$.

\begin{lemma}[Non-vanishing $\clarke \phi_{r}$ around a regular value]\label{lem:cas}
Let $c$ be a regular value of $f_{|X}$. Then there exists a positive constant $\alpha$ such that for any sequences of 
non-negative reals $(r_i)$,
positive reals $(K_i)$ such that $r_i, K_i \to 0^+$, and any sequence $(x_i)$ of points within $\phi^{-1}_{r_i}(0, K_i]$ for all $i \in \N$, we have:
\[  \liminf_{i \to \infty} \Delta(\clarke \phi_{r_i}(x_i)) \geq \alpha. \]
\end{lemma}

\begin{proof}
The map $\phi_{r_i} = d_{X^{-r_i}} + \max(0, f_{r_i} - c)$ is the sum of a Lipschitz function and the 
non-negative part of a $C^2$ function. 
Depending on the sign of $(f_{r_i}(x_i) - c)$, and on whether $x_i$ lies in the interior, the exterior or on the boundary of $X^{-r_i}$, we distinguish seven cases to
 provide a positive lower bound for $\Delta(\clarke \phi_{r_i}(x_i))$. In each of these cases  $d_{X^{-r_i}}$ and $\max(0, f_{r_i} - c)$ contribute in different ways to the Clarke gradient $\clarke \phi_{r_i}(x_i)$, which is reflected by the inclusion 
 \begin{equation*}
 \clarke \phi_{r_i}(x_i) \subset \clarke d_{X^{-r_i}}(x_i) + \clarke (\max(0, f_{r_i} - c)),
\end{equation*}
where the addition on subsets of $\R^d$ is the Minkowski sum $A + B = \{ a + b, a \in A, b \in B \}$. By extracting subsequences, we can assume that the sequence $(x_i)$ lies in one of these cases. They are illustrated in \Cref{fig:cas}, presented more systematically in the following table. Crossed-out cells correspond to the cases where $\Delta(\clarke \phi_{r_i}(x_i)) \leq 0$, which are irrelevant in our analysis.

\vspace{1em}
{
\begin{center} 
\renewcommand{\arraystretch}{1.5}
\setlength{\tabcolsep}{10pt}

\begin{tabular}{c||c|c|c|}
& $f_{r_i}(x_i) < c$
& $f_{r_i}(x_i) = c$
& $f_{r_i}(x_i) > c$
\\ \hline\hline

$x_i \in \interior(X^{-r_i})$
& $\times$
& $\times$
& Case 2
\\ \hline

$x_i \in \partial X^{-r_i}$
& $\times$
& $\times$
& Case 5
\\ \hline

$0 < d_{X^{-r_i}}(x_i) \le r_i$
& Case 7
& \multicolumn{2}{c|}{Case 6}
\\ \hline

$d_X(x_i) > r_i$
& Case 1
& Case 4
& Case 3
\\ \hline
\end{tabular}
\end{center}
}

\vspace{1em}

In fact, we will show that for any such sequence, we have:

\begin{equation*}
    \liminf_{i\to \infty} \Delta(\clarke \phi_{r_i}(x_i)) \geq \min(\mu, \sigma, \kappa) > 0,
\end{equation*}

where

 \begin{itemize}
     \item $\kappa \coloneqq \inf_{f^{-1}(c) \cap X} \norme{\nabla f}$
     is a positive quantity because $c$ is a regular value of $f_{|X}$.
     
     \item $\mu \leq \inf_{t \to 0} \{ \Delta(\clarke d_X(x)) \tq 0 < d_X(x) < t \}$ is positive by hypothesis.

     \item  $\sigma \coloneqq \inf_{x \in \partial X \cap f^{-1}(c)} \Delta(A_x)$
     where $x \mapsto A_x$ is the upper semicontinuous set-valued map defined by:
    \begin{equation*}
    A_x \coloneqq \left ( [0,1] \cdot \clarke d_X(x) + \{ \nabla f(x) \} \right )
     \cup \left ( \clarke d_X(x) + [0,1] \cdot \{ \nabla f(x) \}\right).
    \end{equation*}
     To see why $\sigma$ is positive when $c$ is a regular value of $f_{|X}$, 
    recall from \Cref{th:clarke_cone} the identity 
    \[ \Cone \left (\clarke d_X(x) \right) = \Nor(X,x)\]
    for any point $x \in \partial X$. 
   	The set $\partial X \cap f^{-1}(c)$ is compact, and the map $x \mapsto \Delta(A_x)$ is lower semicontinuous, meaning that its infimum over $\partial X \cap f^{-1}(c)$ is attained. If $\sigma$ were zero, there would be a point $x \in \partial X \cap f^{-1}(c)$ with $\Delta(A_x) = 0$. At such a point the direction of $-\nabla f(x)$ would belong to $\Nor(X,x)$, and thus $c$ would be a critical value of $f_{|X}$, a contradiction. Intuitively, $\sigma$ measures how transverse $X$ is to $f^{-1}(-\infty, c]$.
    
 \end{itemize}

    \begin{figure}[h!]
    \centering
    \includegraphics[width = 1.\textwidth]{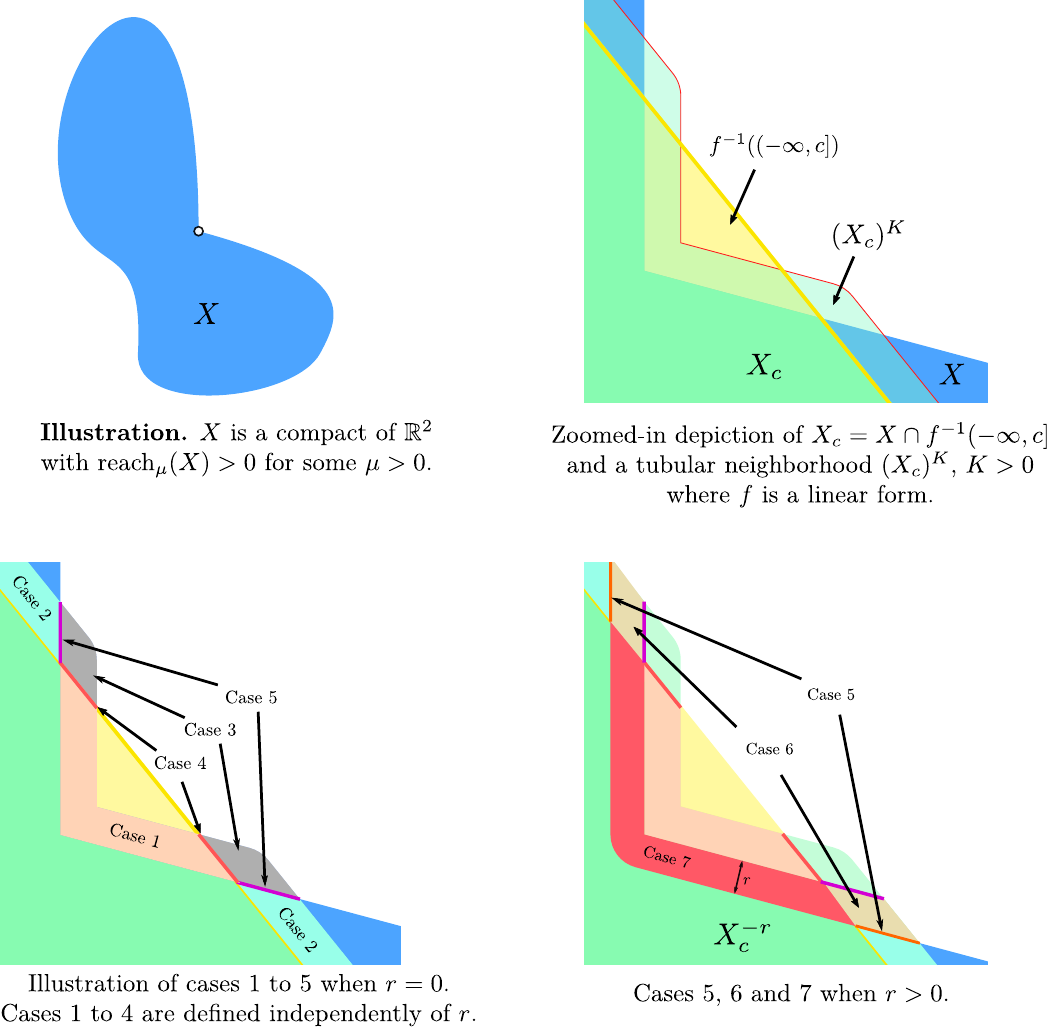}
    \caption{Illustration of the 7 cases of \Cref{lem:cas}.}
    \label{fig:cas}
\end{figure}

We will show that $\clarke \phi_{r_i}(x_i)$ either lies close to $\nabla f(x_i)$, $\clarke d_X(x_i)$ or close to being inside $A_{x_i}$, each being bounded away from zero respectively by the non-vanishing of $\kappa, \mu$ and $\sigma$. 

To ease notation, we write $\nu(x) \coloneqq \frac{x}{\norme{x}}$ and  $\norme{\nabla f_{r_i} - \nabla f}_{\infty, \overline{V}} \eqqcolon \eps_i$ the infinity norm of $\nabla f_{r_i} - \nabla f$ on the closure of some relatively compact, open set $V$ such that $X \subset V \subset \overline{V} \subset U = \mathrm{dom} f$. Note that $\eps_i = O(r_i)$.

\begin{itemize}
    \item[\emph{Case 1.}] $d_{X^{-r_i}}(x_i) > r_i$ and $f_{r_i}(x_i) < c$.\\
    
    Then in a neighborhood of $x_i$, $\phi_{r_i} = d_X$, so that $\clarke \phi_{r_i}(x_i) = \clarke d_X(x_i)$. Furthermore, we have $0 < d_X(x_i) < K_i + d_H(X^{-r_i}, X)$, and the right-hand side tends to $0$ as $i \to \infty$. By the $\mu$-reach hypothesis, and the lower semicontinuity of $y \mapsto \Delta(\clarke d_X(y))$, we thus have
    \begin{equation*}
        \liminf_{i \to \infty} \Delta(\clarke \phi_{r_i}(x_i)) \geq \inf_{\substack{t \to 0^+ \\ 0 < d_X(x) < t}} \Delta(\clarke d_X(x)) \geq \mu > 0.
    \end{equation*}
    
    \item[\emph{Case 2.}] $x_i \in \interior(X^{-r_i})$ and $f_{r_i}(x_i) > c$.\\
    
    Then $\phi_{r_i} = f_{r_i}$ in a neighborhood of $x_i$, so that $\clarke \phi_{r_i}(x_i) = \left \{ \nabla f_{r_i}(x_i) \right \}$ and $0 < f_{r_i}(x_i) - c \leq K_i$. As such, we have the inclusion $\clarke \phi_{r_i}(x_i) \subset \{ \nabla f(x_i) \}^{\eps_i}$, 
    and since $\eps_i$ goes to zero we have
    \begin{equation*}
    \liminf_{i \to \infty} \Delta(\clarke \phi_{r_i}(x_i)) \geq \inf_{\substack{t \to 0^+ \\c < f(x) \leq c + t}} \norme{\nabla f(x)} \geq \kappa > 0.
    \end{equation*}
    
    \item[\emph{Case 3.}] $d_{X^{-r_i}}(x_i) > r_i$ and $f_{r_i}(x_i) > c$.\\
    
    Then $\phi_{r_i} = d_{X^{-r_i}} + f_{r_i} = d_{X}(x_i) - r_i + f_{r_i}$ in a neighborhood of $x_i$, so that $\clarke \phi_{r_i}(x_i) = \clarke d_X(x_i) +  \{ \nabla f_{r_i}(x_i) \}$. 
Recall that $\norme{\nabla f -  \nabla f_{r_i}} \leq \eps_i$ so that 
 $\clarke \phi_{r_i}(x_i)$ lies within
 Hausdorff distance $\eps_i$ to  the set $\clarke d_X(x_i) + \{ \nabla f(x_i) \}$, which is itself a subset of $A_{x_i}$. Therefore we have $\clarke \phi_{r_i}(x_i) \subset (A_{x_i})^{\eps_i}$. Since $\eps_i$ tends to zero, we have
    \begin{equation*}
    \liminf_{i \to \infty} \Delta(\clarke \phi_{r_i}(x_i)) \geq \inf_{\substack{x \in \partial X \\ f(x) = c}} \Delta(A_x) \geq \sigma > 0.
    \end{equation*}
    \item[\emph{Case 4.}] $d_{X^{-r_i}}(x_i) > r_i$ and $f_{r_i}(x_i) = c$.\\
    
    Since $d_{X^{-r_i}}(x_i) > r_i$ we have $\clarke d_{X^{-r_i}}(x_i) = \clarke d_{X}(x_i)$ and  $d_X(x_i) > 0$, and furthermore $d_X(x_i) \to 0$ since $\displaystyle \lim_{r \to 0} X^{-r} = X$. Without loss of generality by extracting subsequence we can assume that $x_i$ converges to a point $x$, which has to belong to {$\partial X \cap f^{-1}(c)$}.
    
    Now $\nabla f_{r_i}(x_i)$ has to be non-zero for $i$ sufficiently large as $\eps_i = O(r_i)$ and
    \[ \liminf_{i \to \infty} \norme{\nabla f(x_i)} \geq \inf_{x \in X \cap f^{-1}(c)} \norme{\nabla f(x)} = \kappa, \] 
    which yields that the set $\{ y \tq f_{r_i}(y) \neq c \}$ has density one at $x_i$ by the local inverse function theorem. As the Clarke gradient can be computed in a set of density one at $x_i$ (see \Cref{eq:clarke_densite}), we have for any $x_i$ where $\nabla f_{r_i}(x_i) \neq 0$:
    \begin{align*}
  \clarke \phi_{r_i}(x_i) = \Conv \left \{ \; \lim_{n \to \infty} \nabla \phi_{r_i}(z_n) \tq z_n \to x_i,  f_{r_i}(z_n) \neq c \right \}. 
    \end{align*}
    We can decompose this set as  
    \[ \clarke \phi_{r_i}(x_i) = \Conv( A_+ \cup A_-), \]
    where
\[
    \begin{array}{rcl}
        A_+ \! & \coloneqq & \! \left \{ \; \lim_{n \to \infty} \nabla \phi_{r_i}(z_n) \; \Bigr | \;  z_n \to x_i, f_{r_i}(z_n) > c \right \} \,\\
        A_-  \! & \coloneqq & \! \left \{ \; \lim_{n \to \infty} \nabla \phi_{r_i}(z_n) \; \Bigr | \; z_n \to x_i, f_{r_i}(z_n) < c \right \} \! . 
    \end{array} \]
    Only $d_{X^{-r_i}}$ contributes to the gradients of $A_-$ whereas $f_{r_i}$ also contributes to $A_+$. Thus any point in $\Conv( A_+ \cup A_-)$ can be expressed as $u + \lambda \nabla f_{r_i}(x)$ where $u \in \clarke d_{X^{-r_i}}(x_i) = \clarke d_{X}(x_i)$ and $\lambda \in [0,1]$.   
	Since $\norme{\nabla f_{r_i} - \nabla f} \leq \eps_i$,
	we obtain $\Conv(A_+ \cup A_-) \subset (A_{x_i})^{\eps_i}$.
    Thanks to the lower semicontinuity of $y \mapsto \Delta(A_y)$, this finally yields, 
    \begin{equation*}
    \liminf_{i \to \infty} \Delta(\clarke \phi_{r_i}(x_i)) \geq \Delta(A_x) \geq \sigma > 0.
    \end{equation*}
    \item[\emph{Case 5.}] $x_i \in \partial X^{-r_i}$ and $f_{r_i}(x_i) > c$.\\
    
    If $r_i > 0$, 
	we have $\clarke d_{X^{-r_i}}(x_i) = [0, 1] \cdot \nu(\xi_{\complementaire{X}}(x_i) - x_i) \subset [0, 1] \cdot \clarke d_X(\xi_{\complementaire{X}}(x_i))$ by point 4.2 in \Cref{lem:eroded}.
Moreover, the map $f_{r_i}$ is smooth around $x_i$ so that its contribution to $\clarke \phi_{r_i}(x_i)$ is $\nabla f_{r_i}(x_i)$. This vector lies at distance at most $\eps_i$ to $\nabla f(x_i)$, which itself lies at distance at most $O(r_i)$ to $\nabla f(x_i)$, so that we obtain  
	\begin{align*}    
    \clarke \phi_{r_i}(x_i) \subset  [0, 1] \cdot \clarke d_X(\xi_{\complementaire{X}}(x_i)) + \{ \nabla f_{r_i}(x_i) \} \subset (A_{\xi_{\complementaire{X}}(x_i)})^{\eps_i + r_i}.
    \end{align*}   
    
    If $r_i = 0$, then $\clarke \phi_{r_i}(x_i) \subset [0,1] \cdot \clarke d_X(x_i) + \{ \nabla f(x_i) \}$ and we obtain similarly
    \[ 
    \clarke \phi_{r_i}(x_i) \subset (A_{x_i})^{\eps_i}.
    \]
    Either way, since $\eps_i, r_i$ go to zero 
    and since any cluster point of $x_i$ belongs to $\partial X \cap f^{-1}(c)$, we have
    \begin{equation*}
    \liminf_{i \to \infty} \Delta(\clarke \phi_{r_i}(x_i)) \geq \inf_{x \in \partial X \cap f^{-1}(c)} \Delta(A_x) \geq \sigma > 0.
    \end{equation*}

\end{itemize}
\vspace{0.5 cm}
    
    Now the remaining cases 6 and 7 fit inside sequences of points $(x,r)$ such that ${0 < d_{X^{-r}}(x) \leq r}$. Recall from \Cref{lem:eroded} that for such $x$
\begin{equation*}
\clarke d_{X^{-r}}(x) \subset \Conv( \clarke d_{X}(\xi_{\complementaire{X}}(x)) \cap \Sphere^{d-1}).
\end{equation*}

\begin{itemize}
    \item[\emph{Case 6.}] 
    $0 < d_{X^{-r_i}}(x_i) \leq r_i $ and $f_{r_i}(x_i) \geq c$\\
    
   We have 
    $\clarke \phi_{r_i}(x_i) \subset \clarke d_{X^{-r_i}}(x_i) + \clarke \max(0, f_{r_i} -c)(x_i)) \subset \Conv \left (\Nor(X,\xi_{\complementaire{X}}(x_i)) \right ) + [0,1] \cdot \nabla f_{r_i}(x_i)$. Now by compactness we can  assume that $x_i$ converges to some point $x$, which must belong to $\partial X \cap f^{-1}(c)$. This leads to
    \begin{equation*}
    \liminf_{i \to \infty} \Delta(\clarke \phi_{r_i}(x_i)) \geq \Delta(A_{x}) \geq \sigma > 0.
    \end{equation*}
    
    \item[\emph{Case 7.}]
    $0 < d_X^{-r}(x_i) \leq r_i $ and $f_{r_i}(x_i) < c$\\
    Then there is no contribution from $\nabla f_{r_i}(x_i)$, so that
    $\clarke \phi_{r_i}(x_i) = \clarke d_{X^{-r_i}}(x_i) \subset \Conv \left (\clarke d_X(\xi_{\complementaire{X}}(x_i)) \cap \Sphere^{d-1} \right )$.
     This yields
    \begin{equation*}
    \liminf_{i \to \infty} \Delta(\clarke \phi_{r_i}(x_i)) \geq \mu > 0.
    \end{equation*}
\end{itemize}
\vspace{1em}
\end{proof}

We are now able to build retractions in neighborhoods of fixed size of both $X_c$ and $X^{-r}_c$ when $r$ is small enough.
They will be used in the proof of \Cref{cor:homotopy_equivalence} to obtain that $X_c$ and its surrogate $X^{-r}_c$ have the same homotopy type.
\begin{lemma}[Deformation retractions around $X_c$ and $X^{-r}_c$]
\label{lem:clarke max}
Let $c$ be a regular value of $f_{|X}$. Using the notations of \Cref{def:surrogate}, there exist $K > 0, M \geq 0$ such that for all $r > 0$ small enough, there exist continuous piecewise-smooth flows \[ C: [0,1] \times \phi^{-1}(-\infty, K] \to \phi^{-1}(-\infty, K] \]
\[C^r: [0,1] \times \phi_r^{-1}(-\infty, K] \to \phi_r^{-1}(-\infty, K] \] such that:
\begin{itemize}
    \item $C(0, \cdot), C^r(0, \cdot)$ are the identity over their respective spaces of definition;
    \item $C(1, \phi^{-1}( -\infty, K]) = X_c$ and $C^r(1, \phi^{-1}_r (-\infty, K]) = X^{-r}_c$;
    \item For any $t \in [0,1]$, $C(t, \cdot)_{|X_c}$, $C^r(t, \cdot)_{|X^{-r}_c}$ are the identity over $X_c$ and $X^{-r}_c$;
    \item $(X_c)^{\frac{K}{M+1}} \subset \phi_r^{-1}(-\infty, K]$ and $(X_c^{-r})^{\frac{K}{M+1}} \subset \phi^{-1}(-\infty, K]$.
\end{itemize}

\end{lemma}

\begin{proof}

Note that $X_c = \phi^{-1}(0)$ and $X^{-r}_c = (\phi_r)^{-1}(0)$. We want to bound $\Delta \circ \clarke \phi_r$ and $\Delta \circ \clarke \phi$ from below 
among sets of the form $\phi^{-1}(0, K]$, $\phi_r^{-1}(0,K]$ for all $r$ small enough to apply \Cref{prop:flot}.

To obtain such a constant $K$, recall that $\phi_0 = \phi$ and let
\begin{equation*}
\omega(s,K) \coloneqq \inf_{\substack{r \in [0,s] \\ x \in \phi^{-1}_{r}(0,K] }} \Delta (\clarke \phi_r(x) ).    
\end{equation*}

\Cref{lem:cas} states that 
\begin{equation*}
    \liminf_{\substack{s \to 0^+ \\ K \to 0^+}} \omega(s,K) > 0.
\end{equation*}

When $K, s > 0$ are small enough, $\omega(s,K)$ is positive,  so that for all $r \in [0,s]$, $\Delta \circ \clarke \phi_r$  is uniformly bounded away from zero in $\phi^{-1}_{r}(0,K]$. Consider $C, C^r$ flows given by \Cref{prop:flot} for some $\eps \in (0, \omega(s,K) )$ respectively for $\phi$ and $\phi_r$, with $a= 0$ and $b =K$. 
Being deformation retractions, $C$ and $C^r$ satisfy the first three properties of our statement; there remains to prove the last one.

Let $M$ be a common Lipschitz constant for the $(\phi_r)_{r \in [0,s]}$
and put $\delta_r = d_H(X_c, X_c^{-r})$ for any $r$ in $(0,s)$. Recall from \Cref{lem:hausdorff_convergence} that $\delta_r$ tends to 0 as $r$ goes to zero.
Since $\phi_r$ is $M$-Lipschitz and takes value zero on $X_c^{-r}$ we have
\begin{equation*}
X_c \subset \phi^{-1}_r(-\infty, M \delta_r],
\end{equation*}
and furthermore for any $\alpha > 0$,
\begin{equation*}
(X_c)^{\alpha} \subset \phi^{-1}_r(-\infty, M (\delta_r + \alpha)].
\end{equation*}
For all $r$ small enough (more precisely for any $r$ such that $K/(M+1) + \delta_r \leq K/M$), we thus have
\begin{equation*}
(X_c)^{K/(M+1)} \subset \phi^{-1}_r(-\infty, K].
\end{equation*}
Since $\phi$ is also $M$-Lipschitz, the same reasoning holds for the inclusion $X^{-r}_c \subset \phi^{-1}(-\infty, K]$.
\end{proof}

\begin{corollary}[Homotopy equivalence]\label{cor:homotopy_equivalence}
Let $X \subset \R^d$ be a complementary regular set and $f : \R^d \to \R$ be a smooth function. Let $c$ be a regular value of $f_{|X}$, let $\eta$ be a smooth function with $\norme{\eta}_{\infty} \leq 1$ and let  $f_r : x \mapsto  f(x + r \eta(x))$.\\

Then for all $r > 0$ small enough, ${X^{-r}_c = X^{-r} \cap f_r^{-1}(-\infty, c]}$ and $X_c$ have the same homotopy type.
\end{corollary}

\begin{proof}
Let $r, K, M$ be as in the previous \Cref{lem:clarke max}. From \Cref{lem:hausdorff_convergence}, we can take $r$ small enough that $d_H(X^{-r}_{c}, X_c) < K/(M+1)$.
This implies that the homotopies $C$, $C^r$ are well-defined on $X^{-r}_c$, $X_c$ respectively.
 Letting $\psi \coloneqq C(1, \cdot)_{|X^{-r}_c} : X^{-r}_c \to X_c$ and $\psi^r \coloneqq C^r(1, \cdot)_{|X_c} : X_c \to X^{-r}_c$, their composition $\psi \circ \psi^r$ is homotopic to $\Id_{X_c}$ via the map
\[
\left \{
\begin{array}{ccl}
X_c \times [0,1] & \to & X_c \\
(x,t) & \mapsto & C(1,C(t,C^r(t,x))).
\end{array} \right.
\] 
Similarly, $\psi^r \circ \psi$ is homotopic to $\Id_{X^{-r}_c}$ via $(t,x) \mapsto C^r(1,C^r(t, C(t,x))).$
\end{proof}

\subsection{Constant homotopy type lemma}
\label{sub:isotopie}
In this section, we prove that the homotopy type of the sublevel sets of a smooth map restricted to a complementary regular set is constant between critical values.
\vspace{1em}

\begin{theorem}[Constant homotopy type between critical values]
\label{th:isotopy}
Let $X \subset \R^d$ be a complementary regular set. 
Let $f: U \to \R$ be a smooth map defined on an open set $U$ containing $X$, and let $a <b \in \R$ be such that $[a,b]$ contains only regular values of $f_{|X}$. 
Then $X_a$ is a deformation retract of $X_b$.
\end{theorem}

This theorem is a direct consequence of the compactness of $[a,b]$ and \Cref{lem:regular}, which we will prove using the following technical lemma.
\vspace{0.2cm}

\begin{lemma}[Regular values of the family $(\phi^{c})_{c \in \R}$ are open.]
\label{lem:technique}
Let $c$ be a regular value of $f_{|X}$ and let $\phi^{s} \coloneqq d_X + \max( f - s, 0)$ for any $s \in \R$.
Then we have:
\begin{equation*}
    \lim_{\substack{\delta \to 0^+\\ K \to 0^+}} \inf \left \{ \Delta(\clarke \phi^{c+a}(x)) \tq x \in (\phi^{c+a})^{-1}(0, K], a \in [-\delta, \delta] \right \} > 0.
\end{equation*}
\end{lemma}

\begin{proof}
We proceed by contradiction. Assuming the inequality is false, there exist two real sequences $a_i \to 0, K_i \to 0^+$, and $(x_i)_{i \in \N}$ a sequence in $\R^d$ such that:
\[ 
\begin{array}{ccc}
     \forall i \in \N, \, 0 < \phi^{c + a_i}(x_i) \leq K_i & \text{ and } & \displaystyle \lim_{i \to \infty} \Delta( \clarke \phi^{c + a_i}(x_i)) = 0.
\end{array}
\]

We use the same distinction of sequences of $\phi^{-1}_{c+a_i}(0,K_i]$ into cases as in the proof of \Cref{lem:cas}. We distinguish 5 cases to compute $\clarke \phi^{c+a_i}$, as depicted in the following table.
\vspace{1em}
{
\begin{center} 
\renewcommand{\arraystretch}{1.5}
\setlength{\tabcolsep}{10pt}

\begin{tabular}{c||c|c|c|}
& $f(x_i) < c + a_i$
& $f(x_i) = c + a_i$
& $f(x_i) > c + a_i$
\\ \hline\hline

$x_i \in \interior(X)$
& $\times$
& $\times$
& Case 2
\\ \hline

$x_i \in \partial X$
& $\times$
& $\times$
& Case 3
\\ \hline

$d_{X}(x_i) > 0$
& Case 1
& Case 5
& Case 4
\\ \hline
\end{tabular}
\end{center}
}

\vspace{1em}
\begin{itemize}
    \item[\emph{Case 1.}] $f(x_i) < c + a_i$ and $d_X(x_i) > 0$.
        \vspace{0.3em}
        
	In that case, in a neighborhood of $x_i$, we have $\phi^{c+a_i} = d_X$, so that $\clarke \phi^{c+a_i}(x_i) = \clarke d_X(x_i) $. Since $d_X(x_i)  \leq K_i \to 0$, we have:
\[ \liminf_{i \to \infty} \Delta( \clarke \phi^{c + a_i}(x_i)) \geq \liminf_{t \to 0^+} 
\{ \Delta(\clarke d_X(x)), 0 < d_X(x) \leq  t \} \geq \mu > 0. \]
    
    \item[\emph{Case 2.}] $x_i \in \interior(X)$ and $f(x_i) > c + a_i$.
    \vspace{0.3em}
    
	In that case, in a neighborhood of $x_i$, we have $\phi^{c+a_i} = f$, so that     
     $\clarke \phi^{c+a_i}(x_i) = \{ \nabla f(x_i) \}$. Moreover, $f(x_i)$ tends to $c$, so that \[ \liminf_{i \to \infty} \Delta( \clarke \phi^{c + a_i}(x_i)) \geq \inf_{x \in f^{-1}(c)} \norme{\nabla f(x)} \geq \sigma > 0. \]
    
    \item[\emph{Cases 3, 4}] 
    Respectively $f(x_i) > c + a_i, x_i \in \partial X$ and $f(x_i) > c + a_i, d_X(x_i) > 0$. 
        \vspace{0.3em}
    
    In these cases, in a neighborhood of $x_i$ we have $\phi^{c+a_i} = d_X + f$ and $f$ is smooth, yielding 
    \[\clarke \phi^{c+a_i}(x_i) = \clarke d_X(x_i) + \{ \nabla f(x_i) \} \subset A_{x_i},\]

    \item[\emph{Case 5}]
    $f(x_i) = c + a_i, d_X(x_i) > 0$.
    
     \vspace{0.3em}
     In this case we have the inclusion 
    \[ \clarke \phi^{c+a_i}(x_i) \subset \clarke d_X(x_i) +  [0,1] \cdot \{ \nabla f(x_i) \} \subset A_{x_i}.\]
    
    In cases 3, 4, 5, we have the inclusion $\clarke \phi^{c+a_i}(x_i) \subset A_{x_i}$. As in the proof of \Cref{lem:cas}, the map $y \mapsto \Delta(A_y)$ is 
    lower
    semicontinuous. Any cluster point of $x_i$ has to belong to $\partial X \cap f^{-1}(c)$. Since $c$ is a regular value, we have: 
    \[ \liminf_{i \to \infty} \Delta( \clarke \phi^{c + a_i}(x_i)) \geq \inf_{x \in \partial X \cap f^{-1}(c)} \Delta(A_x) \geq \kappa > 0. \]
\end{itemize}
\end{proof}

\begin{lemma}[Local deformation retractions]\label{lem:regular}
Let $X$ be complementary regular, $f : \R^d \to \R$ smooth and let $c$ be a regular value of $f_{|X}$. Then there is a positive $\delta$ such that for any ${-\delta \leq a \leq b \leq \delta}$, $X_{c + a}$ is a deformation retract of $X_{c+b}$.
\end{lemma}

\begin{proof}

By \Cref{lem:technique} there exist $\sigma, \delta, K > 0$ such that for every $a \in [-\delta, \delta]$ we have 
\begin{equation*}
    \Delta( \clarke \phi^{c+a}(x)) \geq \sigma \text{ for all } x \text { in } (\phi^{c+a})^{-1}(0,K].
\end{equation*}
Thus for every $\alpha \in [-\delta, \delta]$, by \Cref{prop:flot} (taking e.g. $\eps = \sigma/2$) there exists a deformation retraction $C_{c+\alpha}(\cdot, \cdot)$ of $(\phi^{c+\alpha})^{-1}[0,K]$ onto $(\phi^{c+\alpha})^{-1}(0)$, whose trajectories make $\phi^{c+\alpha}$ decrease.
By elementary computations for every $a < b \in [-\delta, \delta]:$
\begin{equation*}
\phi^{c+a}(X_{c+b}) \subset  [0, b-a] \subset [0,2\delta]
\end{equation*}
meaning that $X_{c+b} \subset (\phi^{c+a})^{-1}[0,K]$ when $2\delta \leq K$. Since $\phi^{c+b} \leq \phi^{c+a}$, and since trajectories of $C_{c+a}$ make $\phi^{c+a}$ decrease, so that for $\delta > 0$ small enough and any $t \in [0,1]$, we have
\begin{equation*}
C_{c+a}(t,X_{c+b}) \subset (\phi^{c+a})^{-1}[0,2\delta] \subset (\phi^{c+a})^{-1}[0, K] \subset (\phi^{c+b})^{-1}[0, K].
\end{equation*}
Consequently, the composition $C_{c+b}(s, C_{c+a}(t,x))$ is well-defined for any $t,s \in [0,1]$ and $x \in X_{c+b}$ and is continuous in each of these variables. Now letting $i$ be the inclusion $X_{c+a} \to X_{c+b}$ and $\psi \coloneqq C_{c+a}(1, \cdot) : X_{c+b} \to X_{c+a}$, one clearly has $\psi \circ i = \Id_{X_{c+a}}$. The map $i \circ \psi$ is homotopic to $\Id_{X_{c+b}}$ via the homotopy
\begin{equation*}
\left \{
\begin{array}{ccl}
X_{c+b} \times [0,1] & \to & X_{c+b}\\
(x,t) & \mapsto & C_{c+b}(1, C_{c+a}(t,x)).
\end{array} \right.
\end{equation*} 
\vspace{0.5em}
\end{proof}

\subsection{Handle attachment around critical values}
\label{sub:handle}

We now study the evolution of the topology of the sublevel set filtration of a Morse function and prove the handle attachment lemma.
We begin by showing that non-degenerate critical points are isolated. 

\begin{proposition}[The critical points of a Morse function are isolated]
\label{prop:isol2}
Let $X \subset \R^d$ be a set with positive reach or a complementary regular set and let $f : U \to \R$ be a smooth function
defined on an open set $U$ containing $X$. Then the set of critical points of $f_{|X}$ (as defined in \Cref{def:crit}) is closed, and in this set the non-degenerate critical points are isolated.
As a corollary, if $f_{|X}$ is Morse, it admits a finite number of critical points.
\end{proposition}

\begin{proof}
First we prove that the set of critical points is closed. Let $x_i$ be a sequence of critical points of $f_{|X}$ converging to a point $x \in X$. Either $\nabla f(x) = 0$, and $x$ is critical, or assume that $ \nabla f(x) \neq 0$. In this case, by continuity of $\nabla f$, for $i$ large enough $\nabla f(x_i) \neq 0$. Since the $x_i$ are critical this implies $(x_i, - \nabla f(x_i)/\norme{\nabla f(x_i)}) \in \Nor(X)$. Since this set is closed and $\nabla f/{\norme{\nabla f}}$ is continuous, we have $(x, -\nabla f(x)/\norme{\nabla f(x)}) \in \Nor(X)$, so that $- \nabla f(x) \in \Nor(X,x)$.

We now turn our attention to the isolatedness of non-degenerate critical points.
Let $x$ be a non-degenerate critical point
of $f_{|X}$. If $x \in \interior(X)$, the Hessian $H_x f_{|X} = H_x f$ is non-degenerate so that by a classical result of Morse theory, $\nabla f$ is non-zero in some neighborhood of $x$ included in $\interior(X)$.

Otherwise $x \in \partial X$. Since $\nabla f(x) \neq 0$ by non-degeneracy, $\nabla f$ does not vanish in a neighborhood of $x$ by continuity so that there is no critical point of $f_{|X}$ belonging to $\interior(X)$ in that neighborhood. It now remains to prove that $x$ is isolated among critical points in $\partial X$. Assuming it were not the case, suppose that there is a sequence $x_i$ in $\partial X$ of critical points of $f_{|X}$ all distinct from $x$ and converging to $x$. Since $\nabla f \neq 0$ in the compact set $f^{-1}(c) \cap X$, we can assume by continuity of $\nabla f$ that
\[ \exists K > 0, \forall i \in \N, \norme{\nabla f(x_i)} \geq K. \]
In particular, the map $\nabla f / \norme{\nabla f}$ is well-defined and smooth around any point or cluster point of the sequence $(x_i)$.

 Since $x_i$ is critical for every $i \in \N$, the unit vector $n_i \coloneqq - \nabla f(x_i)/\norme{\nabla f(x_i)}$ lies in $\Nor(X,x_i)$. The sequence $(x_i, n_i)$ lies in $\Nor(X)$ and converges to $(x,n)$ where $n \coloneqq - \nabla f(x)/\norme{\nabla f(x)}$. 
Extracting a subsequence, we can assume that ${(x_i -x, n_i - n)/\norme{(x_i -x, n_i - n)}}$ converges to some $(u,v) \in \Tan(\Nor(X), (x,n))$. Since $x$ is non-degenerate, $\Tan(\Nor(X),(x,n))$ is a vector space and both $u$ and $-u$ belong to $\pi_0(\Tan(\Nor(X),(x,n)) \subset \Tan(X,x)$, yielding $\scal{u}{n} = 0$. Moreover, the second fundamental form of $X$ at $(x,n)$ in the direction $u$ is given by: 
\begin{equation*}
    \sff_{x,n}(u,u) = \scal{u}{v}\!. 
\end{equation*} 
Since $\nu \coloneqq -\normalise{\nabla f}$ is Lipschitz around $x$, we have $\norme{n_i - n} = \norme{\nu(x_i) - \nu(x)} = O(\norme{x_i -x})$, ensuring $\norme{(x_i - x, n_i - n)} = O(\norme{x_i -x})$, so that $u$ is non-zero.
This implies that $\normalise{x_i -x}$ converges to $\normalise{u}$.
The first order expansion of $n_i = \nu(x_i)$ gives
\begin{equation*}
    n_i - n = \norme{x_i - x}D_x \nu \left ( \normalise{u} \right )+ o(\norme{x_i - x}).
\end{equation*}

If $D_x \nu (u) = 0$, we have $\norme{n_i - n} = o(\norme{x_i -x})$  meaning that $v = 0 = D_x \nu (u)$ and $\sff_x(u,u) = 0$. Otherwise, $\norme{n_i - n} \displaystyle \sim C \norme{x_i - x}$ for some $C > 0$. By elementary computations this also yields $D_x \nu (u) = v$ and we thus have in any case
\begin{equation*}
    D_x \nu (u) = v.
\end{equation*}
Now we can write the first order expansion of $\nabla f(x_i) + \norme{\nabla f(x_i)} n_i$,
which on the other hand by definition of $n_i$ is zero:
\begin{align*}
    0 = & \nabla f(x_i) + \norme{\nabla f(x_i)} n_i \\
    = & \norme{x_i - x} \left ( H_x f(u) + \norme{\nabla f(x)} D_x \nu (u) -  n \scal{n}{H_x f(u)} \right ) + o(\norme{x_i -x}). 
\end{align*}
 Taking the scalar product of this vector with $u$
and letting $x_i$ go to $x$
yields:
\begin{equation*}
    H_x f_{|X}(u,u) = H_x f(u,u) + \norme{\nabla f(x)} \scal{u}{v} = 0.
    \label{eq:degen}
\end{equation*}
This contradicts the non-degeneracy of $H_x f_{|X} $ in the direction $u$ which belongs to $\pi_0(\Tan(\Nor(X)), (x,n)) \setminus \{ 0 \}$.

\end{proof}

When $c$ is a critical value of $f_{|X}$ with only one corresponding critical point $x$, we choose the map $\eta : \R^d \to \R$ as follows.

\begin{definition}[Choice of surrogates when there is at most one critical point per level set]
\label{def:notation}
Let $X \subset \R^d$ be a complementary regular set and let $f : U \to \R$ be a smooth function
defined on an open set $U$ containing $X$. If $c \in \R$ is such that $f^{-1}(c)$ contains only one critical point $x$ of $f_{|X}$ which is non-degenerate, we put for any $r > 0$:

\[
     \gamma^c_r \coloneqq y \mapsto y - r \normalise{\nabla f(x)} \qquad \qquad \quad 
     f_{r,c} \coloneqq f \circ \gamma^c_r : y \mapsto f \left ( y - r \normalise{\nabla f(x)} \right ) \]
\vspace{0.1cm}

When the value $c$ is clear from the context, we write $\gamma_r$ and $f_r$ instead to ease notation.
\end{definition}

The following two lemmas focus on the properties of the critical points of $f_{r_{|X^{-r}}}$.

\begin{lemma}[Local correspondence between critical points of $f_{|X}$ and $f_{r|X^{-r}}$]\label{lem:isol_index}
Let $X$ be a complementary regular subset of $\R^d$. Assume $x \in \partial X$ is a non-degenerate critical point of $f_{|X}$ and let $\mathrm{ind}_x$ be the index of the Hessian of $f_{|X}$ at $x$. 
Then 
$x^r \coloneqq x + r\frac{\nabla f(x)}{\norme{\nabla f(x)}}$
is a critical point of $f_{r |X^{-r}}$ such that $f_r(x^r) = f(x)$ for all $0 < r < 
\reach(\complementaire{X})$. When $r$ is small enough, $x^r$ is a non-degenerate critical point of $f_{r | X^{-r}}$, whose Hessian at point $x_r$ has index 
\[ \mathrm{ind}^r_x \coloneqq \mathrm{ind}_x  + \text{number of infinite curvatures} \text{ of } \Nor(X) \text{ at } \left ( x , -\normalise{\nabla f(x)} \right ). \]
\end{lemma}

\begin{proof}
Let $n = -\normalise{ \nabla f(x)} \in \Nor(\complementaire{X},x)$. Keep in mind that $f_r(x) = f(x + r n)$.

The pair $(x,n) \in \Nor(X)$ is regular by non-degeneracy of $f$ at $x$. Denote by $(\kappa_i)_{1 \leq i \leq d-1 }$ the principal curvatures (defined in \Cref{prop:tangent_spaces}) of $X$ at $(x,n)$ sorted in ascending order and put $m \coloneqq \max \{ i \, | \,  \kappa_i < \infty \}$. From there we follow the reasoning of Fu \cite{fuCurvatureMeasuresGeneralized1989}. When $0 < r < \reach(\complementaire{X})$, $X^{-r}$ is a $C^{1,1}$-domain and the regularity of the pair $(x,n)$ in $\Nor(X)$ guarantees that the Gauss map  $x \in \partial (\complementaire{X^{-r}}) \mapsto n(x) \in \Sphere^{d-1}$ is differentiable at $x - rn$. Per \Cref{lem:eroded} we have the following linear correspondence between tangent spaces:
\[ \Tan(\Nor(\complementaire{(X^{-r}})), (x - rn, -n)) = \{ (\tau + r \sigma, \sigma) \tq (\tau, \sigma) \in \Tan(\Nor(\complementaire{X}),(x,-n)) \}. \]
Since $\Nor(X^{-r}) = \{ (z,-\nu) \,| \, (z,\nu) \in \Nor(\complementaire{(X^{-r})}\}$ we have:
\begin{equation}
\label{eq:tan_eroded}
\begin{aligned}
 \pi_0(\Tan(\Nor(X^{-r}), (x - rn,n))) = & \{ \tau + r \sigma \tq (\tau, \sigma) \in \Tan(\Nor(\complementaire{X}),(x,-n)) \} \\
 = & \{ \tau - r \sigma \tq (\tau, \sigma) \in \Tan(\Nor(X),(x,n)) \}.
 \end{aligned}
\end{equation}
From \Cref{prop:tangent_spaces} there exists a basis of the space $\Tan(\Nor(X), (x,n))$ of the form 
\[ \left (\frac{b_i}{\sqrt{1 + \kappa^2_i}}, \frac{\kappa_i b_i}{\sqrt{ 1 + \kappa^2_i}} \right )_{1 \leq i \leq d-1} \] 
where the family $(b_i)$ forms an orthonormal basis of the subspace of $\R^d$ orthogonal to $n$ and the principal curvatures $\kappa_i$ belong to $(-\infty, + \infty]$; we choose this basis to be indexed in non-decreasing order with respect to the $\kappa_i$.
 Since the space
 $\pi_0(\Tan(\Nor(X^{-r}), (x - rn,n)))$ 
is identifiable with the classical tangent space of the $C^{1,1}$ domain $X^{-r}$, it has dimension $(d-1)$ for any $r \in (0, \reach(\complementaire{X}) )$, so that by the explicit description of \Cref{eq:tan_eroded} any finite principal curvature $\kappa_i$ must be at most $\reach(\complementaire{X})^{-1}$. 

Now compute an explicit expression of the Hessian: for any $\tau - r \sigma, \tau' - r \sigma'$ in $\pi_0(\Tan(\Nor(X^{-r}), (x - rn,n)))$:
\begin{align*}
    H_{x-rn} f_{r|X^{-r}}&  (\tau- r \sigma, \tau'- r \sigma')\\
 = & \, H_{x-rn} f_r (\tau- r \sigma,\tau' - r\sigma') + \norme{\nabla f_r(x^r)} \sff_{x-rn} (\tau- r\sigma, \tau' - r\sigma') \\ 
    = & \, H_{x} f(\tau- r \sigma,\tau' - r\sigma') + \norme{\nabla f(x)}\scal{\tau- r\sigma}{\sigma'}. \\ 
\end{align*}

A similar computation was already carried out in 4.6 of \cite{fuCurvatureMeasuresGeneralized1989}.
The above explicit formulas invite us to consider the union of the two families $((1-r \kappa_i)b_i)_{1 \leq i \leq m}$, $(r b_i)_{m+1 \leq i \leq d-1}$ as a convenient orthonormal basis of $\pi_0(\Tan(\Nor(X^{-r}), (x - rn,n)))$.
The $(i,j)$ coefficient of the matrix $M$ of the bilinear form $H_{x-rn} f_{r|X^{-r}}$ at $1 \leq i,j \leq m$ in this basis is
\[
\mathmakebox[1.1\linewidth][c]{
    \displaystyle 
\left \{
\begin{array}{l r}
 (1 - r \kappa_i)(1 - r \kappa_j)H_x f(b_i, b_j) + \delta_{i,j}\norme{\nabla f(x)}( 1 - r \kappa_i) &  1 \leq i, j \leq m \\
 r(1 - r \kappa_i) H_x f(b_i,b_j) & 1 \leq i \leq m, m+1 \leq j \leq d-1 \\
 r^2H_x(b_i,b_j) - r \norme{\nabla f(x)}\delta_{i,j} & m+1 \leq i,j \leq d-1
\end{array} \right.
}
\]
where $\delta_{i,j} = 1$ if $i=j$ and $0$ otherwise.

Identifying coefficients in front of the $r$-monomials, there are square matrices $A_1, A_2, A_3$ of size $m$, a square matrix $B$ of size $(d-m-1)$ and two rectangular matrices $C_1, C_2$ such that the Hessian $H_{x-rn} f_{r | X^{-r}}$ has the form 
\begin{equation}
M_1 = 
\begin{pmatrix}
A_1 + rA_2 + r^2A_3 & r C_1 + r^2 C_2 \\
r C_1^{T} + r^2 C_2^{T} & -r \norme{\nabla f(x)}Id+ r^2 B\\
\end{pmatrix}
\end{equation}
where $A_1$ is the matrix of $H_x f_{|X}$ in the basis $(b_i)_{1 \leq i \leq m}$. It is the same computation as in \cite{fuCurvatureMeasuresGeneralized1989} except that we end up with a minus sign in front of the identity in the lower right corner.
To compute the index of $M_1$, observe that for $r$ small enough and $0 \leq t \leq 1$ the upper left corner of the matrix
\begin{align*}
M_t \coloneqq &
\begin{pmatrix}
A_1 + t (rA_2 + r^2A_3) & t(r C_1 ^{T} + r^2 C_2 ^{T}) \\
t(r C_1 + r^2 C_2) & -r \norme{\nabla f(x)}Id+ t r^2 B\\
\end{pmatrix} 
\end{align*}
is invertible. From a classical formula in linear algebra for block determinants, a matrix of the form
\[ 
\begin{pmatrix}
A & C^T \\
C & B
\end{pmatrix}
\]
has determinant $\det(A)\det(B - CA^{-1}C^{T})$ when $A$ is invertible. Applied to the matrix $M_t$, this formula reads
\vspace{0.1em}
\begin{equation*}
\mathmakebox[1.1\linewidth][c]{
 \det(M_t) =  r^{d-m-1}\det(A_1 + tr A(r)) \det(-\norme{\nabla f(x)}I_d + t [rB + r C(r)(A_1 + trA(r)^{-1}) C(r)^{T}]),
 }
\end{equation*}
\vspace{0.1em}
where $A(r), C(r)$ are matrices polynomial in $r$. For $r$ small enough the determinant of $M_t$ is non-zero for any $t \in [0,1]$, so that there is a continuous path of matrices with non-vanishing determinant between $M = M_1$ and $M_0$. Consequently $M_1$ has the same index as $M_0$, which is
that of $A_1$ plus the dimension of the identity matrix in the lower right corner. 
\end{proof}

\begin{lemma}[Critical points of $f_{r_{|X^{-r}}}$ when $r$ is small enough]
\label{lem:isol3}
Let $X \subset \R^d$ be a complementary regular set. 
Let $f : U \to \R$ be a smooth function defined on an open set $U$ containing $X$ and assume that there is exactly one critical point $x \in \partial X$ of $f_{|X}$ at a level $f^{-1}(c)$, and that it is non-degenerate.
 Then there exist $\delta, R  > 0$ such that every $r \in [0, R]$, $c$ is the only critical value of $f_{r|X^{-r}}$ in $[c - \delta,c + \delta]$ and such that $x^r = x + r \normalise{\nabla f(x)}$ is the only critical point of $f_{r|X^{-r}}$ with value $c$.
\end{lemma}

\begin{proof}
Direct computations show that $f_r(x^r) = f(x)$ and $\nabla f_r(x^r) = \nabla f(x)$. For $r > 0$ small enough, we have {$\Nor(X^{-r},x^r) = - \Cone ( \nabla f(x) )$} (exactly as $x_1$ in \Cref{fig:motivation_critique}) and $x^r$ is a critical point of $f_{r|X^{-r}}$. Assuming the claim of \Cref{lem:isol3} is false, there are sequences $\delta_i >0, r_i \geq 0$ converging to $0$, and $y_i$ a sequence in $X$ such that:
\vspace{0.1cm}
\begin{center}
\begin{tabular}{l c r}
\textbullet\ $y_i \neq x^{r_i}$ \quad
&
\quad \textbullet\ $c-\delta_i
\leq f_{r_i}(y_i)
\leq c+\delta_i$
&
\textbullet\ $y_i$ is critical for $f_{r|X^{-r_i}}$.

\end{tabular}
\end{center}
\vspace{0.1cm}
Since there is only one critical point in $f^{-1}(c)$, which is non-degenerate, $\nabla f$ is non-zero on $X \cap f^{-1}(c)$, and by continuity it is also non-zero on some set of the form $X \cap f^{-1}(c - \eps, c + \eps)$ where $\eps > 0$ is small enough. We can thus assume that $\nabla f_{r_i}(y_i)$ is non-zero for $i$ large enough. Up to the extraction of a subsequence we may assume that $\nabla f_{r_i}(y_i)$ is non-zero for all $i$. Since $y_i$ is critical, by \Cref{lem:eroded} we have $y_i \in \partial (X^{-r_i})$, so that $y_i = x_i - r n_i$ where $x_i = \xi_{\complementaire{X}}(y_i)$ and
\[ n_i = -\normalise{\nabla f_{r_i}(y_i)} \in \Nor(X^{-r}, y_i) \subset \Nor(X, x_i). \]

Now let $\overline{x}$ be a cluster point of the sequence $y_i$ and assume after extraction of a subsequence that $y_i$ converges to $\overline{x}$, which has to belong to $ \partial X$. By continuity of the projection map, this implies that $x_i$ converges to $\xi_{\complementaire{X}}(\overline{x}) = \overline{x}$. Observe that  $n_i = -\normalise{\nabla f_{r_i}(y_i)}$ converges to $n = -\normalise{\nabla f(\overline{x})}$, and that $f(y_i) \to f(\overline{x}) = c$.
Since $X$ is complementary regular, let $\mu > 0$ be such that $\reach_{\mu}(X) > 0$.
By the second property of \Cref{cor:d_zero_mu}, $n_i \in \Nor(X,x)$ can be expressed as $\normalise{u_i}$ for some $u_i$ in $\clarke d_X(x_i)$ with $\norme{u_i} \geq \mu$. By the upper semicontinuity of the Clarke gradient, after passing to a subsequence $u_i$ converges to some $u \in \clarke d_X(\overline{x})$, whose norm is at least $\mu$. Thus $n$ converges to $\normalise{u}$, which belongs to $\Cone( \clarke d_X( \overline{x})) = \Nor(X,\overline{x})$, so that
\begin{equation}
\label{eq:crit_lim}
- \normalise{\nabla f(\overline{x})} \in \Nor(X,\overline{x}),
\end{equation}
that is, $\overline{x}$ is a critical point of value $c$. By the uniqueness assumption on critical points in $f^{-1}(c)$, $\overline{x} = x $; so that $y_i$ converges to $x$.

If we assume that $x_i = x$ for all $i$, then $y_i = x + r_i n_i$. Since $y_i \neq x^{r_i}$, $n_i$ and $n = - \nabla f(x)/\norme{\nabla f(x)}$ are not equal, and we can also assume that $\normalise{n_i - n}$ converges to some unit vector $v' \in \R^d$ by extracting a subsequence.
Then we would have 
\begin{align*}
n_i - n = \, & - \normalise{\nabla f(x + r_i(n_i - n))} + \normalise{\nabla f(x)} \\
 =\,  & -r_i \norme{n_i - n} (D_x \nu)(v) + o(r_i\norme{n_i - n}) \\
 =\, & o(\norme{n_i -n}),
\end{align*}  
which is absurd. We can thus assume without loss of generality that $x_i$ is different from $x$ for all $i \in \N$. Arguing exactly as in the proof of \Cref{prop:isol2}, we can extract a subsequence such that the sequence $\normalise{(x_i - x, n_i - n)}$ converges to $(u,v) \in \Tan(\Nor(X),(x,n))$. The same computations yield that the restricted Hessian $H_x f_{|X} = H_x f + \norme{\nabla f(x)} \sff_{x,n} $ is degenerate in the direction $u \in
\pi_0( \Tan(\Nor(X)),(x,n)) \setminus \{ 0 \} $.
\end{proof}

Thanks to the two previous lemmas and the homotopy equivalence between the surrogate sublevel sets $X^{-r}_t$ and $X_t$ of \Cref{cor:homotopy_equivalence}, we prove the handle-attachment lemma in the case there is exactly one critical point at the critical value $c$, which is non-degenerate.

\begin{theorem}
[Handle attachment around unique critical values]
\label{th:cellule}
Let $X \subset \R^d$ be complementary regular and $f : U \to \R$ be a smooth function
defined on an open set $U$ containing $X$. Assume $f_{|X}$ has only one critical point $x$ in $f^{-1}(c)$ which is non-degenerate. Then for any $\varepsilon > 0$ small enough, $X_{c+\varepsilon}$ has the homotopy type of $X_{c-\varepsilon}$ with a $\lambda_x$-cell attached, where 
\begin{align*}
 \lambda_x \coloneqq & \text{ index } \text{ of the Hessian of }  f_{|X} \text{ at } x \\ + &
\text{ number of infinite curvatures of } \Nor(X) \text{ at } \left (x,-\normalise{\nabla f(x)} \right )
\end{align*}
when $x \in \partial X$, and $
 \lambda_x \coloneqq \text{ index } \text{ of the Hessian of }  f \text{ at } x$
 when $x \in \interior{(X)}$.
\end{theorem}

\begin{proof}
 When $x \in \interior(X)$, the classical handle attachment lemma for smooth maps applies. We turn our attention to $x \in \partial X$.
By \Cref{lem:isol3}, when $\eps, r > 0$ are small enough, there is only one critical point $x_r$ in $f_{r|X^{-r}}^{-1}( c - \eps, c + \eps)$. $X^{-r}$ is the upper-level set of the $C^{1,1}$ map $d_{\complementaire{X}}$ at a regular value, so that classical $C^{1,1}$ Morse theory (Theorem 3.9 in \cite{fuCurvatureMeasuresGeneralized1989}) applies:
 $X^{-r}_{c+\eps}$ has the homotopy type of $X^{-r}_{c -\eps}$ with a cell attached around $x^r$ should the Hessian of $f_{r|X^{-r}}$ be non-degenerate at $x^r$; in this case the dimension of the cell is the index of the Hessian. Thanks to \Cref{lem:isol_index}, this Hessian is indeed non-degenerate and its index is $\lambda_x$ for all $r > 0$ small enough. Now by \Cref{cor:homotopy_equivalence}, when $r > 0$ is small enough, $X^{-r}_{c + \eps}$ and $X_{c +\eps}$ are homotopy equivalent, and so are $X^{-r}_{c- \eps}$ and $X_{c-\eps}$.
This is summarized in \Cref{fig:diagramme_commutatif}.

\begin{figure}[h!]
    \centering
    \includegraphics[width = 1.1\textwidth]{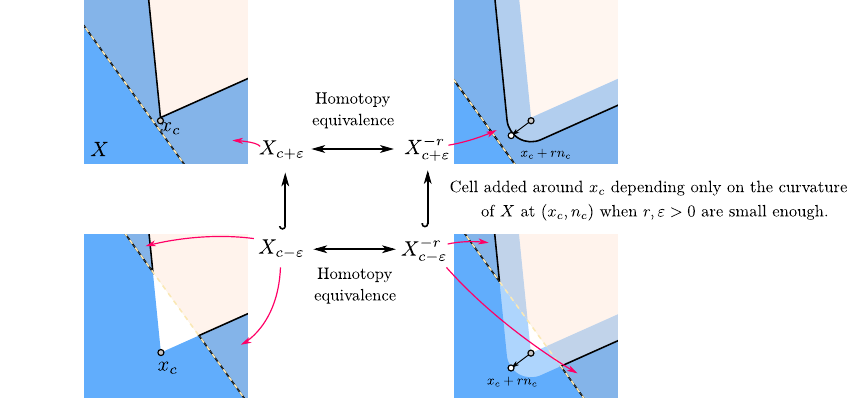}
    \caption{Commutative diagram in the proof of \Cref{th:cellule}.}
    \label{fig:diagramme_commutatif}
\end{figure}
\end{proof}

Finally, we prove that the previous result holds when there might be several critical points of $f_{|X}$ sharing the same critical value. 

\begin{theorem}[Morse theory for complementary regular sets]{}
Let ${X \subset \R^d}$ be a complementary regular set,
$f : U \to \R$ be a smooth function defined on an open set $U$ containing $X$, and let $c$ be a critical value of $f_{|X}$ so that every critical point of $f_{|X}$ in $f^{-1}(c)$ is non-degenerate. \\

Then this level set contains a finite number $\{ x_1, \dots, x_{p_c} \}$ of critical points. Furthermore denoting their indices (defined in \Cref{th:cellule}) by $\lambda_{x_1}, \dots,  \lambda_{x_{p_c}}$, for all $\eps > 0$ small enough $X_{c+\eps}$ has the homotopy type of $X_{c-\eps}$ with exactly $p_c$ cells of respective dimension $\lambda_{x_1}, \dots \lambda_{x_{p_c}}$ attached.
\end{theorem}

\begin{proof}
 We have already obtained in \Cref{lem:isol3} that when $x \in \partial X$ is the unique critical point in $f^{-1}(c)$ and is non-degenerate, $x^r = x + r \normalise{\nabla f(x)} \in X^{-r}$ must be the only critical point with associated critical value close to $c$ of the map $f$ translated by $- r \normalise{\nabla f(x)}$. When there are several critical points of this type in $f^{-1}(c)$, define another map $f_r$ by locally translating $f$ by $- r \normalise{\nabla{ f(x_i)}}$ around each critical point $x_i \in f^{-1}(c) \cap \partial X$. 
Interpolating these translations we can ensure that \Cref{cor:homotopy_equivalence} is satisfied, that is, $X_{c +\eps}$ and $X^{-r}_{c+ \eps}$ share the same homotopy type when $\eps, r$ are small enough. We give below a rigorous construction following this idea, illustrated in \Cref{fig:multi_attach}.

\begin{figure}[h!]
    \centering
    \includegraphics[width = 0.9\textwidth]{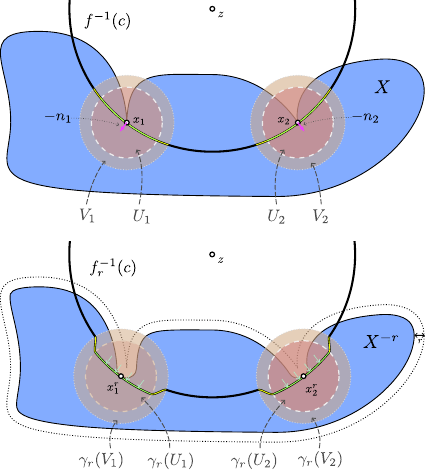}
     
    \caption{ \textbf{Construction of $f_r$ and $\gamma_r$ at a critical value with multiple critical points.}\\
    $X \subset \R^2$ a complementary regular set, $f = x \mapsto \norme{x - z}^2$. \\
    \textbf{Top:} $x_1, x_2$ are two critical points of $f_{|X}$ of value $c = f(x_1) = f(x_2)$ and $n_i = - \normalise{\nabla f(x_i)}$.\\
    \textbf{Bottom:} The level set $f_r^{-1}(c)$ is the circle $f^{-1}(c)$ deformed by local translations $x \mapsto x + r n_i$ in the neighborhood of each $x_i$.}
    \label{fig:multi_attach}
\end{figure}
Let $c$ be a critical value of $f_{|X}$. Since non-degenerate critical points are isolated (\Cref{prop:isol2}) and $X$ is compact, the level set $f^{-1}(c)$ must contain a finite number of points critical for $f_{|X}$.  Let $x_1, \dots, x_k$ be all the critical points inside $f^{-1}(c) \cap \partial X$ and $x_{k+1}, \dots x_p$ be the critical points in $\interior(X)$. Put $n_i \coloneqq -\frac{\nabla f(x_i)}{\norme{ \nabla 
 f(x_i)}}$ if $i \leq k$ and $n_i =0$ otherwise, and let $x^r_i = x_i - r n_i$. Let $n(x)$ be the function mapping $x$ to the $n_i$ associated to the closest critical point $x_i$ of $x$.
Let $U_i \subset V_i$ be respectively closed and open balls containing $x_i$ such that $\overline{V_i} \cap \overline{V_j} = \emptyset$ when $j \neq i$. Let $\eta_c$ be a smooth function on $\R^d$ with values in $[0,1]$ such that $\eta_c$ is constant of value $1$ inside each $U_i$ and $0$ outside $\bigcup V_i$. 
The map $n_c: y \mapsto \eta_c(y) n(y)$ is well-defined and smooth provided the $U_i$ are small enough. When $r$ is small enough, the map $\gamma_r : y \mapsto y + r n_c(y)$ has the following behavior:
\begin{itemize}
\item[(1)] in $U_i$, $\gamma_r$ is the translation by $r n_i$;
\item[(2)] in $X \setminus (\bigcup V_i)$, $\gamma_r$ is the identity;
\item[(3)] in each $V_i$, $\gamma_r$ is a smooth transition between the translation behavior in $U_i$ and the identity.
\end{itemize}

Now define $f_r$ to be $f$ locally translated around the critical points:
\[ f_r \coloneqq f \circ \gamma_r : y \mapsto f(y + rn_c(y)). \]

Since $n_c$ is smooth with $\norme{n_c} \leq 1$, \Cref{cor:homotopy_equivalence} applies with $X^{-r}_t = X^{-r} \cap f_r^{-1}(-\infty, t]$, so that for every regular value $t$ of $f_{|X}$, for every $r > 0$ small enough $t$ is a regular value of $f_{r|X^{-r}}$ and $X^{-r}_t$ has the same homotopy type as $X_t$.\\
To obtain our desired result, it thus suffices to show the property (P): for $\eps, r > 0$ small enough, the only critical points of 
$f_{r|X^{-r}}$ with critical values in $[c - \eps, c+\eps]$ are the $x^r_i$ and that these are non-degenerate with index $\lambda_i$.
 If this property holds, in a fashion similar to the previous theorem, for such $r, \eps$ by $C^{1,1}$ Morse theory $X^{-r}_{c + \eps}$ has the homotopy type of $X^{-r}_{c - \eps}$ with cells of dimension $\lambda_{i}$ attached  around each $x^r_{i}$, and our theorem follows from \Cref{cor:homotopy_equivalence}.

From \Cref{lem:isol_index}  we already know
 that the $(x^r_i)_{1 \leq i \leq p}$ are non-degenerate critical points of $X^{-r}$ for $f_{r {|X^{-r}}}$ with corresponding index $(\lambda_{x_i})_{1 \leq i \leq p}$.
Moreover, $\gamma_{r}$ restricted to each $U_i$ is simply the translation by $rn_i$ and $x_i$ is the unique critical point of value $c$ within $U_i$. 
 This is exactly the setting of \Cref{lem:isol3}, so that property $(P)$ holds for sequences in $\gamma_{r_i}(U)$ where $U \coloneqq \cup_i U_i$.
Assuming $(P)$ does not hold in the compact set $X \setminus U$, there exist positive sequences $r_i, \delta_i \to 0$, a sequence $y_i \in X \setminus U$ such that 

\vspace{0.1cm}
\begin{center}
\begin{tabular}{l r}
\quad \textbullet\ $c-\delta_i
\leq f_{r_i}(y_i)
\leq c+\delta_i$
&
\textbullet\ $y_i$ is critical for $f_{r|X^{-r_i}}$.

\end{tabular}
\end{center}
\vspace{0.1cm}
As in \Cref{lem:isol3}, on $(X \setminus U) \cap f^{-1}(c)$, $\norme{\nabla f}$ is non-zero, so that by continuity we can assume there exists $K, \eps > 0$ such that $\norme{\nabla f(y_i)} \geq K$ for all $i$. Since $y_i$ is critical for $(f_r)_{|X^{-r}}$, $y_i$ must belong to $\partial(X^{-r})$.
Any cluster point $\overline{x}$ of the sequence $y_i$ satisfies $\nabla f(\overline{x}) = \lim_{i \to \infty} \nabla f_{r_i}(y_i)$, and by the same reasoning as the paragraph before \Cref{eq:crit_lim}, the fact that 
\[ -\nabla f_{r_i}(y_i) \in \Nor(X^{-r},y_i) \]
implies that 
\[ -\nabla f(\overline{x}) \in \Nor(X,\overline{x}), \]
 so that $\overline{x}$ must be a critical point of $f_{|X}$ in $\overline{X \setminus U}$ of value $c$, which is a contradiction. 
\vspace{1em}

\end{proof}

\begin{corollary}[Theorem 1.1]
\label{cor:MR}
Let ${X \subset \R^d}$ be a complementary regular set, $f :U \to \R$ be a smooth 
defined on an open set $U$ containing $X$, and assume that $f_{|X}$ is Morse.
Then $f_{|X}$ has a finite number of critical points, and
\begin{itemize}
    \item If $[a,b]$ does not contain any critical value of $f_{|X}$, $X_a$ is a deformation retract of $X_b$.
    \item 
     If $c$ is a critical value, for all $\eps > 0$ small enough $X_{c+\eps}$ has the homotopy type of $X_{c-\eps}$ with exactly $p_c$ cells of respective dimension $\lambda_{x_1}, \dots \lambda_{x_{p_c}}$ attached, where $\{ x_1, \dots, x_{p_c} \}$ is the set of critical points of $f_{|X}$ in $f^{-1}(c)$,  where 
\begin{align*}
 \lambda_x \coloneqq & \text{ index } \text{ of the Hessian of }  f_{|X} \text{ at } x \\ + &
\text{ number of infinite curvatures of } \Nor(X) \text{ at } \left (x,-\normalise{\nabla f(x)} \right )
\end{align*}
when $x \in \partial X$, and $
 \lambda_x \coloneqq \text{ index } \text{ of the Hessian of }  f \text{ at } x$
 when $x \in \interior{(X)}$.
     \item If $c$ is a regular value of $f_{|X}$, $X_c$ has the homotopy type of a CW-complex obtained by the attachment, in ascending order of $f(x)$, of a cell of dimension $\lambda_{x}$ for each critical point $x$ with $f(x) \leq c$.
\end{itemize}
\end{corollary}
\begin{proof}
The fact that $f_{|X}$ admits a finite number of critical points follows from the fact that $X$ is compact and that non-degenerate critical points are isolated (\Cref{prop:isol2}).
The constant homotopy type lemma (first point in the list) is \Cref{th:isotopy}. The second point is the previous \Cref{th:cellule}, and the final point follows easily by induction on the ranked critical values of $f_{|X}$ using the two previous points.
\end{proof}

\begin{remark}[Morse functions are plentiful]
\label{rem:plentiful}
When $X$ is a complementary regular set, it is straightforward from previous computations and the fact that $\Nor(X,x) = - \Nor( \complementaire{X},x)$ that $f_{|X}$ is Morse if and only if $(-f)_{|\complementaire{X}}$ is. We thus know from Fu \cite{fuCurvatureMeasuresGeneralized1989} that almost every linear form restricted to $X$ is Morse. Similar computations 
from geometric measure theory also show that as a function of $x$, maps $d^2_x : z \mapsto \norme{z -x}^2$ restricted to $X$ are Morse almost everywhere in $\R^d$ with respect to the Lebesgue measure
- see \cite{PersistentIntrinsicVolumes} for the explicit computation.
\end{remark}

\subsection{Connections to persistent homology}
\label{sub:persistence}
The \textit{persistent homology} of a map $f$ taking real values studies the topological evolution of the sublevel set filtration of $f$ by considering the families $H_i(f^{-1}(-\infty, t])_{t \in \R}$ for each degree $i$. It does so by determining the two filtration values at which a topological feature \textit{is born} and at which it \textit{dies}. Pairing the birth and death of topological features yields the so-called \textit{persistence diagram} in degree $i$ of $f$, which is a multi-set of $\{ (x,y),  x \leq y \} \subset \R^2$. Their disjoint union (over $i$) forms the persistence diagram of $f$, denoted by $\dgm(f)$.
Persistence theory provides several tools to compare persistence diagrams, with various applications in data analysis \cite{FredBertrand,Carlsson} where it is used to understand both the geometry and topology of high-dimensional datasets, but also in more theoretical domains of mathematics such as symplectic geometry \cite{Polterovitch} and geometric measure theory \cite{PersistentIntrinsicVolumes}. 

By giving a sense of Morse functions over a complementary regular set, one can effectively apply the persistence framework to them. 
 Indeed, let $X \subset \R^d$ be complementary regular and let $f$ be a smooth map such that $f_{|X}$ is Morse. Then for each degree $0 \leq i \leq d$, we determine birth and deaths by the following procedure:
\begin{itemize}
    \item There are finitely many topological events in the family of vector spaces $(H_i(X_t))_{t \in \R}$, where the homology is taken over any field - this stems from the \textit{constant homotopy type} property.
    \item For any critical point $x$ of $f$ with index $d_x$, either $\dim H_{d_x}(X_t)$ increases by one at $t = f(x)$, corresponding to the \textit{birth} of a $(d_x)$-dimensional feature (i.e., the appearance of a new cycle in the filtration), or $\dim H_{d_x -1}(X_t)$ decreases by one, corresponding to the \textit{death} of a $(d_x -1)$-dimensional feature (i.e, a cycle that gets filled in the filtration).
 \end{itemize} 
While the way births and deaths are paired is given by persistence theory, $f_{|X}$ being Morse guarantees that the persistence diagram $\dgm(f_{|X})$ is well-defined.

An important tool of persistence theory is the \textit{Persistent Homology Transform} (PHT) \cite{PHT}, which to a set $Y \subset \R^d$ associates the family of persistent diagrams $\dgm(h_{\nu| Y})$ of the height functions $h_{\nu} : x \mapsto \scal{x}{\nu}$ where $\nu$ ranges among $\Sphere^{d-1}$. Since the maps 
$h_{\nu|X}$ are Morse for almost every $\nu$ by \Cref{rem:plentiful}, we have the following:

\begin{proposition}[Persistent homology transform over complementary regular sets]
Let $X \in \R^d$ be a complementary regular set. Then the diagram $\dgm(h_{\nu |X})$ is well-defined for almost every $\nu \in \Sphere^{d-1}$.
\end{proposition}

Another important construction is the \textit{Euler Characteristic Transform} (ECT) \cite{ECT}, which to $X$ associates the function $ECT_X : (\nu, t) \mapsto \chi(X \cap h_{\nu}^{-1}(-\infty, t])$. When $X \subset \R^d$ is a complementary regular set, we have already seen in \Cref{rem:plentiful} that for almost all $\nu \in \Sphere^{d-1}$, $h_{\nu |X}$ is Morse. It follows that for every regular value $t$ of this map, the set $X \cap h_{\nu}^{-1}(-\infty ,t]$ has the homotopy type of a CW-Complex, so its Euler characteristic $\chi(X \cap h_{\nu}^{-1}(-\infty ,t])$ is well-defined.

\begin{proposition}[Euler Characteristic transform of complementary regular sets]
Let $X \subset \R^d$ be a complementary regular set. Then the map
\[ ECT_X : (\nu, t) \mapsto \chi( X \cap h_{\nu|X}^{-1}(-\infty, t]) \]
is well-defined for almost every $(\nu, t) \in \Sphere^{d-1} \times \R$.
\end{proposition} 

For compact sets definable in an o-minimal structure, the injectivity of the ECT follows from an explicit inverse construction, first given in \cite{SchapiraInverse}. However, complementary regular sets need not be definable in an o-minimal structure, so this argument does not apply here.
Nevertheless, the ECT remains injective in our context, as we now show using  the uniqueness theorem for normal cycles of Fu (Theorem 3.2 of \cite{fuCurvatureMeasuresSubanalytic1994}). 

\begin{proposition}[Injectivity of the ECT for complementary regular sets]
Let $X, Y \subset \R^d$ be two complementary regular sets such that $ECT_X = ECT_Y$ almost everywhere, i.e., such that for almost all $(\nu,t) \in \Sphere^{d-1} \times \R$, 
\[\chi( X \cap h^{-1}_{\nu}(- \infty, t]) = \chi( Y \cap h^{-1}_{\nu}(-\infty, t]). \] Then $X = Y$.
\end{proposition}
\begin{proof}
For a complementary regular set $Z$, consider the current $N_Z = \rho_{\#}(N_{\complementaire{Z}})$ where $N_{\complementaire{Z}}$ is the current given by integration over the normal bundle $\Nor(\complementaire{Z})$, and  $\rho_{\#}$ indicates the push-forward of a current by the antipodal map $\rho : (x,n) \mapsto (x, -n)$. From 9.10 in Section 9.2 of \cite{Ja-2019}, $N_Z$ is the normal cycle of $Z$, that is, a $(d-1)$-current of $\R^d \times \Sphere^{d-1}$ used to describe second order information of $Z$. The uniqueness theorem of Fu states, using constructions from geometric measure theory, that the Euler characteristic of a set intersected with almost all half-space determines its normal cycle in the sense that two sets whose ECT coincide almost everywhere must have the same normal cycles. Thanks to the explicit construction of $N_Z$ given above for $Z=X$ and $Z=Y$, $N_X = N_Y$ implies $\Nor(\complementaire{X}) = \Nor(\complementaire{Y})$, so that $\partial X = \partial Y$ by taking the projection onto the first factor of $\R^d \times \Sphere^{d-1}$. Since $X, Y$ are compact Lipschitz domains, their interior must also coincide and we have $X = Y$.
\end{proof}

\bibliographystyle{ieeetr}
\bibliography{ref}

\end{document}